\documentclass[11pt,a4paper,reqno]{amsart}
\usepackage{amssymb,amsmath}
\usepackage{xcolor}
\usepackage{bbm}
\usepackage{dsfont}
\usepackage{bigints}
\usepackage{mathtools}
\usepackage{mathrsfs}
\usepackage{hyperref}
\usepackage{enumitem}
\usepackage[autostyle]{csquotes}
\usepackage[mathscr]{euscript}
\usepackage{fouridx}
\usepackage{mathtools}
\usepackage{cases}

\usepackage{todonotes}
\renewcommand{\mathcal}[1]{\mathscr{#1}}

\usepackage{mathtools}

\newcommand{\1}{\mathds{1}}% requires package dsfont

\newcommand\dela[1]{}

\newcommand{\cO}{\mathscr{O}}

\newcommand{\topp}{\mathrm{top}}
\newcommand{\strong}{\mathrm{strong}}

\newcommand{\esssup}{\mathrm{esssup}}
\newcommand{\toup}{\nearrow}

\newcommand{\Za}{Z}

\newcommand{\Xa}{X}
\newcommand{\Ya}{Y}
\newcommand\proj{\Pi}
\newcommand\nA{A}
\newcommand\Leb{\mathrm{Leb}}
\newcommand{\LebAuxProbSpace}{L^2(\Omega, \mathscr{F},\mathbb{P};E)}
\newcommand{\LpAuxProbSpace}{L^p(\Omega, \mathscr{F},\mathbb{P};E)}
\newcommand{\AuxProbSpace}{(\Omega, \mathscr{F},\mathbb{P})}
\newcommand{\tP}{\mathbb{P}}

\newcommand{\tE}{\mathbb{E}}
\newcommand{\tF}{\mathbb{F}}

\newcommand{\embed}{\hookrightarrow }

\numberwithin{equation}{section}
\newcommand\new[1]{}

\allowdisplaybreaks

\advance\textwidth by 2.9cm

\advance\oddsidemargin by -1.6cm
\advance\evensidemargin by -1.6cm

\newtheorem{theorem}{Theorem}[section]
\newtheorem{definition}[theorem]{Definition}
\newtheorem{proposition}[theorem]{Proposition}
\newtheorem{lemma}[theorem]{Lemma}
\newtheorem{corollary}[theorem]{Corollary}
\newtheorem{remark}[theorem]{Remark}

\newcommand{\rV}{H_A}%{\mathrm{V}}
\newcommand{\rE}{\mathrm{E}}

\newcommand\restr[2]{{
		\left.\kern-\nulldelimiterspace
		#1
		\vphantom{\big|}
		\right|_{#2}
}}

\newcommand{\threepartdef}[6]
{
	\left\{
	\begin{array}{lll}
		#1 & \mbox{if } #2 \\
		#3 & \mbox{if } #4 \\
		#5 & \mbox{if } #6
	\end{array}
	\right.
}

\newcommand{\lb}{\langle}
\newcommand{\rb}{\rangle}

\newcommand{\Dom}{\mathscr{D}}

\makeatletter
\def\@tocline#1#2#3#4#5#6#7{\relax
	\ifnum #1>\c@tocdepth % then omit
	\else
	\par \addpenalty\@secpenalty\addvspace{#2}%
	\begingroup \hyphenpenalty\@M
	\@ifempty{#4}{%
		\@tempdima\csname r@tocindent\number#1\endcsname\relax
	}{%
		\@tempdima#4\relax
	}%
	\parindent\z@ \leftskip#3\relax \advance\leftskip\@tempdima\relax
	\rightskip\@pnumwidth plus4em \parfillskip-\@pnumwidth
	#5\leavevmode\hskip-\@tempdima
	\ifcase #1
	\or\or \hskip 1em \or \hskip 2em \else \hskip 3em \fi%
	#6\nobreak\relax
	\dotfill\hbox to\@pnumwidth{\@tocpagenum{#7}}\par
	\nobreak
	\endgroup
	\fi}
\makeatother

\begin{document}

\bibliographystyle{plain}
\pagenumbering{arabic}

\title[Energy critical 2-D stochastic wave equation]{Local solution to an energy critical 2-D stochastic wave equation with
	exponential nonlinearity in a bounded domain}
\author{Zdzis{\l}aw Brze\'zniak and Nimit Rana}
\address{Department of Mathematics \\
The University of York \\
Heslington, York, YO105DD, UK} \email{zdzislaw.brzezniak@york.ac.uk}
\address{Fakult\" at f\"ur Mathematik \\
	Universit\"at Bielefeld\\
	Universit\"atsstra\ss e 25, 33615 Bielefeld, Germany} \email{nrana@math.uni-bielefeld.de}
%\today

\begin{abstract}
  We prove the existence and the uniqueness of a local maximal solution to an $H^1$-critical stochastic wave equation with multiplicative noise on a smooth bounded domain $ \mathscr{D} \subset \mathbb{R}^2$ with exponential nonlinearity. First, we derive the appropriate deterministic and stochastic Strichartz inequalities in suitable spaces and, then use them in arguments based on fixed point method to show the local well-posedness result. We also present an explosion result for the constructed unique local maximal solution.
\end{abstract}

\date{}
\maketitle

\tableofcontents

\section{Introduction}\label{sec-Intro}

In this paper we consider a nonlinear wave equation subject to random forcing, called the stochastic nonlinear  wave equation (SNLWE). Due to its numerous applications to physics, relativistic quantum mechanics and oceanography, SNLWEs have been thoroughly studied under various sets of assumptions, see for example \cite{Brz+Millet_2014}-\cite{Brz+O_2007}, \cite{C1972}-\cite{DS2009B}, \cite{KZ2001}-\cite{O2004}, \cite{O2006}-\cite{O2010II}, \cite{PZ2000}-\cite{P2002} and references therein. The case that has so far attracted the most attention seems to be of the stochastic wave equation with initial data belonging to the energy space $H^1(\mathbb{R}^d) \times L^2(\mathbb{R}^d)$. For such equations, the nonlinearities can be of polynomial type, for instance the following SNLWE
 \begin{equation}\label{O2010IIWE}
	u_{tt} - \Delta u = -u |u|^{p-1} + |u|^q \dot{W}, ~\textrm{ s.t. }~ u(0) = u_0, ~ \partial_t u(0) =u_1,
\end{equation}
with the suitable exponents $p,q \in (0,\infty)$; see a series of papers by Ondrej\'at \cite{O2004}, \cite{O2006}-\cite{O2010II}. Another extensively studied important case is when the initial data $(u_0,u_1)$ is in $L^2(\mathbb{R}^d) \times H^{-1}(\mathbb{R}^d)$ (possibly with weights), see \cite{PZ2000}, \cite{P2002} for more details. Similar problems on a  bounded domain have been investigated  in \cite{Brz+O+S_2016}, \cite{C2006} and \cite{O2004}.

In the deterministic case, see for instance \cite{S1992},  the question of solvability of \eqref{O2010IIWE} without noise, when the initial data belongs to  $H^1(\mathbb{R}^d) \times L^2(\mathbb{R}^d)$, has been investigated in the following three cases:
(i) subcritical, i.e. $p<p_c$; (ii) critical, i.e. $p=p_c$; and (iii) supercritical, i.e. $p > p_c$, where $p_c = \frac{d+2}{d-2}$. In particular, for $d=2$, any polynomial nonlinearity is subcritical. Therefore,  an exponential nonlinearity is a legitimate choice for a critical one. Nonlinearities of exponential type have been studied in many physical models, e.g. a model of self-trapped beams in plasma, see \cite{LLT1977}, and mathematically in   \cite{Baraket_2004}, \cite{C1979}, \cite{IMM2006}-\cite{IJ2011} and \cite{NO1999}. With the help of suitable Strichartz estimates, the existence of global solutions have been proved, in \cite{IMM2006}-\cite{IJ2011},  in the cases when the initial energy is strictly below or at the threshold given by the sharp Moser-Trudinger inequality. Moreover, an instability result has been shown when the initial energy is strictly above the threshold.

Our aim here is to extend the existing studies to the wave equation with exponential nonlinearity subject to randomness. In this way, we generalise the above mentioned results of Ondrej\'at for two dimensional domains, by allowing the exponential nonlinearites, as well as the results of Ibrahim, Majdoub, and Masmoudi and others by allowing randomness. To be precise, we are  interested in the following stochastic nonlinear wave equation on a smooth bounded  domain ${\cO} \subset \mathbb{R}^2$,
\begin{equation}\label{AbsWe}\left\{
\begin{array}{ll}
u_{tt} +\nA u + F(u) = G(u) \dot{W} \quad \textrm{ in } [0,\infty) \times {\cO} \\
u(0) = u_0, ~~  u_t(0) = u_1 \quad \textrm{ on } {\cO},
\end{array}
\right.
\end{equation}
where  $\nA$ is either $-\Delta_{D}$ or $-\Delta_{N}$, i.e. $-\nA$ is the Laplace-Beltrami operator with the Dirichlet or the Neumann boundary conditions, respectively; $(u_0,u_1) \in \Dom(\nA^{\frac12 }) \times L^2({\cO})$; $W = \{W(t): t \geq 0  \}$  is a $K$-cylindrical Wiener process for  a   real separable Hilbert space $K$; $F$ and $G$ are locally Lipschitz maps with some growth properties. In particular, the functions $F(u)$ and $G(u)$ are allowed to be of the form $\pm u \bigl(e^{4 \pi u^2} -1 \bigr)$ and hence our results  cover the recent results obtained in  \cite{IJ2011}.
Detailed and precise assumptions on the model are stated in the subsequent sections. We would like to stress that, to the best of our knowledge, the present paper is the first one to study the wave equations in two dimensional domain with an exponential nonlinearity and an additive or multiplicative noise. As compared with the deterministic paper \cite{IJ2011} we prove a counterpart of its Theorem 1.4, i.e. we prove the existence of a unique local solution with a smallness condition on the gradient of initial position. The difficulty we encounter in the present paper is that although the nonlinear function is defined on the whole space $\Dom(A^{\frac{1}{2}})$, it's values belong to the space $L^2$ only for the elements of the space $\Dom(A^{\frac{1}{2}})$ whose norm is sufficiently small.
We emphasize that our proof of stochastic Strichartz estimate, see Theorem \ref{thm-StocStrichartzEstimate}, simplifies and clarifies the one from \cite{Brz+Millet_2014}.
Since the proof of the existence and the uniqueness presented here is obtained by means of appropriate Strichartz estimates and these estimate are  different for the  full domain case, we will address the question of solvability of \eqref{AbsWe} on $\mathbb{R}^2$ in a forthcoming paper.

The organization of the present paper is as follows. In Section \ref{sec-Notation}, we introduce our notation and provide the required definitions used in the paper.  In Sections \ref{sec-inHomStrichartz} and \ref{sec-Strichartz stochastic}, we derive the required inhomogeneous and stochastic Strichartz estimates, respectively, by the methods introduced in \cite{Burq+L+P_2008}-\cite{Burq+P_2009} and \cite{Brz+Millet_2014}. Section \ref{sec-localExistence} is devoted to the estimates which are sufficient to apply the Banach Fixed Point Theorem in a suitable space. We study the approximated version of problem  \eqref{AbsWe} in Section \ref{sec-approximation problem}.  The proof of the existence and uniqueness of a local maximal solution is given in Section \ref{sec-main}. Here we also formulate the results about the explosion for the constructed unique local maximal solution.
In Appendix \ref{sec-stopped processes}, we provide a rigorous  justification of our definition of a local mild solution. In Appendix \ref{sec-pointwise evaluation}, we formulate a result about pointwise evaluation of $L^p$-valued Bochner integrals. In Appendix \ref{sec-def sol} we state, without proof, an equivalence of two natural definitions of a mild solution for stochastic PDE \eqref{AbsWe}. We conclude the paper with Appendix \ref{sec-Gronwall}  and \ref{sec-brze-carr:app-A} in which we prove, respectively, a slight simplification of the stochastic Gronwall Lemma \cite[Lemma~5.3]{GHZ2009} and a generalization of an existence of a Lipschitz extension result \cite[Corollary 3]{Brz+C_2003}.

\section{Notation and conventions}\label{sec-Notation}
In this section we introduce  notation and some basic estimates that we use throughout the paper. We write $a \lesssim b$ if there exists a universal constant $c >0$, independent of $a,b$, such that $a \leq cb$, and we write $a \simeq b$ when $a \lesssim b$ and $b \lesssim a$.  In case we want to emphasize the dependence of $c$ on some parameters $a_1,\ldots,a_k$, then we write, respectively, $\lesssim_{a_1,\ldots,a_k}$ and $\simeq_{a_1, \ldots, a_k}$. For any two Banach spaces $X, Y$, we denote by $\mathcal{L}(X, Y)$ the space of linear bounded operators $L : X \to Y$.

To state the definitions of required spaces here, we denote by $E$ and $H$ a separable Banach and Hilbert space, respectively.

\subsection{Function spaces and interpolation theory}
In the next few basic definitions and remarks, which are included here for the reader's convenience, from function spaces and interpolation theory we  borrow the notation from \cite{Triebel_1978B}.

By $L^q({\cO})$, for $q \in [1,\infty)$ and a bounded smooth domain ${\cO}$ of $\mathbb{R}^2$, we denote the classical real Banach space of all (equivalence classes of) $\mathbb{R}$-valued $q$-integrable Lebesgue measurable functions on ${\cO}$. The norm in $L^q({\cO})$ is given by
\begin{equation*}
	\Vert u \Vert_{L^q({\cO})} := \left( \int_{{\cO}} |u(x)|^q  \, dx  \right)^{\frac{1}{q}}, \quad u \in L^q({\cO}).
\end{equation*}

By $L^\infty({\cO})$ we denote the real Banach space of all (equivalence classes of) Lebesgue measurable essentially bounded $\mathbb{R}$-valued functions defined on ${\cO}$ with the norm \begin{equation*}
	\Vert u \Vert_{L^\infty({\cO})} := \esssup~ \{ | u(x)| : x \in {\cO} \}, \quad u \in L^\infty({\cO}).
\end{equation*}

For any $T > 0$, the space $\mathrm{C}\bigl([0,T]; H \bigr)$  of all $H$-valued continuous functions $u : [0,T] \to H$ endowed with the norm
\begin{equation*}
	\Vert u \Vert_{\mathrm{C}([0,T]; H)} := \sup_{t \in [0,T]} \Vert u(t)\Vert_{H}, \quad u \in \mathrm{C}\bigl([0,T]; H \bigr).
\end{equation*}
is the real Banach space.

For every $p \in [1,\infty)$, we  define the space   $L^p\bigl(0,T; E \bigr)$ as the vector space  of all (equivalence classes of) $E$-valued strongly measurable functions
$u : [0,T] \to E$  such that $\int_{0}^{T} \Vert u(t)\Vert_E^p  \, dt < \infty$. The space $L^p\bigl(0,T; E \bigr)$  endowed with the norm
\begin{equation*}
	\Vert u \Vert_{L^p(0,T;E)}  := \left( \int_{0}^{T} \Vert u(t)\Vert_E^p  \, dt  \right)^{\frac{1}{p}}, \quad u \in L^p(0,T;E),
\end{equation*}
is a real Banach space.

For $\alpha \in(0,1)$, by $\dot{\mathrm{C}}^\alpha({\cO})$, the homogeneous $\alpha$-H\"older space, we mean the set of continuous
functions $u$ whose H\"older  coefficient
%
$$\Vert u \Vert_{\dot{\mathrm{C}}^\alpha({\cO})} := \sup_{x \neq y} \frac{|u(x) - u(y)|}{|x-y|^\alpha},$$
is finite. Note that the H\"older coefficient serves as a seminorm. The inhomogeneous $\alpha$-H\"older space is
$\mathrm{C}^\alpha({\cO}) = \dot{\mathrm{C}}^\alpha({\cO}) \cap L^\infty({\cO})$  endowed with the norm $\Vert u \Vert_{\dot{\mathrm{C}}^\alpha({\cO})} + \Vert u \Vert_{L^\infty({\cO})} $.

For any $s \in \mathbb{R}$ and $q \in (1,\infty)$, the Sobolev space $H^{s,q}({\cO})$ is defined as the restriction of $H^{s,q}(\mathbb{R}^2)$ (see e.g. \cite[Definition 2.3.1/1]{Triebel_1978B}) to ${\cO}$ with norm
\begin{equation}\label{eqn-H^s,q(D) space}
	\| f \|_{H^{s,q}({\cO})} := \inf_{\substack{g \upharpoonright {\cO} =f \\ g \in H^{s,q}(\mathbb{R}^2)}} \| g \|_{H^{s,q}(\mathbb{R}^2)}, \quad f \in H^{s,q}({\cO}).
\end{equation} Here $g \upharpoonright {\cO}$ is the restriction in the sense of distribution.  We denote the closure of $\mathrm{C}_0^\infty({\cO})$, the set of smooth functions defined over ${\cO}$ with compact support, in  $H^{s,q}({\cO})$ by  ${H}_0^{s,q}({\cO})$.

Throughout the whole paper, we denote by $\nA$ the Dirichlet or the Neumann$-$Laplacian on Hilbert space $L^2({\cO})$ with domain, respectively,
\begin{align}
& \Dom(-\Delta_{D}) = H^{2,2}({\cO}) \cap {H}_0^{1,2}({\cO}), \nonumber\\
\textrm{ and } \quad & \Dom(-\Delta_{N}) =  \left\{f \in H^{2,2}({\cO}) : \partial_\nu f \upharpoonright \partial {\cO} =0  \right\}. \nonumber
\end{align}
Here $\nu$ denotes the outward normal unit vector to $\partial {\cO}$.
It is well-known, see e.g. \cite{RT1997B}, that the Dirichlet Laplacian $(-\Delta_{D},\Dom(-\Delta_{D}))$ is a positive self-adjoint operator on $L^2({\cO})$ and there exists an orthonormal basis $\{e_j \}_{j \in \mathbb{N}}$ of $L^2({\cO})$ which consists of eigenvectors of $-\Delta_{D}$. If we denote the corresponding eigenvalues by $\{ \lambda_j^2 \}_{j \in \mathbb{N}}$, then we have
\begin{align}
-\Delta_{D} e_j = \lambda_j^2 e_j; \quad e_j \in \Dom(-\Delta_{D}), \forall j \geq 1 \textrm{; } \quad 0 < \lambda_1^2 \leq \lambda_2^2  \leq \cdots \textrm{ and } \lambda_n^2 \xrightarrow[n \to \infty]{} \infty. \nonumber
\end{align}
In the case of the Neumann Laplacian, $(-\Delta_{N},\Dom(-\Delta_{N}))$ is a non-negative self-adjoint operator on $L^2({\cO})$ and there exists an orthonormal basis $\{e_j \}_{j \in \mathbb{N}}$ of $L^2({\cO})$ which consists of eigenvectors of $-\Delta_{N}$. Moreover, if we denote the corresponding eigenvalues by $\{ \lambda_j^2 \}_{j \in \mathbb{N}}$, then we have
\begin{align}
-\Delta_{N} e_j & = \lambda_j^2 e_j; \quad e_j \in \Dom(-\Delta_{N}), \forall j \geq 1; \quad \lambda_n^2 \xrightarrow[n \to \infty]{} \infty, \nonumber\\
& \textrm{ and } \quad 0 = \lambda_1^2 = \lambda_2^2 = \cdots=\lambda_{m_0}^2 < \lambda_{m_0 +1}^2  \leq \lambda_{m_0+2}^2 \leq \cdots, \nonumber
\end{align}  for some $m_0 \in \mathbb{N}$.
Since we work with both the operators simultaneously, we denote the pair of operator and its domain by $(\nA,\Dom(\nA))$ and make the distinction where required.

From the functional calculus of self-adjoint operators, see for instance \cite{EZ1995B}, it is known that, the power $\nA^s$ of operator $\nA$, for every $s \in \mathbb{R}$, is well-defined and self-adjoint. It is also known that, for any $s \in \mathbb{R}$, $\Dom(\nA^{s/2})$, where $A=-\Delta_D$ or $A=-\Delta_N$, with the following norm
\begin{align}\nonumber
\| u \|_{\Dom(\nA^{s/2})} & := \bigl( \sum_{j \in \mathbb{N}} (1+\lambda_j^2)^{s} |\langle u, e_j \rangle_{L^2({\cO})}|^2  \bigr)^{\frac12 }, \nonumber
\end{align}
is a Hilbert space. For $s \in (0,2)$ the space $\Dom(\nA^{s/2})$ is equal to the following complex interpolating space, refer \cite[2.5.3/(13)]{Triebel_1978B},
\begin{equation*}
\Dom(\nA^{s/2}) = \left[ L^2({\cO}), \Dom(\nA)  \right]_{s/2}.
\end{equation*}

To derive the Strichartz estimate in suitable spaces, we also need to consider the Dirichlet or the Neumann$-$Laplacian on Banach space $L^q({\cO})$, $q \in (1,\infty)$,  denoted by $\nA_{D,q}$ and respectively $\nA_{N,q}$,  with domains, respectively,
\begin{align}
& \Dom(\nA_{D,q}) = H^{2,q}({\cO}) \cap {H}_0^{1,q}({\cO}) , \label{DirDom-q} \\
 \textrm{ and } \quad & \Dom(\nA_{N,q}) =  \left\{f \in H^{2,q}({\cO}) : \partial_\nu f \upharpoonright \partial {\cO} =0  \right\}.  \label{NeuDom-q}
\end{align}
Note that $\nA_{D,2}=-\Delta_{D}$ and $\nA_{N,2}=-\Delta_{N}$.

Under some reasonable assumptions on the regularity of the domain ${\cO}$, one can show that both of these operators have very nice analytic properties. In particular both have bounded imaginary powers with exponent strictly less than $\frac{\pi}{2}$ (and  thus both $-\nA_{D,q}$ and  $-\nA_{N,q}$ generate analytic semigroups on the space $L^q({\cO})$). As in \cite{Triebel_1978B}, one can define the fractional powers $\nA_{B,q}^{r/2}$, where as below $B=D$ or $B=N$. The domains  $\Dom(\nA_{B,q}^{r/2})$ of  these operators can be identified as certain subsets of the Sobolev spaces $H^{r,q}({\cO})$, see Lemma \ref{lem-EquiSpace1} below.

Next, we fix the notation for the required subspaces of $H^{s,q}({\cO})$ which are determined by differential operators.
Fix $k \in \mathbb{N}$ and let $B_j$, $j=1,\ldots,k$,  be differential operators on $\partial {\cO}$ defined by
\begin{equation*}
	B_jf(x) = \sum_{|\alpha| \leq m_j} b_{j,\alpha}(x) D^\alpha f(x), ~~ b_{j,\alpha} \in \mathrm{C}^\infty(\partial {\cO}).
\end{equation*}
 Then $\{B_j \}_{j=1}^k$ is said to be a normal system if and only if
\begin{equation*}
	0 \leq m_1 < m_2 < \cdots < m_k,
\end{equation*}
and for every vector $\nu_x$ which is normal to $\partial {\cO}$ at $x$ the following holds
\begin{equation*}
	\sum_{|\alpha| =m_j} b_{j,\alpha}(x) \nu_x^\alpha \neq 0, ~~ j=1,\ldots,k,
\end{equation*}
where for $\alpha \in \mathbb{N}^d$ and $y \in \mathbb{R}^d$, $y^\alpha = \Pi_{i} y_i^{\alpha_i}$.
\begin{definition}\label{def-HT1978B-Defn4.3.3/2}
	Let $\{B_j\}_{j=1}^k$ be a normal system as defined above for some $k \in \mathbb{N}$. For $s>0, q \in (1,\infty)$, we set
\begin{equation*}
	H_{\{B_j\}}^{s,q}({\cO}) : = \left\{ f \in H^{s,q}({\cO}) : B_jf \upharpoonright \partial {\cO}=0 \textrm{ whenever } m_j < s-\frac{1}{q} \right\}.
	\end{equation*}
\end{definition}
By taking the suitable choice of normal system $\{B_j\}$ in the Definition \ref{def-HT1978B-Defn4.3.3/2}, for $s>0$ and $q \in (1,\infty)$, we define  \begin{equation}
H_{D}^{s,q}({\cO}) := \left\{ f \in H^{s,q}({\cO}) : f \upharpoonright \partial {\cO}=0   \textrm{ if  } s >\frac{1}{q}  \right\}, \nonumber
\end{equation} and \begin{equation}
H_{N}^{s,q}({\cO}): = \left\{ f \in H^{s,q}({\cO}) : \nu_x \cdot \nabla f \upharpoonright \partial {\cO}=0 \textrm{ if } s> 1 + \frac{1}{q} \right\}. \nonumber
\end{equation}
Since the ${H}_0^{1,q}({\cO})$ space can also be defined by using the $f \upharpoonright \partial {\cO}=0$ condition appearing in \eqref{DirDom-q} and the Neumann boundary condition appearing in \eqref{NeuDom-q} can be written as $\nu_x \cdot \nabla f \upharpoonright \partial {\cO}=0$, we expect to have some relation between the spaces $H_{B}^{s,q}({\cO})$  and $\Dom(\nA_{B,q}^{s/2})$ where $A=-\Delta_B$ with $B=D$ or $B=N$. The next stated result, which is standard in the theory of interpolation spaces, see \cite[Theorem 4.3.3]{Triebel_1978B}, provides a suitable range of $s$ for which the function spaces $H_{B}^{s,q}({\cO})$ and $\Dom(\nA_{B,q}^{s/2})$ are equivalent.
\begin{lemma}\label{lem-EquiSpace1}
	With our notation from this section, we have the following \begin{enumerate}
		\item  For $s \in (0,2) \setminus \left\{1+\frac 1 q \right\}$,
\begin{equation*}
			H_{N}^{s,q}({\cO})  = \Dom(\nA_{N,q}^{s/2}).
		\end{equation*}
		\item  For $s \in (0,2) \setminus \left\{\frac 1 q \right\}$, \begin{equation*}
		H_{D}^{s,q}({\cO})  = \Dom(\nA_{D,q}^{s/2}).
		\end{equation*}
	\end{enumerate}
\end{lemma}
We close this subsection with the following well-known identity \begin{equation*}
\Dom(\sqrt{-\Delta_{D}}) = {H}_0^{1,2}({\cO}) \textrm{ and } \Dom(\sqrt{-\Delta_{N}}) = H^{1,2}({\cO}).
\end{equation*}

\subsection{Stochastic analysis}\label{subsec-Stochastic analysis}

Now we state a few required definitions from the theory of stochastic analysis, refer \cite{Brz_1997} and \cite{Brz+P_2001} for more details.  Throughout the whole paper we assume that  $(\Omega, \mathscr{F},\tP, \mathbb{F})$, where $\mathbb{F} := \{ \mathscr{F}_t: t \geq 0 \}$,  is  a filtered probability space which satisfies the \textbf{usual hypothesis}, that is, the filtration $\mathbb{F}$ is right continuous and the $\sigma$-field $\mathscr{F}_0$ contains all $\tP$-null sets of $\mathscr{F}$, see \cite[Definition I.1.1]{Metivier_1982}.

As the noise, we consider a cylindrical $\mathbb{F}$-Wiener process on a real separable Hilbert space $K$, see \cite[Definition 4.1]{Brz+P_2001}. Let us recall that $E$ is a separable Banach space. We denote by $L^p(\Omega, \mathscr{F}, \tP;E)$, for $p \in [1,\infty)$, the Banach space of all (equivalence classes of) $E$-valued random variables equipped with the norm \begin{equation*}
	\| X\|_{L^p(\Omega)} = \bigl(\tE \left[\| X\|_E^p \right]\bigr)^{\frac{1}{p}}, \quad  X \in L^p(\Omega, \mathscr{F}, \tP;E),
\end{equation*}
where $\tE$ is the expectation operator w.r.t. $\tP$.
\begin{definition}\label{def-gammaRadOper}
	For any $K$, a separable Hilbert space, the set $\gamma(K,E)$ of all $\gamma$-radonifying operators consists of all bounded operators $\Lambda: K \to E$ such that the series $\sum_{j=1}^{\infty}\beta_j \Lambda(f_j)$ converges in $\LebAuxProbSpace$  for some (or any) orthonormal basis $\{f_j\}_{j \in\mathbb{N}}$ of $K$ and some (or any) sequence $\{\beta_j\}_{j \in \mathbb{N}}$ of i.i.d. $N(0,1)$ real random variables on probability space $\AuxProbSpace$. We set
	\begin{equation} \label{eqn-norm gamma(K,E)}
	\Vert \Lambda \Vert_{\gamma(K,E)} := \left(   \tE \bigg\Vert \sum_{j \in \mathbb{N}} \beta_j \Lambda(f_j)\bigg\Vert_{E}^2 \right)^{\frac{1}{2}}.
	\end{equation}
\end{definition}
One may prove that $\Vert \cdot \Vert_{\gamma(K,E)}$ is a norm, and  $(\gamma(K,E), \Vert \cdot \Vert_{\gamma(K,E)})$ is a separable Banach space, see \cite[Section 4 in Chapter 1]{Neidhardt_1978}. Note that if $K= \mathbb{R}$, then $\gamma(\mathbb{R},E)$ can be identified with $E$.

\subsection{Various types of measurability}\label{sec-measurability}

Before we continue we need to recall some basic definitions about various notions of measurability.
For this purpose let us fix $(Z,\mathscr{Z})$ be a measurable space and $X$ be a Banach space. Let us denote the Borel $\sigma$-algebra of $X$ by $\mathscr{B}(X)$. The names we use are slightly untypical, since we use the letter $\mathscr{Z}$ to underline the $\sigma$-field in question.

\begin{definition}\cite[Sec 1.1., p. 2]{HVVW2016B}\label{def-BorelMeas}
Suppose that $\mathscr{X}$ is a $\sigma$-field of subsets of $X$. A function $f: Z \to X$ is $\mathscr{Z}$/$\mathscr{X}$  measurable, iff the pre-image $f^{-1}(B)$ belongs to $\mathscr{Z}$ for every set  $B \in \mathscr{X}$.
\\	In particular,  a function $f: Z \to X$ is $\mathscr{Z}$-Borel iff it is  $\mathscr{Z}$/$\mathscr{B}(X)$  measurable.
\end{definition}
\begin{definition}\cite[Defn 1.1.3]{HVVW2016B}
	A function $f: Z \to X$ is called $\mathscr{Z}$-simple iff it is of the form $f = \sum_{i=1}^{N} \mathds{1}_{A_i} \otimes x_i$ for some $N \in \mathbb{N}$,  $A_i \in \mathscr{Z}$ and $x_i \in X$ for all $i=1, \ldots,N$.  Here $(\mathds{1}_{A_n} \otimes x_n )(s) = \mathds{1}_{A_n}  (s) x_n$.
\end{definition}

\begin{definition}\cite[Defn 1.1.4]{HVVW2016B}\label{def-strongMeas}
	A function $f: Z \to X$ is said to be strongly $\mathscr{Z}$-measurable iff there exists a sequence of $\mathscr{Z}$-simple functions $f_n : Z \to X$ such that $\lim\limits_{n \to \infty} f_n = f$ pointwise on $Z$.
\end{definition}

\begin{remark}\cite[Corollary 1.1.10]{HVVW2016B}\label{rem-measurability}
	If $X$ is a separable Banach space,  then a function $f: Z \to X$  is strongly $\mathscr{Z}$-measurable iff it  is $\mathscr{Z}$/$\mathscr{B}(X)$  measurable.
\end{remark}

When $X$ is a space of linear operators, the situation becomes more involved. In particular,
the pronoun ``strongly" can be   used in two different ways. Hence we have to be careful.
It is well-known,  see   \cite[Lemma 44, p. 51]{Neidhardt_1978}, that many classical $\sigma$-fields of subsets of the space  $\gamma(K,E)$ coincide.
In particular, the Borel $\sigma$-field $\mathscr{B}(\gamma(K,E))$, i.e. the  $\sigma$-field generated by $\topp(\gamma(K,E))$, i.e. the topology  on $\gamma(K,E)$ which is induced by the  norm \eqref{eqn-norm gamma(K,E)} on
 $\gamma(K,E)$,  is equal to the strong $\sigma$-field $\mathcal{S}(\gamma(K,E))$, i.e. the  $\sigma$-field generated by the following family  of strong open sets
 \begin{equation}\label{eqn-strong set}
 \strong\topp(\gamma(K,E)):=\bigl\{ \{L \in \gamma(K,E): Lk \in A \}, \;\; k\in K, \; A \in \topp(E)\bigr\},
 \end{equation}
 where $\topp(E)$ is the topology on $E$ which is induced by the  norm on
 $E$. In other words,
 \begin{equation}\label{eqn-Neidhard Lemma 44}
 \mathscr{B}(\gamma(K,E))=\mathcal{S}(\gamma(K,E)).
   \end{equation}

To formulate  our main result in this part of our work we need the following definitions of strong and Borel measurability   of an $\gamma(K, E)$-valued function.
\begin{definition}\cite[Defn 1.1.27]{HVVW2016B}\label{def-opValStrongMeas}
	A function $f : Z \to \gamma(K, E)$ is called $\mathscr{Z}$-double-strongly measurable if for all $k \in K$ the $E$-valued function $f(\cdot)k : Z \ni z\mapsto f(z)(k) \in E$ is strongly $\mathscr{Z}$-measurable according to   Definition \ref{def-strongMeas}.
\end{definition}

\begin{definition}\label{def-opValBorelMeas}
	A function $f : Z \to \gamma(K, E)$ is called strongly $\mathscr{Z}$-Borel measurable if for all $k \in K$ the $E$-valued function $f(\cdot)k : Z \ni z\mapsto f(z)(k) \in E$ is $\mathscr{Z}$-Borel  according to   Definition \ref{def-BorelMeas}.
\end{definition}

\begin{proposition}\label{prop-Neidhard Lemma 44}
For a function $f:Z \to \gamma(K,E)$  the following conditions are equivalent.
\begin{trivlist}
\item[(i)] $f$ is  $\mathcal{Z}$/$\mathscr{B}(\gamma(K,E))$ measurable, see Definition \ref{def-BorelMeas};
\item[(ii)] $f$ is  strongly $\mathscr{Z}$-measurable,  see Definition \ref{def-strongMeas};
\item[(iii)] $f$ is strongly $\mathscr{Z}$-Borel measurable,  see Definition \ref{def-opValBorelMeas};
\item[(iv)] $f$ is  $\mathscr{Z}$-double strongly  measurable,  see Definition \ref{def-opValStrongMeas};
\item[(v)] $f$ is $\mathcal{Z}$/$\mathcal{S}(\gamma(K,E))$  measurable, see Definition \ref{def-BorelMeas}.
\end{trivlist}
\end{proposition}

It is well-known, that such a result does not hold if the space $\gamma(K,E)$ is replaced by the space $\mathcal{L}(K,E)$ of all bounded and linear operators from $K$ to $E$.

\begin{proof}[\textbf{Proof of Proposition \ref{prop-Neidhard Lemma 44}}]
Since $X:= \gamma(K,E)$ and $E$  are  separable Banach spaces, by Remark \ref{rem-measurability} we infer that
(i)$\Leftrightarrow$ (ii) and  (iii)$\Leftrightarrow$ (iv). Moreover, (iv)$\Leftrightarrow$ (v). For this, let us observe that $k\in K$, then $f(\cdot)k=f\circ i_{k}$, where $i_k:\gamma(K,E)\ni L \mapsto L(k)\in E$ is the evaluation map.
Thus, for every $A\in \topp(E)$ and every $k \in K$, $f^{-1} (\{L \in \gamma(K,E): L\circ i_k \in A \})=(f(\cdot)k)^{-1}(A)$. Hence, the equivalence (iv)$\Leftrightarrow$ (v) follows. We conclude the proof by observing that the equivalence (i)$\Leftrightarrow$ (v)
is a consequence of Neidhardt's \cite{Neidhardt_1978} result \eqref{eqn-Neidhard Lemma 44}. This completes the proof.
\end{proof}

Let us recall that the progressive $\sigma$-field
$\mathscr{B}\mathbb{F}(\mathbb{R}_+ \times \Omega)$  consists of all sets  $A \subseteq \mathbb{R}_+ \times \Omega$ such that for any $t \in \mathbb{R}_+$, the set $A \cap ([0,t] \times \Omega) \in \mathscr{B}([0,t]) \otimes \mathscr{F}_t$, see \cite[Definition I.4.7]{RY1999B} and \cite[6.0.4 and 6.0.5]{Wentzel_1981}. Similarly we can define the $\sigma$-field $\mathscr{B}\mathbb{F}([0,T] \times \Omega)$.

Specializing Proposition \ref{prop-Neidhard Lemma 44} to $Z=[0,T]\times \Omega$ and $\mathcal{Z}=\mathscr{B}\mathbb{F}([0,T]\times \Omega)$, we infer that all five possible definitions  of $\mathscr{B}\mathbb{F}([0,T] \times \Omega)$-measurability  of   processes taking values in $\gamma(K,E)$ are equivalent. In what follows we will simply use words ``a progressively measurable $\gamma(K,E)$-valued process".

Let us point out that there are even more versions of this definition. For instance,
a process $\eta:[0,T]\times \Omega \to \gamma(K,E)$  is  progressively Borel  measurable iff
$\eta^{-1}(A) \in \mathscr{B}\mathbb{F}([0,T] \times \Omega)$ for every $A \in \mathscr{B}(\gamma(K,E) )$.

In a more general case, when $\gamma(K,E)$  is replaced by a separable Banach space $X$, we can only use
Remark \ref{rem-measurability} to assert that a process $\xi:[0,T]\times \Omega \to X$
is Borel progressively  measurable iff it is   strongly progressively measurable.

In addition to measurable and strongly measurable functions, one can consider a notion of strongly measurable functions with respect to a measure. More precisely, we have, see \cite{Yosida_1980} and \cite[Definitions 1.1.13 and 1.1.14]{HVVW2016B}.

\begin{definition}\label{def-measure measurability}
Assume that $(Z,\mathcal{Z},\mu)$ is a measure space and $X$ is a Banach space.
\begin{trivlist}
\item[(i)] A function $f: Z \to X$ is called $\mu$-simple iff it is of the form $f = \sum_{i=1}^{N} \mathds{1}_{A_i} \otimes x_i$ for some $N\in \mathbb{N}$,  $A_i \in \mathscr{Z}$ with $\mu(A_i)<\infty$  and $x_i \in X$ for all $i=1, \ldots, N$.
\item[(ii)] A function $f: Z \to X$ is said to be strongly $\mu$-measurable iff there exists a sequence of $\mu$-simple functions $f_n : Z \to X$ such that $ f_n \to  f$, $\mu$-almost everywhere.
\end{trivlist}
\end{definition}

The next result is borrowed from \cite[Proposition 1.1.16 and Remark 1.1.18]{HVVW2016B}.
\begin{proposition} \label{prop-equivalence measurabilities}
Assume that $(Z,\mathcal{Z},\mu)$ is a $\sigma$-finite measure space and $X$ is a separable Banach space. Then the following three assertions are equivalent.
\begin{trivlist}
\item[(i)] $f$ is strongly $\mu$-measurable;
\item[(ii)] there exists a strongly $\mathcal{Z}$-measurable function $\tilde{f}:Z \to X$ such that
$f=\tilde{f}$,  $\mu$-almost everywhere;
\item[(iii)] $f$ is strongly $\mathcal{Z}_{\mu}$-measurable, where $\mathcal{Z}_{\mu}$ is the completion of $\mathcal{Z}$ with respect to ${\mu}$.
\end{trivlist}
\end{proposition}

Let us also recall the following definition of Bochner spaces.

\begin{definition}\label{def-L^p space of strongly measurable functions}
Assume that $(Z,\mathcal{Z},\mu)$ is a measure space, $X$ is a Banach space and $p\in [1,\infty)$. By $\mathcal{L}^p(Z,\mathcal{Z};X)=
\mathcal{L}^p(Z;X)$ we denote the vector space of all strongly $\mu$-measurable functions $f:Z \to X$ such that
\begin{equation}\label{eqn-L^p norm}
  \int_Z \vert f(z) \vert_{X}^p\, \mu(dz)<\infty.
\end{equation}
 By $\tilde{\mathcal{L}}^p(Z,\mathcal{Z};X)=
\tilde{\mathcal{L}}^p(Z;X)$ we denote the vector space of all strongly $\mathcal{Z}$-measurable functions $f:Z \to X$ such that
condition \eqref{eqn-L^p norm} holds. \\
By ${L}^p(Z,\mathcal{Z};X)=
{L}^p(Z;X)$ we denote the vector space of all equivalence classes of functions from $\mathcal{L}^p(Z;X)$.
By $\tilde{{L}}^p(Z,\mathcal{Z};X)=
\tilde{{L}}^p(Z;X)$ we denote the vector space of all equivalence classes of functions from $\tilde{\mathcal{L}}^p(Z;X)$.
\end{definition}
Proposition \ref{prop-equivalence measurabilities} implies that the spaces ${L}^p(Z,\mathcal{Z};X)$ and $\tilde{{L}}^p(Z,\mathcal{Z};X)$
are naturally isometrically isomorphic.   In particular,  for every element  $g$ of   ${L}^p(Z,\mathcal{Z};X)$, there exists a strongly $\mathcal{Z}$-measurable function $\tilde{f}:Z \to X$ such that $[\tilde{f}]=g$. Endowed with the classical norm, the Bochner  space ${L}^p(Z;X)$ is a Banach space.

When $X=\gamma(K,E)$, where $K$ and $E$ are separable Hilbert and, respectively, Banach spaces, it follows from the above and Proposition \ref{prop-Neidhard Lemma 44} that
the space ${L}^p(Z,\mathcal{Z};X)$  can be defined as the space of all equivalence classes of $\mathcal{Z}$/$\mathscr{B}(\gamma(K,E))$-measurable functions such that
\begin{equation}\label{eqn-L^p norm-gamma(K,E)}
  \int_Z \vert f(z) \vert_{\gamma(K,E)}^p\, \mu(dz)<\infty.
\end{equation}
In fact, the $\mathcal{Z}$/$\mathscr{B}(\gamma(K,E))$-measurability can be replaced by each of the five versions of measurability listed in
Proposition \ref{prop-Neidhard Lemma 44}.

\subsection{Stopping times}\label{sub-stopping times}

In this section, as throughout the whole paper,  we assume that  $(\Omega, \mathscr{F},\tP, \mathbb{F})$, where $\mathbb{F} := \{ \mathscr{F}_t: t \geq 0 \}$,  is  a filtered probability space which satisfies the \textbf{usual hypothesis}, that is, the filtration $\mathbb{F}$ is right continuous and the $\sigma$-field $\mathscr{F}_0$ contains all $\tP$-null sets of $\mathscr{F}$, see \cite[Definition I.1.1]{Metivier_1982}.
According to  \cite[Definition I.4.1]{Metivier_1982}, a function $\tau: \Omega \to [0,\infty]$ is an $\mathbb{F}$-stopping time iff for every $t \in \mathbb{R}_+$, the set $\{\tau \leq t\}:=\{ \omega \in \Omega: \tau(\omega) \leq t\}$ belongs to the $\sigma$-field $\mathcal{F}_t$.

\begin{definition}\label{def-accessible stopping time}
A stopping time  $\tau$ is called accessible, see e.g. \cite[section 2.1, p. 45]{Kunita-90}, iff there exists an increasing sequence of stopping times $\{\tau_n\}_{n \in \mathbb{N}}$ with the following properties: \begin{enumerate}
		\item $\lim\limits_{n \to \infty}\tau_n = \tau$, $\tP$-a.s.,
		\item for every $n$, $\tau_n < \tau$, $\tP$-a.s. on $\{\tau >0\}$.
\end{enumerate} For such sequence we write $\tau_n \nearrow \tau$. Such a sequence $\{\tau_n\}_{n \in \mathbb{N}}$ will be called an announcing sequence for the accessible stopping time $\tau$.
\end{definition}

Let us point out that in \cite[Definition IV.5.4]{RY1999B} a process which we call accessible is called predictable. On the other hand, Metivier in \cite[Definition I.4.9]{Metivier_1982} gives a different definition of a predictable stopping time. Fortunately, according to \cite[Theorem I.6.6]{Metivier_1982}, $\tau$ is a predictable stopping time according to \cite[Definition I.4.9]{Metivier_1982} if and only if
$\tau$ is an accessible according to our definition, provided the usual hypothesis satisfies, see \cite[Definition I.1.1]{Metivier_1982}. Let us point out that our standing assumption is that the
the  filtered probability  space $(\Omega, \mathscr{F},\tP, \mathbb{F})$ satisfies the usual hypothesis.
Therefore, in our paper, the notions of accessible and predictable stopping times are equivalent and this allows us to use  later on \cite[Proposition I.4.14]{Metivier_1982}.

For any given stopping time $\tau$, we set
\begin{align}
\label{eqn-Omega_tau}
	\Omega_t(\tau) &:= \{ \omega \in \Omega: t < \tau(\omega)\}, \\
\label{eqn-stochastic interval open}
[0,\tau) \times \Omega  &:= \{(t,\omega) \in [0,\infty) \times \Omega: 0 \leq t < \tau(\omega)  \},\\
 [0,\tau] \times \Omega  &:= \{(t,\omega) \in [0,\infty) \times \Omega: 0 \leq t \leq  \tau(\omega)  \}.
\label{eqn-stochastic interval closed}
\end{align}
Note that the sets $ [0,\tau) \times \Omega$ and $ [0,\tau] \times \Omega$ are sometimes denoted by $ [[0,\tau))$ and $ [[0,\tau]]$ respectively. They are also sometimes called  "stochastic intervals".
It is useful to observe that for two stopping times $\tau_1$ and $\tau_2$ and $t>0$,  the following equality holds:
\[\Omega_t(\tau_1) \cap \Omega_t(\tau_2)=\Omega_t(\tau_1 \wedge \tau_2),
\]
where $(\tau_1 \wedge \tau_2)(\omega):=\min\{ \tau_1(\omega), \tau_2(\omega)\}$, $\omega \in \Omega$.

To prove the uniqueness of a local solution we need the following criteria of equivalent processes.

\begin{definition}\label{defn-EquivSP}
	Assume that $X$ is a separable Banach space. A local $X$-valued stochastic process is a function $\xi : [0,\tau) \times \Omega \to X$, where $\tau$ is an accessible stopping time. \\
Suppose that  $\tau$ is an accessible stopping time with an announcing sequence $\{\tau_n\}_{n \in \mathbb{N}}$.
Then a local stochastic process $\xi: [0,\tau) \times \Omega \to X$ is called $\mathbb{F}$-progressively measurable, see e.g. \cite{Brz+H+R_2020-liquid crystals penalised},  iff   for every $n$, the  stopped $X$-valued process
\[\xi(\cdot \wedge \tau_n):= [0,\infty) \times \Omega \ni (t,\omega) \mapsto \xi(t \wedge \tau_n(\omega),\omega)\in X \] is $\mathbb{F}$-progressively measurable. If the filtration $\mathbb{F}$ is unambiguous  from the context, then we often skip it from using. \\
	Two local stochastic processes $\xi_i : [0,\tau_i) \times \Omega \to X$, $i=1,2$ are called equivalent, we will write $(\xi_1,\tau_1) \sim (\xi_2,\tau_2)$, if and only if  $\tau_1=\tau_2$, $\tP$-a.s. and for any $t >0$ the following holds
\begin{equation*}
	\xi_1(\cdot,\omega) = \xi_2 (\cdot,\omega) ~~ \textrm{ on } [0,t],
	\end{equation*}
for almost all $\omega \in \Omega_t(\tau_1) \cap \Omega_t(\tau_2)$.
\end{definition}

Assume that $p \in [1,\infty)$.	By $\mathcal{M}_{loc}^p(\mathbb{R}_+,E)$,  we denote the space of all progressively measurable  $E$-valued processes $\xi : \mathbb{R}_+ \times \Omega \to E$ for which there exists a sequence $\{\tau_n\}_{n \in \mathbb{N}}$ of bounded stopping times such that $\tau_n \nearrow \infty, ~ \tP$-a.s. and
\begin{equation*}
	\tE\left[ \int_{0}^{\tau_n} \Vert \xi(t)\Vert_E^p \, dt \right] < \infty, ~~ \forall n \in \mathbb{N}.
	\end{equation*}

Assume that $p\in [1,\infty)$, $q \in [1,\infty)$, and  $T>0$. By	 $\mathcal{M}^{p}(\mathcal{L}^{q}([0,T],E))$
we denote the space of all progressively measurable  $E$-valued processes $\xi : [0,T]  \times \Omega \to E$  such that
 \begin{equation*}
	\tE\left[ \left( \int_{0}^{T} \Vert \xi(t)\Vert_E^q \, dt \right)^{\frac{p}{q}} \right] < \infty.
	\end{equation*}

The space $\mathcal{M}^{p}(\mathcal{L}^{\infty}([0,T],E))$ is defined analogously with the ``norm" $\bigl(\int_{0}^{T} \Vert \xi(t)\Vert_E^q \, dt \bigr)^{\frac{1}{q}}$ being replaced by  $\esssup_{t \in [0,T]}\Vert \xi(t)\Vert_E$.

As usual, see e.g. \cite[Definition IV.2.1]{RY1999B}, by	 $\mathbb{M}^{p}({L}^{q}([0,T],E))$
we denote the space of equivalence classes of elements of  $\mathcal{M}^{p}(\mathcal{L}^{q}([0,T],E))$.
Let us note that  $\mathbb{M}^{p}({L}^{q}([0,T],E))$ is a closed subspace of, typically not equal to, $L^p(\Omega,\mathscr{F};{L}^{q}([0,T],E))$.
To simplify the notation, we will often use the notation ${M}^{p}([0,T],E)$ to denote $\mathbb{M}^{p}({L}^{p}([0,T],E))$.

\section{Inhomogeneous Strichartz estimates}\label{sec-inHomStrichartz}
In this section we will prove the deterministic Strichartz type estimate, see Theorem \ref{thm-InHomStrichartz} below, which is a generalization of \cite[Theorem 1.2]{IJ2011}  and is essential to tackle, both, the Dirichlet and the Neumann boundary case.

Recall that in our setting, the operator $(\nA,\Dom(\nA))$ possesses a complete orthonormal system of eigenvectors $\{ e_j \}_{j \in \mathbb{N}}$ in $L^2({\cO})$. We have denoted the corresponding eigenvalues by $\lambda_j^2$. From the functional calculus of self-adjoint operators, it is known that $\{(e_j, \lambda_j)\}_{j \in \mathbb{N}  }$ is a sequence of the associated eigenvector and eigenvalue pair for $\sqrt{\nA}$. For any integer $\lambda \geq 0 $, $\Pi_\lambda$ is defined as the spectral projection of $L^2({\cO})$ onto the subspace spanned by $\{ e_j\}_{j \in \mathbb{N}}$ for which $\lambda_j\in [\lambda,\lambda+1)$, i.e.
\begin{equation*}
\Pi_\lambda u = \sum_{j=1}^{\infty} \mathds{1}_{[\lambda,\lambda+1)} (\lambda_j) \langle u,e_{j} \rangle_{L^2({\cO})} e_j, \quad u \in L^2({\cO}).
\end{equation*}

At this juncture, it is relevant to note that the proof of the Strichartz estimate in deterministic setting, see e.g. \cite{Burq+L+P_2008} and \cite{Burq+P_2009}, is based on the following estimate in the Lebesgue spaces of the spectral projector $\Pi_\lambda$, refer \cite{SS2007} for the proof.
\begin{theorem}\label{thm-HSCS2007-Thm1.1}
	For any smooth bounded domain ${\cO}\subset \mathbb{R}^2$, the following estimate holds for all $u \in L^2({\cO})$
	\begin{equation}\label{eqn-Pi_lambda}
		\|\Pi_\lambda u\|_{L^q({\cO})}\leq C \lambda^{\rho}\|u\|_{L^2({\cO})},
	\end{equation}
 where
 \begin{equation}\label{eqn-rho}
 	\rho :=  \left\{   \begin{aligned}
 		\frac{2}{3} \bigl(\frac{1}{2}-\frac{1}{q}\bigr) \quad \textrm{ if } 2 \leq q \leq 8, \\
 		2 \bigl(\frac{1}{2}-\frac{1}{q}\bigr) - \frac{1}{2} \quad \textrm{ if } 8 \leq q \leq \infty.
 	\end{aligned} \right.
 \end{equation}
\end{theorem}

Since the  Strichartz estimates below derived, for the homogeneous and inhomogeneous wave equation, holds for both the Dirichlet and the Neumann case, from now onwards, to shorten the notation, we denote $A_{B,q}$ and $A_{B,2}$, respectively, by $A_q$ and $A$.
\begin{theorem}[Deterministic Strichartz Estimates]\label{thm-InHomStrichartz}
	Let us assume that $T>0$. Then there exists a positive constant $C_T$, which is increasing w.r.t. $T$ and may also depend on $p,q,r$, such that the following holds: if $u$ satisfies the following linear inhomogeneous wave equation
	\begin{equation}\label{eqn-wave}\left\{
	\begin{array}{ll}
	u_{tt} - \Delta u = f \textrm{ in } (0,T)\times {\cO} \\
	u(0,\cdot) = u_0(\cdot), ~~ u_t(0,\cdot) = u_1(\cdot),
	\end{array}
	\right.
	\end{equation}	
with either boundary condition \begin{align}
	& \textrm{Dirichlet : }  u \upharpoonright (0,T) \times \partial {\cO}= 0,  \nonumber\\
	& \textrm{Neumann : } \partial_{\nu} u \upharpoonright (0,T) \times \partial {\cO} =0, \nonumber
	\end{align}
	where $\nu$ is the outward normal unit vector to $\partial {\cO}$ and $f \in L^1(0,T; L^2({\cO}))$, then
	\begin{equation}\label{InHomStrichartzIneq}
	\Vert u \Vert_{L^p(0,T;\Dom(\nA_q^{\frac{1-r}{2}}))}  \leq C_T \left[\Vert u_0 \Vert_{\Dom(\nA^{\frac12 })} + \Vert u_1 \Vert_{L^2({\cO})} + \Vert f \Vert_{L^1(0,T;L^2({\cO}))} \right],
	\end{equation}
	for all $(p,q,r)$ which satisfy
	\begin{equation}\label{InHomStrichartzIneqCondn}
	2 \leq q \leq p < \infty,  ~~  \textrm{ and }  r = \left\{   \begin{aligned}
		\frac{5}{6} - \frac{1}{p} - \frac{2}{3q} \quad \textrm{ if } 2 \leq q \leq 8, \\
		1 - \frac{1}{p} - \frac{2}{q} \quad \textrm{ if } 8 \leq q < \infty.
		\end{aligned} \right.
	\end{equation}
\end{theorem}

\begin{remark}\label{rem-L^infty}
In addition to the Strichartz estiamtes \eqref{InHomStrichartzIneq} one also has the classical (including those on the velocity $u_t$) estimates, see \cite{Ball_1977} for the inhomogeneous part,
	\begin{equation}\label{ineq-L^inftyH}
	\Vert u \Vert_{\mathrm{C}([0,T];\Dom(\nA^{\frac12 }))} +\Vert u_t \Vert_{\mathrm{C}([0,T];H)}  \leq \bar{C}_T \left[\Vert u_0 \Vert_{\Dom(\nA^{\frac12 })} + \Vert u_1 \Vert_{L^2({\cO})} + \Vert f \Vert_{L^1(0,T;L^2({\cO}))} \right].
	\end{equation}
We can assume that the constants $\bar{C}_T$  increase w.r.t. $T$.
In a standard way, this inequality can be lifted to more regular data as follows, for any $k\geq 0$,
\begin{equation}\label{ineq-L^infty-regular}
	\Vert u \Vert_{\mathrm{C}([0,T];\Dom(\nA^{k+\frac12}))}+ \Vert u_t \Vert_{\mathrm{C}([0,T];\Dom(\nA^{k}))}  \leq \bar{C}_T \left[\Vert u_0 \Vert_{\Dom(\nA^{k+\frac12})} + \Vert u_1 \Vert_{\Dom(\nA^{k })} + \Vert f \Vert_{L^1(0,T;\Dom(\nA^{k })} \right].
	\end{equation}
\end{remark}

\begin{remark} Let us observe that if for $T>0$, $C_T$ denotes the smallest constant for which  inequality \eqref{InHomStrichartzIneq} holds for all data $u_0,u_1$ and $F$ from appropriate spaces, then the function
\[
(0,\infty) \ni T\mapsto C_T \in (0,\infty),
\]
is non-decreasing (or weakly increasing as some people call). Similar results hold true for $\bar{C}_T$.
\end{remark}
\begin{remark}\label{rem-Ibrahim-Jrad}
Ibrahim and Jrad proved, see inequality (3.5) in the proof of \cite[Theorem 1.2]{IJ2011}, that
	
\begin{equation}\label{eqn-Ibrahim-Jrad}
	\Vert u \Vert_{L^8(0,T;H^{\frac38,8}({\cO})}  \leq C_T \left[\Vert u_0 \Vert_{\Dom(\nA^{\frac12 })} + \Vert u_1 \Vert_{L^2({\cO})} + \Vert F \Vert_{L^1(0,T;L^2({\cO}))} \right],
	\end{equation}
for $A=-\Delta_{D}$. Since the substitution of $p=q=8$ into \eqref{InHomStrichartzIneq} gives \eqref{eqn-Ibrahim-Jrad}, our result generalizes \cite[Theorem 1.2]{IJ2011}. Note that in \cite{IJ2011} the space $H^{\frac38,8}({\cO})$ is denoted, inconsistently   with current approach, by $W^{\frac38,8}({\cO})$. Finally, let us point out that the inequality (1.8) with the H\"older space $\dot{C}^{\frac18}({\cO})$  in \cite[Theorem 1.2]{IJ2011}
is a consequence of the Sobolev embedding, see e.g. \cite[Theorem 2.8.1 (e) and Definition 1 (d) in 2.3.1]{Triebel_1978B}.
\end{remark}

\begin{proof}[\textbf{Proof of Theorem \ref{thm-InHomStrichartz}}] Without loss of generality we assume that $T=2\pi$. The proof is divided into two cases. In the first case, we derive the Strichartz estimate for the homogeneous problem (i.e. $F=0$) and then, in second case, we prove the inhomogeneous one (i.e. $F \neq 0$) by using the homogeneous estimate from first case.

\textbf{The first case: Estimate for the homogeneous problem.} In this case, the Duhamel formula gives
	\begin{equation}\label{InHomStrichartz-t1}
	u(t) = \cos(t{\sqrt{\nA}})u_0+\frac{\sin(t{\sqrt{\nA}})}{{\sqrt{\nA}}}u_1,
	\end{equation}
	where, from the functional calculus for self-adjoint operators, for each $t$, $\cos(t{\sqrt{\nA}})$ and $\frac{\sin(t{\sqrt{\nA}})}{\sqrt{\nA}}$ are well-defined bounded operators on $L^2({\cO})$. Moreover, we have
	\begin{equation*}
	\cos(t{\sqrt{\nA}}) = \frac{e^{it{\sqrt{\nA}}}+e^{-it{\sqrt{\nA}}}}{2} .
	\end{equation*}
Let $\mathcal{L}_{\pm}(t)u_0 := e^{\pm it\sqrt{\nA}}u_0$ be the solution $u$ of $\partial_t u = \pm i \sqrt{\nA} u$ such that $u(0) = u_0$. In other words, $\mathcal{L}_{\pm}=\bigl(\mathcal{L}_{\pm}(t))_{t \geq 0}$ is $C_0$-group with the generator $\pm i\sqrt{\nA}$.
Using the Minkowski's inequality we get
\begin{align}\label{InHomStrichartz-t3}
	\Vert u\Vert_{L^p(0,T;\Dom(\nA_q^{\frac{1-r}{2}}))} & \lesssim \Vert e^{it\sqrt{\nA}}u_0 \Vert_{L^p(0,T;\Dom(\nA_q^{\frac{1-r}{2}}))} \\
	&\quad  + \Vert e^{-it\sqrt{\nA}}u_0 \Vert_{L^p(0,T;\Dom(\nA_q^{\frac{1-r}{2}}))}  \quad + \bigg\Vert \frac{\sin(t{\sqrt{\nA}})}{{\sqrt{\nA}}} u_1 \bigg\Vert_{L^p(0,T;\Dom(\nA_q^{\frac{1-r}{2}}))}.
\nonumber	\end{align}
Therefore, it is enough to estimate, as done in the following Steps 1-4, the $L^p(0,T;\Dom(\nA_q^{\frac{1-r}{2}}))$-norm of $e^{it\sqrt{\nA}}u_0 $ and $\frac{\sin(t{\sqrt{\nA}})}{{\sqrt{\nA}}} u_1$. We will write the variables in subscript, wherever required, to avoid any confusion.

\noindent \textbf{Step 1. } Here we show that \begin{equation}\label{a1}
\|e^{itB}u_0\|_{L_t^p(0,2\pi; L_x^q({\cO}))}\leq C \|u_0\|_{\Dom(\nA^{r/2})},
\end{equation}
where $B$ is the following  \enquote{modification} of ${\sqrt{\nA}}$ operator by considering only the integer eigenvalues, i.e.
\begin{align}\label{eqn-B} 	B(e_j)=[\lambda_j]e_j, ~~ j \in \mathbb{N}. \nonumber
\end{align}
	The notation $[\cdot]$ stands for the integer part and $e_j$ is an eigenfunction of $\nA$ associated to the eigenvalue $\lambda_j^2$. Before moving further we prove the boundedness property of the operator $B-\sqrt{\nA}$.

\begin{lemma}\label{lem\nA-B}
	For every $r \in [0,1]$, the operator $B-\sqrt{\nA}$ is bounded  on $\Dom(\nA^{\frac r2 })$.
\end{lemma}
\begin{proof}[Proof of Lemma \ref{lem\nA-B}]
		Let us fix $r \in [0,1]$. Observe that by definition of $B$ we have for every $u  \in \Dom(\nA^{\frac{r}{2}})$,
\begin{equation*}
		(B-\sqrt{\nA})u = \sum_{j \in\mathbb{N}} \{ \lambda_j \} \langle u,e_j \rangle_{L^2({\cO})} e_j,
		\end{equation*}
		where $\{ \lambda_j \} := \lambda_j - \left[ \lambda_j\right]$ is the fractional part of $\lambda_j$. Then
		\begin{align}
		\| (B-\sqrt{\nA}) u\|_{L^2({\cO})}^2 & \leq \sum_{j \in \mathbb{N}}  \{ \lambda_j\}^2~ |\langle u,e_j \rangle_{L^2({\cO})}|^2 \leq  \| u\|_{L^2({\cO})}^2  \leq  \| u\|_{\Dom(A^{\frac{r}{2}})}^2. \nonumber
		\end{align}
		Moreover,
		\begin{align}
		\| \nA^{\frac{r}{2}}(B-\sqrt{\nA}) u\|_{L^2({\cO})}^2 =  \sum_{j \in \mathbb{N}} \lambda_j^{2r} \{ \lambda_j\}^2~ |\langle u,e_j \rangle_{L^2({\cO})}|^2 \leq \sum_{j \in \mathbb{N}} \lambda_j^{2r}  |\langle u,e_j \rangle_{L^2({\cO})}|^2  = \| \nA^{\frac{r}{2}} u\|_{L^2({\cO})}^2. \nonumber
		\end{align}
		Hence, by the definition of norm in $\Dom(\nA^{\frac{r}{2}})$ we have
		\begin{align}
		\| (B-\sqrt{\nA}) u\|_{\Dom(\nA^{\frac{r}{2}})}^2 = \| (B-\sqrt{\nA}) u\|_{L^2({\cO})}^2 + \|\nA^{\frac{r}{2}}(B-\sqrt{\nA}) u\|_{L^2({\cO})}^2 \lesssim \| u\|_{\Dom(\nA^{\frac{r}{2}})}^2.  \nonumber
		\end{align}
\end{proof}

\noindent In continuation of the proof of \eqref{a1}, since $u_0 \in L^2({\cO})$, we can write
$$	u_0=\sum_{j \in \mathbb{N}} \langle u_0, e_j \rangle_{L^2({\cO})} e_j =: \sum_{j \in\mathbb{N}} u_j e_j.	$$	
By functional calculus for self-adjoint operators,
\begin{equation}
	e^{itB}u_0(x) = \sum_{j \in \mathbb{N}}e^{it[\lambda_j]}u_j e_j(x) =: \sum_{k \in \mathbb{N}} u_k(t,x), \nonumber
\end{equation}
where, since $[\lambda_j] = k$,
\begin{align}
	u_k(t,x) & = \sum_{j \in \mathbb{N}} \mathds{1}_{[k,k+1)} (\lambda_j) ~ e^{itk}u_j e_j(x) = e^{itk}\Pi_k u_0(x). \nonumber
\end{align}

Assume that   $0 < \theta < \frac{1}{2}$. Let $C := -D_{tt} + I$ on $L^2(0,2 \pi)$ with the periodic boundary conditions. Then, see \cite{Triebel_1978B},
 $H^{\theta,2} = \Dom(C^{\theta/2}) $ with equivalent norms. The norm $\Dom(C^{\theta/2}) $ is Hilbertian with the corresponding inner product
$$\langle u,v \rangle_{\Dom(C^{\theta/2})} :=  \langle C^{\theta/2} u,C^{\theta/2} v \rangle_{L^2(0,2\pi)}. $$
We claim that the sequence $\{e^{i\cdot k} \}_{k \in\mathbb{N}}$, after normalization, is an orthonormal basis (ONB) in $\Dom(C^{\theta/2}) $.
Indeed, since the sequence $\{e^{i\cdot k} \}_{k \in\mathbb{N}}$, consisting of  eigenvectors of $C$, is an ONB of  $L^2(0,2\pi)$, we infer that for all $j,k\in \mathbb{N}$,
 \begin{align*}
	\langle e^{i\cdot j} ,e^{i\cdot k} \rangle_{\Dom(C^{\theta/2})}  = \langle C^{\theta/2} e^{i\cdot j}  ,C^{\theta/2} e^{i\cdot k}  \rangle_{L^2(0,2\pi)} = (1+|j|^2)^{\theta/2} (1+|k|^2)^{\theta/2}   \langle e^{i\cdot j}  , e^{i\cdot k}  \rangle_{L^2(0,2\pi)}.
\end{align*}
Hence the claim follows.\\
Moreover, there exists $C_\theta>0$ such that

\begin{equation}\label{eqn-norm of e^itk}
 \left\|e^{i\cdot k}\right\|_{H^{\theta,2}(0,2\pi)}  \leq C_\theta  (1+k)^{\theta},\;\; k\in \mathbb{N}.
\end{equation}
Therefore, for a fixed $x \in D$, we have
\[
e^{itB}u_0(x)=\sum_{k \in \mathbb{N}} \frac{e^{itk}}{ \Vert e^{i\cdot k} \Vert } \Vert e^{i\cdot k} \Vert \Pi_k u_0(x),\;\; t\in (0,2\pi).
\]
Hence,
\begin{align*}
\Vert e^{i\cdot B}u_0(x)\Vert_{\Dom(C^{\theta/2})}^2=\sum_{k \in \mathbb{N}}   \Vert e^{it k} \Vert_{\Dom(C^{\theta/2})}^2 \vert \Pi_k u_0(x) \vert^2.
\end{align*}
Therefore,
\begin{align*}
\Vert e^{i\cdot B}u_0(x)\Vert_{\Dom(C^{\theta/2})}^2 &\leq \sum_{k \in \mathbb{N}} \Vert e^{i\cdot k} \Vert_{\Dom(C^{\theta/2})}^2 \vert  \Pi_k u_0(x) \vert^2
\leq C_\theta^2 \sum_{k \in \mathbb{N}}    (1+k)^{2\theta} \vert  \Pi_k u_0(x) \vert^2.
\end{align*}
Using the equivalence of the two norms we deduce that
\begin{align*}
\Vert e^{i\cdot B}u_0(x)\Vert_{H_t^{\theta,2}}^2 &\leq  C^\prime_\theta \sum_{k \in \mathbb{N}}    (1+k)^{2\theta} \vert  \Pi_k u_0(x) \vert^2.
\end{align*}
Thanks to the 1D Sobolev embedding and Lemma \ref{lem-EquiSpace1}, we have $$H^{\frac{1}{2} - \frac{1}{p},2}(0,2\pi) \hookrightarrow L^p(0,2\pi) \textrm{ for all } p \in [2,\infty), $$
where the space $H^{\theta,2}(0,2\pi)$, $\theta>0$,  is defined by formula \eqref{eqn-H^s,q(D) space}, with set ${\cO}$ replaced by the interval
$(0,2\pi)$. Consequently we argue as follows:
\begin{align}\label{a1-t3}
	& \Vert e^{itB}u_0 \Vert_{L_x^q({\cO};L_t^p(0,2\pi))}^2
	=\bigl(\int_{{\cO}} \Vert e^{itB}u_0(x) \Vert_{L_t^{p}(0,2\pi)}^q \, dx \bigr)^{\frac{2}{q}}
	\nonumber \\
	&\lesssim \bigl(\int_{{\cO}} \Vert e^{itB}u_0(x) \Vert_{H_t^{{\frac{1}{2} - \frac{1}{p}},2}(0,2\pi)}^q \, dx \bigr)^{\frac{2}{q}}  = \bigl\Vert \;\|e^{itB}u_0(x)\|_{H_t^{{\frac{1}{2} - \frac{1}{p}},2}(0,2\pi)}^2\bigr\Vert_{L_x^{\frac{q}{2}}({\cO})}.
\end{align}
Hence, by applying the last Claim with $\theta =\frac{1}{2} - \frac{1}{p}$ in view of  the Parseval formula we deduce that, for fixed $x$,
\begin{align}\label{a1-t4}
	\Vert e^{itB}u_0(x) \Vert_{H_t^{\frac{1}{2} - \frac{1}{p},2}(0,2\pi)}^2 &  \lesssim  \sum_{k \in \mathbb{N}}(1+k)^{1 - \frac{2}{p}}\|u_k(t,x)\|_{L_t^2(0,2\pi)}^2.
\end{align}
Combining the estimate \eqref{a1-t4} and \eqref{a1-t3} followed by Minkowski's inequality and Theorem \ref{thm-HSCS2007-Thm1.1} we obtain
\begin{align}\label{a1-t5}
	\Vert & e^{itB}u_0 \Vert_{L_x^q({\cO};L_t^p(0,2\pi))}^2  \lesssim  \left\| \sum_{k \in \mathbb{N}}(1+k)^{1 - \frac{2}{p}}\|u_k\|_{L_t^2(0,2\pi)}^2 \right\|_{L_x^{\frac{q}{2}}({\cO})}\nonumber\\
	& \leq \sum_{k \in \mathbb{N}} (1+k)^{1 - \frac{2}{p}}\| u_k\|_{L_t^2(0,2\pi; L_x^q({\cO}))}^2  \lesssim \sum_{k \in \mathbb{N}} (1+k)^{1 - \frac{2}{p}}\| \Pi_k u_0\|_{L_x^q({\cO})}^2 \nonumber\\
	&\lesssim \sum_{k \in \mathbb{N}} (1+k)^{1 - \frac{2}{p}}k^{2\rho}  \|\Pi_k u_0\|_{L_x^2({\cO})}^2   \lesssim   \sum_{k \in \mathbb{N}} (1+k)^{1-\frac{2}{p}+ 2\rho} \sum_{j \in \mathbb{N}} \mathds{1}_{[k,k+1)} (\lambda_j) |\langle u_0, e_j \rangle_{L^2({\cO})}|^2 \nonumber\\
	& = \sum_{j \in \mathbb{N}} (1+[\lambda_j])^{1-\frac{2}{p}+ 2\rho} |\langle u_0, e_j \rangle_{L^2({\cO})}|^2
	 = \Vert u_0 \Vert_{\Dom(B^{r})}^2 \simeq \Vert u_0 \Vert_{\Dom(\nA^{r/2})}^2,
\end{align}
where, from $\rho$ in Theorem \ref{thm-HSCS2007-Thm1.1}, we have\footnotemark,
\begin{equation*}
r := \frac{1}{2} - \frac{1}{p} + \rho = \left\{   \begin{aligned}
\frac{5}{6} - \frac{1}{p} - \frac{2}{3q} \quad \textrm{ if } 2 \leq q \leq 8, \\
1 - \frac{1}{p} - \frac{2}{q} \quad \textrm{ if } 8 \leq q < \infty.
\end{aligned} \right.
\end{equation*}
\footnotetext{Note that $r< \frac34$ in the case $2 \leq q \leq 8$ and $r<1$ in the complimentary case $8 \leq q < \infty$.}

Here it is important to highlight that, the equivalence $\Vert u_0 \Vert_{\Dom(B^{r})} \simeq \Vert u_0 \Vert_{\Dom(\nA^{r/2})}$ holds in the last step of \eqref{a1-t5}, because  $\Dom(\nA) = \Dom(B^2)$ and the spaces $\Dom(B^{r})$  and $\Dom(\nA^{r/2})$, for $r \in [0,1]$, are equal to the complex interpolation spaces, between $L^2({\cO})$ and, respectively,  $\Dom(B^2)$ and $\Dom(A)$, see \cite[Theorem 4.3.3]{Triebel_1978B}.

Next, since $p \geq q$, by the Minkowski inequality we obtain the following desired result
	\begin{equation*}
		\Vert  e^{itB}u_0 \Vert_{L_t^p(0,2\pi;L_x^q({\cO}))} \lesssim \Vert u_0 \Vert_{\Dom(\nA^{r/2})},
	\end{equation*}
what also implies that the operator $e^{itB}$ is continuous from $\Dom(\nA^{r/2})$ to $L_t^p(0,2\pi;L_x^q({\cO}))$.\\

\noindent\textbf{Step 2.} In this step we extend inequality \eqref{a1} to operator $\mathcal{L}_{+}$, i.e. we show that
\begin{equation}\label{a2}
	\|\mathcal{L}_{+}(\cdot) u_0\|_{L_t^p(0,2\pi; L_x^q({\cO}))}\leq C\|u_0\|_{\Dom(\nA^{r/2})}.
	\end{equation}
	Let $v(t)=e^{it{\sqrt{\nA}}}u_0$. It is clear that $v$ satisfies \begin{equation*}\left\{
	\begin{array}{ll}
	(\partial_t-iB)v=(-iB+i{\sqrt{\nA}})v \\
	v|_{t=0}=u_0,
	\end{array}
	\right.
	\end{equation*}
	and, therefore, according to the Duhamel formula
	\begin{equation}\label{a2-t0}
	v(t)=e^{itB}u_0+\int_0^t e^{i(t-s)B}(-iB+i{\sqrt{\nA}})v(s)\, ds.	
	\end{equation}
	If we denote $e^{i(t-s)B}(-iB+i{\sqrt{\nA}})v(s,x)$ by $z(s,t,x)$ and $(-iB+i{\sqrt{\nA}})v(s,x)$ by $w(s,x)$, then using the Minkowski inequality, followed by estimate \eqref{a1} and Lemma \ref{lem\nA-B}, we argue as follows:
	\begin{align}\label{a2-t00}
	& \left\| \bigl( \int_{0}^{2\pi} \left[ \int_0^t |z(s,t,x)|  \, ds\right]^p \, dt \bigr)^{\frac{1}{p}}\right\|_{L_x^q({\cO})}  \leq \bigl( \int_{{\cO}} \bigl(  \int_{0}^{2\pi} \| z(s,t,x)\|_{L_t^p(0,2\pi)} \, ds  \bigr)^q \, dx \bigr)^{\frac{1}{q}}  \nonumber\\	
	&  \leq  \int_{0}^{2\pi} \bigl( \int_{{\cO}}  \| z(s,t,x)\|_{L_t^p(0,2\pi)}^q \, dx  \bigr)^{\frac{1}{q}} \, ds   \lesssim	\int_{0}^{2\pi} \left\| v(s,x)\right\|_{\Dom(\nA^{r/2})} \, ds.
	\end{align}
	By putting together \eqref{a2-t0} and \eqref{a2-t00} we obtain
	\begin{align}\label{a2-t1}
	\|v(t,x)\|_{L_x^q({\cO}; L_t^p(0,2\pi))} & \leq   \|e^{itB}u_0(x)\|_{L_x^q({\cO}; L_t^p(0,2\pi))} + \int_{0}^{2\pi} \left\| v(s,x)\right\|_{\Dom(\nA^{r/2})} \, ds \nonumber\\	
	& \leq	\|u_0\|_{\Dom(\nA^{r/2})} +  \int_{0}^{2\pi} \left\| v(s,x)\right\|_{\Dom(\nA^{r/2})} \, ds.		
	\end{align}
	Now, from the boundedness of $e^{it\sqrt{\nA}}$ on $\Dom(\nA^{r/2})$, we infer that \begin{align}\label{a2-t1a}
		\sup_{t \in [0,2\pi]} \left\| e^{it\sqrt{\nA}}u_0 \right\|_{\Dom(\nA^{r/2})} \leq C \| u_0 \|_{\Dom(\nA^{r/2})}.
	\end{align}
	Combining \eqref{a2-t1a} and \eqref{a2-t1} we get \begin{align}
	\|v(t,x)\|_{L_x^q({\cO}; L_t^p(0,2\pi))}& \leq \|u_0\|_{\Dom(\nA^{r/2})} + \int_{0}^{2\pi} \|u_0\|_{\Dom(\nA^{r/2})}\, ds  \lesssim \|u_0\|_{\Dom(\nA^{r/2})}. \nonumber
	\end{align} Hence, again, as an application of the Minkowski inequality we get \eqref{a2} and finish with the proof of Step 2.\\
	
	\noindent \textbf{Step 3.} Here, by using the well-known consequence of Agmon-Douglis-Nirenberg regularity results for the elliptic operators, refer \cite{ADN1959-1964},  we prove the required estimate of the first term in \eqref{InHomStrichartz-t3}, in particular, we show
	\begin{equation}\label{a3}
	\| \mathcal{L}_{+}(\cdot) u_0 \|_{L_t^p(0,2\pi; \Dom(\nA_q^{\frac{1-r}{2}}))} \lesssim \| u_0 \|_{\Dom(\nA^{\frac12 })}.
	\end{equation}

	We start the proof by recalling the following consequence of the Agmon-Douglis-Nirenberg regularity results for the elliptic operators. The operators
	 \begin{align}
	-\Delta_D +I:H^{2,q}({\cO}) \cap H_{D}^{1,q}({\cO})=H^{2,q}({\cO}) \cap {H}_0^{1,q}({\cO}) \to L^q({\cO}), \nonumber
	\end{align}
	and
 	\begin{align}
	-\Delta_N +I:H^{2,q}({\cO})\cap H_{N}^{1,q}({\cO})  \to L^q({\cO}), \nonumber
	\end{align}
	are isomorphisms. These operators will, respectively, be denoted by $A_{D,q}+I$ and $A_{N,q}+I$, or simply by $A_q+I$. Suppose that $u_0 \in D(A^k)$  for sufficiently large $k \in \mathbb{N}$ so  that $Au_0 \in \Dom(\nA^{r/2})$. Then, since the operators $A$ and  $\mathcal{L}_{+}$ commute, we infer that  for all $t \in [0,T]$,
	\begin{align}
	\| \mathcal{L}_{+}(t) u_0 \|_{H^{2,q}({\cO})} & \simeq  \| (\nA+I) \mathcal{L}_{+}(t) u_0 \|_{L^q({\cO})}  = \|  \mathcal{L}_{+}(t)((\nA+I)  u_0) \|_{L^q({\cO})}. \nonumber
	\end{align}
	 Consequently by \eqref{a2}  we get \begin{align}\label{a3-t3}
	& \| \mathcal{L}_{+}(\cdot) u_0 \|_{L_t^p(0,2\pi;H^{2,q}({\cO}))} \lesssim \|(\nA+I) u_0\|_{\Dom(\nA^{r/2})} \sim  \|u_0\|_{\Dom(A^{(r+2)/2})}.
	\end{align}
	Thus, complex interpolation between \eqref{a2} and \eqref{a3-t3} with $\theta= \frac{1-r}{2}$ gives the desired following estimate
	$$\| \mathcal{L}_{+}(\cdot) u_0 \|_{L^p(0,2\pi;\Dom(\nA_q^{\frac{1-r}{2}}))}  \lesssim \| u_0 \|_{\Dom(\nA^{\frac12 })}.$$
	Hence we have completed the proof of Step 3. \\

	\noindent \textbf{Step 4}: Here we incorporate the term with $u_1$, in \eqref{InHomStrichartz-t1}, and complete the proof of the homogeneous Strichartz estimate.
	
	Recall that $\lambda_1=0$ for the Neumann condition and $\lambda_1 >0$ in the Dirichlet case. As mentioned before, we denote by $m_0$ the dimension of eigenspace corresponding to zero eigenvalue. It is known that $m_0=0$ for $\nA= -\Delta_{D}$ and a positive finite integer when $\nA=\Delta_{N}$. To proceed with the proof of this Step, as in \cite{Burq+P_2009}, we single out the contribution of zero eigenvalue and decompose $L^2({\cO})$ into the direct sum of  a finite dimensional space  $\ker A$ and the space orthogonal to $\ker A$, which we denote by $L^{2,+}({\cO})$. Let us observe that if ${\cO}$ is connected, then $\ker A$ is a one dimensional vector space consisting of constant functions. Mathematically, it means, for all $u_1 \in L^2({\cO})$,
	\begin{align}\nonumber
	u_1  = \sum_{j=1}^{m_0} \langle u_1, e_j\rangle_{L^2({\cO})} e_j + \sum_{k > m_0} \langle u_1,e_k \rangle_{L^2({\cO})} e_k    =: \proj u_1 + (\mathds{1}-\proj)u_1 . \nonumber
	\end{align}
	Note that the term $\proj u_1$ does not exist in the Dirichlet condition. Then we argue as follows:
	\begin{align}\label{a4-t1}
	\frac{\sin(t \sqrt{\nA})}{\sqrt{\nA}} u_1 & = \frac{\sin(t \sqrt{\nA})}{\sqrt{\nA}} \proj u_1 + \frac{\sin(t \sqrt{\nA})}{\sqrt{\nA}} (\mathds{1}-\proj) u_1= t \proj u_1 + \frac{\sin(t \sqrt{\nA})}{\sqrt{\nA}} (\mathds{1}-\proj) u_1,
	\end{align}  where the last step holds due to the following argument
	\begin{align}\label{a4-t2}
	& \frac{\sin(t \sqrt{\nA})}{\sqrt{\nA}} u_1 = \sum_{j \in \mathbb{N}} \frac{\sin(t \lambda_j)}{\lambda_j} \langle u, e_j \rangle_{L^2({\cO})} e_j \nonumber\\
	& = t \sum_{j \in \mathbb{N}} \mathds{1}_{\{0\}}(\lambda_j) \langle u, e_j  \rangle_{L^2({\cO})} e_j + \sum_{j \in \mathbb{N}}  \mathds{1}_{(0,\infty)}(\lambda_j) \frac{\sin(t \lambda_j)}{\lambda_j} \langle u, e_j  \rangle_{L^2({\cO})} e_j \nonumber\\
	& = t \proj u_1 + \frac{\sin(t \sqrt{\nA})}{\sqrt{\nA}} (\mathds{1}-\proj)u_1.
	\end{align}
	Now,  since $\bigl(\sqrt{\nA}\bigr)^{-1}$ is isometry from $L^{2,+}({\cO})$ into $\Dom(\nA^{\frac12 })$, by invoking \eqref{a3} on $\bigl(\sqrt{\nA}\bigr)^{-1}$ $((\mathds{1}-\proj) u_1)$ we get
	\begin{align}\label{a4-t3}
	\| \mathcal{L}_{+}(\cdot) \bigl(\sqrt{\nA}\bigr)^{-1}((\mathds{1}-\proj) u_1)) \|_{L^p(0,2\pi;\Dom(\nA_q^{\frac{1-r}{2}}))} & \lesssim \| \bigl(\sqrt{\nA}\bigr)^{-1}((\mathds{1}-\proj) u_1) \|_{\Dom(\nA^{\frac12 })} \nonumber\\
	& = \| (\mathds{1}-\proj) u_1 \|_{L^2({\cO})}.
	\end{align}
	We mention that all the computations we have done so far in Steps 1-4 would work if we replace  $\mathcal{L}_{+}$ by $\mathcal{L}_{-}$. Combining \eqref{a4-t1} and \eqref{a4-t3} we obtain \begin{align}
	& \left\| \frac{\sin(t \sqrt{\nA})}{\sqrt{\nA}} u_1 \right\|_{L^p(0,2\pi;\Dom(\nA_q^{\frac{1-r}{2}}))} \lesssim \left\| t \proj u_1 \right\|_{L^p(0,2\pi;\Dom(\nA_q^{\frac{1-r}{2}})} \nonumber\\
	& \hspace{ 6 cm} +  \left\| \mathcal{L}_{+}(\cdot) \bigl(\sqrt{\nA}\bigr)^{-1}((\mathds{1}-\proj) u_1) \right\|_{L^p(0,2\pi;\Dom(\nA_q^{\frac{1-r}{2}})} \nonumber\\
	& \lesssim \left\| t \proj u_1 \right\|_{L^p(0,2\pi;\Dom(\nA_q^{\frac{1-r}{2}}))} + \| (\mathds{1}-\proj) u_1 \|_{L^2({\cO})} \lesssim \|\proj u_1 \|_{\Dom(\nA_q^{\frac{1-r}{2}})} + \left\| u_1 \right\|_{L^2({\cO})} \lesssim \left\| u_1 \right\|_{L^2({\cO})}. \nonumber
	\end{align} This finishes the proof of Step 4 and, in particular, the first case.\\
	
	\textbf{Second case: when $L^1(0,2\pi;L^2({\cO})) \ni F \neq 0$:}  Due to the Duhamel formula \begin{align}
	u(t) = \cos(t \sqrt{\nA}) u_0 + \frac{\sin(t \sqrt{\nA})}{\sqrt{\nA}}u_1 + \int_{0}^{t} \frac{\sin((t-s) \sqrt{\nA})}{\sqrt{\nA}}F(s)\, ds. \nonumber
	\end{align}
	Applying the first case and using the calculation of \eqref{a2-t1} and \eqref{a4-t2} we get \begin{align}
	\| u \|_{L^p(0,2\pi;\Dom(\nA_q^{\frac{1-r}{2}}))} &  \lesssim \| u_0 \|_{\Dom(\nA^{\frac12 })} + \| u_1\|_{L^2({\cO})} \nonumber\\
	& \quad+ \int_{0}^{2\pi} \left\| \frac{\sin((t-s) \sqrt{\nA})}{\sqrt{\nA}}F(s) \right\|_{L^p(0,2\pi;\Dom(\nA_q^{\frac{1-r}{2}}))} \, ds \nonumber\\			
	& \lesssim \| u_0 \|_{\Dom(\nA^{\frac12 })} + \| u_1\|_{L^2({\cO})} + \int_{0}^{2\pi} \| F(s) \|_{L^2({\cO})} \, ds. \nonumber
	\end{align} Hence we have proved the Theorem \ref{thm-InHomStrichartz}.
\end{proof}

\section{Stochastic Strichartz estimates}\label{sec-Strichartz stochastic}
This section is devoted to prove a stochastic Strichartz inequality, which is sufficient to apply the Banach Fixed Point Theorem in the proof of a local well-posedness result for problem \eqref{AbsWe}, see Theorem \ref{thm-local existence} in Section \ref{sec-main}.

\subsection{Main assumptions}\label{subsec-mainAssump}
Here we describe the main assumptions  we consider in Sections  \ref{sec-Strichartz stochastic} -- \ref{sec-main}. Let us set
\begin{equation}\label{eqn-spaces}
H = L^2({\cO}), \;\; \rV = \Dom(\nA^{\frac12 });\;\;\; E = \Dom(\nA_q^{\frac{1-r}{2}}), \,
\end{equation}
where
$q \in [2,\infty)$,  $r \in [0,1]$ and $p=p(q,r)$ satisfies the equality \eqref{InHomStrichartzIneqCondn}.
Let us define the following Banach space.
For fixed $T>0$,  we put
\begin{align}
\label{eqn-X_T}
\Xa_T &:= \mathrm{C}\bigl([0,T];\rV\bigr),
\\
\label{eqn-Y_T}
\Ya_T &:= L^p(0,T;\rE),
\\
\label{eqn-Z_T}
\Za_T &:=  \mathrm{C}\bigl([0,T];\rV\bigr) \cap L^p(0,T;\rE),
\\
\label{eqn-Q_T}
Q_T&:= \mathrm{C}([0,T];H)
\end{align}

Obviously,  all these four spaces  are (separable) Banach spaces with naturally defined  norms, i.e.
\begin{align}
\label{eqn-X_T-norm}
\Vert \xi \Vert_{\Xa_T}&: = \sup_{t \in [0,T]}\Vert u(t) \Vert_{\rV},
\\
\label{eqn-Y_T-norm}
\Vert \xi \Vert_{\Ya_T}&: = \left(\int_0^T \vert \xi(s) \vert_{\rE}^p \, ds \right)^{1/p},
\\
\label{eqn-Z_T-norm}
\Vert \xi \Vert_{\Za_T} &:=  \left(\sup_{t \in [0,T]}\Vert u(t) \Vert_{\rV}^p+ \int_0^T \vert \xi(s) \vert_{\rE}^p \, ds \right)^{1/p},
\\
\label{eqn-Q_T-norm}
\Vert \xi \Vert_{Q_T} &:= \sup_{t \in [0,T]}\Vert u(t) \Vert_{H}.
\end{align}

By  $\mathbb{M}^{p}(\Za_T)$  we denote the Banach space of (equivalence classes) of all  $E$-valued progressively measurable processes $\{ u(t) ,t \in [0,T]\}$ having a continuous $\rV$-valued modification and  satisfying
\begin{equation}
\label{eqn-M^p(Z_T)}
\Vert \xi \Vert_{\mathbb{M}^{p}(\Za_T)}^p := \tE \bigl[ \Vert \xi \Vert_{\Za_T}^p \bigr] = \tE \bigl[ \Vert \xi \Vert_{\mathrm{C}([0,T];\rV)}^p + \Vert \xi \Vert_{L^p(0,T;\rE)}^p \bigr] < \infty.
\end{equation}

\subsection{Martingales}\label{subsec-martingales}
In order to define the It{\^o} type integrals for a Banach space valued stochastic process, we restrict ourself to, the so called, $M$-type 2 Banach spaces which are defined as follows.
\begin{definition}
A Banach space $E$ is of M-type $2$ iff there exists a constant $L:=L_2(E) >0$ such that for any $E$-valued martingale $\{ M(n), \;n \in \mathbb{N} \}$ the following holds:
\begin{equation*}
\sup_{n}\tE(\Vert M(n) \Vert_E^2) \leq L \sum_{n=0}^{\infty} \tE \bigl(\Vert M(n) - M(n-1) \Vert_E^2 \bigr),
\end{equation*}
where, as usual, we set $M(-1) = 0$.	
\end{definition}

For an interval $I \subseteq \mathbb{R}$, we say that an $E$-valued process $\{M(t), t \in I\}$ is an $E$-valued martingale iff $M(t) \in L^1(\Omega,\mathscr{F},\tP;E)$ for $t \in I$ and \begin{equation*}
	\tE\bigl(M(t) \vert \mathscr{F}_s \bigr) = M(s), \quad \tP\textrm{-a.s., for all } s \leq t \in I.
\end{equation*}

\subsection{Burkholder inequality}\label{subsec-Burkholder}

To prove the main result of this section we need the following consequence of the Kahane-Khintchin inequality and the It{\^o}-Nisio Theorem, see \cite{HVVW2016B}.  For any $\Lambda \in \gamma(K,E)$, by the It{\^o}-Nisio Theorem, the series $\sum_{j=1}^{\infty}\beta_j \Lambda(f_j)$ is $\tP$-a.s. convergent in $E$, where $\{f_j \}_{j \in \mathbb{N}}$ and $\{\beta_j \}_{j \in \mathbb{N}}$ are as in Definition \ref{def-gammaRadOper}.  Then, as an application of the Kahane-Khintchin inequality, X. Fernique \cite{Fernique_1970} proved that, for any $p \in [1,\infty)$,  there exists a positive constant $C(p,E)$ such that,
\begin{equation}\label{L2-LpConvergenceRel}
(C(p,E))^{-1} ~ \Vert \Lambda \Vert_{\gamma(K,E)} \leq  \left(\tE \bigg\Vert \sum\limits_{j \in \mathbb{N}} \beta_j \Lambda(f_j) \bigg\Vert_{E}^p  \right)^{\frac{1}{p}}   \leq C(p,E)~ \Vert \Lambda \Vert_{\gamma(K,E)}.
\end{equation}
This inequality tells that the convergence in $\LebAuxProbSpace$ can be replaced by a condition of convergence in $\LpAuxProbSpace$ for some (or any) $p \in [1,\infty)$. Furthermore, we need the following version of the Burkholder inequality which holds in our setting, refer \cite{O2004Diss} for the proof.
\begin{theorem}[Burkholder inequality]\label{thm-Burkholder inequality}
	Let $E$ be a $M$-type 2 Banach space and $p \in [1,\infty)$. If   $\xi \in \mathcal{M}_{loc}^2(\mathbb{R}_+, \gamma(K,E))$, then the following conditions   hold.
\begin{trivlist}
\item[(i)] There exists an  $E$-valued continuous process $\{\Xi(t)$, $t\in [0,T]\}$, such that
	\begin{equation}\label{eqn-Ito=rv}
\Xi(t)= \int_{0}^{t} \xi(s) \, dW(s) \;\tP \mbox{-a.s., for every }t\in [0,T].
	\end{equation}
The random variable $\Xi(t)$ will be denoted, unless a danger of ambiguity, by $\int_{0}^{t} \xi(s) \, dW(s)$.
\item[(ii)]  There exists a constant $B_p(E) > 0$, independent of $\xi$, such that for each  stopping time $\tau > 0$,
	\begin{equation}\label{BurkholderIneq}
	\tE \left[ \sup_{t \geq 0}  \bigg\Vert \int_{0}^{t \wedge \tau} \xi(s) \, dW(s) \bigg\Vert_E^p\right]\leq B_p(E) ~\tE \left[\int_{0}^{\tau} \Vert \xi(t)\Vert_{\gamma(K,E)}^2 \, dt \right]^{\frac{p}{2}}.
	\end{equation}
\end{trivlist}
Finally, \eqref{BurkholderIneq} holds  provided only that $\xi=\{\xi(t), 0\leq t<\tau\}$ is a local progressively measurable process.

\end{theorem}
\begin{remark}\label{remar-ito integral}[\textbf{Warning!}]
The It\^o integral $\int_{0}^{t}  \xi(s) \, dW(s)$ is, by definition, an element of $L^p(\Omega,\mathscr{F}_t;E)$, thus an equivalence class of
a certain class of $\mathscr{F}_t$/$\mathscr{B}(E)$-measurable functions. In what follows we will use a formally imprecise  formulation as
$\Xi(t)= \int_{0}^{t} \xi(s) \, dW(s)$ $\tP$-a.s. used in \eqref{eqn-Ito=rv},  instead of the correct, but awkward, one that $\Xi(t)$ belongs to  $\int_{0}^{t}  \xi(s) \, dW(s)$.
\end{remark}
We ask the reader to refer \cite[Corollary 3.7]{Brz+Millet_2014} for a proof of the following result. It is important to mention here that the range of $p$ assumed in the statement of \cite[Corollary 3.7]{Brz+Millet_2014}  is incorrect as from the proof in \cite{Brz+Millet_2014}  it is clear that the result only holds true if $p \geq 2$.
\begin{corollary}\label{cor-BurkholderIneqCorr}
	Let $E$ be a M-type 2 Banach space and $p \in [2,\infty)$. Then there exists a constant $\hat{B}_p(E)$, depending on $E$, such that for every $T \in (0,\infty]$ and every $L^p(0,T; E)$-valued progressively measurable process $\{\zeta(s), s \in [0,T)\}$,
	\begin{equation}\label{eqn-BurkholderIneqCorr}
	\tE \left[   \bigg\Vert \int_{0}^{T} \zeta(s) \, dW(s) \bigg\Vert_{L^p(0,T;E)}^p \right] \leq \hat{B}_p(E) ~ \tE \left[\int_{0}^{T} \Vert \zeta(s) \Vert_{\gamma(K,L^p(0,T;E))}^2 \, ds \right]^{\frac{p}{2}}.
	\end{equation}
\end{corollary}
For a $\gamma(K,H)$-valued random variable  $\xi$, let us define a $\gamma(K, L^p(0,T; E))$-valued process $\Xi=\{\Xi_r : r \in [0,T] \}$  by
 \begin{equation}\label{eqn-Xi}
\Xi_r:= \left\{  [0,T] \ni t \mapsto \mathds{1}_{[r,T]}(t) \frac{\sin((t-r) \sqrt{\nA})}{\sqrt{\nA}} \xi(r)   \right\} \in \gamma(K,L^p(0,T;E)),\;\; r \in [0,T].
\end{equation}
We will need  the following auxiliary result.

\begin{lemma}\label{lem-AuxRes1}
	Assume that  $T>0$. Then there exists a positive constant $C(p,T,H)$ such that for every     $\gamma(K,H)$-valued progressively measurable process $\xi$,  the $\gamma(K, L^p(0,T; E))$-valued process  $\{\Xi_r : r \in [0,T] \}$,  defined by formula \eqref{eqn-Xi},
	is progressively measurable and,
	\begin{equation}\label{AuxRes1-eqn1}
	\Vert \Xi_r \Vert_{\gamma(K,L^p(0,T; E))} \leq C(p,T,H)~ \Vert \xi(r) \Vert_{\gamma(K,H)}, \;\; \mbox{for each $r \in [0,T]$}.
	\end{equation}
\end{lemma}
\begin{proof}[\textbf{Proof of Lemma \ref{lem-AuxRes1}}]
	Let us consider a sequence $\{ \beta_j \}_{j \in \mathbb{N}}$ of i.i.d. $N(0,1)$ random variables on probability space $\AuxProbSpace$, and an orthonormal basis $\{ f_j \}_{j \in \mathbb{N}}$  of the separable Hilbert space $K$. We first observe that the random variable $\Xi_r$   is well-defined because by Theorem \ref{thm-InHomStrichartz}, for each $r \in [0,T]$ and $x \in H$, the solution of the following homogeneous wave equation
	\begin{equation*}
	\left\{
	\begin{array}{ll}
	u_{tt} - \Delta u = 0 \mbox{ on } [r,r+T]  \\
	u(r) =0, \quad u_t(r) = x,
	\end{array}
	\right.
	\end{equation*}	
	belongs to $L^p(r,r+T;E)$.
	In particular,
	\begin{equation*}
	\mathds{1}_{[r,T]}(\cdot ) \frac{\sin((\cdot-r) \sqrt{\nA})}{\sqrt{\nA}} x \in L^p(0,T;E),
	\end{equation*}
	and the map
	\begin{equation}\label{eqn-Lambda_r}
	\Lambda_r:  H \ni x \mapsto \mathds{1}_{[r,T]}(\cdot ) \frac{\sin((\cdot-r) \sqrt{\nA})}{\sqrt{\nA}} x \in L^p(0,T;E),
	\end{equation}
	is linear and continuous.  Moreover, we have $\sup_{r \in [0,T]} \Vert \Lambda_r \Vert<\infty$.
	By the above argument and \eqref{eqn-Xi}, we infer  that
	\begin{equation*}
		\Xi_r(\omega)= \Lambda_r  \circ [\xi(r,\omega)], \;\;\; (r,\omega) \in [0,T]\times \Omega.
	\end{equation*}
	Next, for each $r\in [0,T]$ and $x \in H$, the deterministic Strichartz estimate \eqref{InHomStrichartzIneq} yields
	\begin{align}\label{lem-AuxRes1-1}
	 \| \Lambda_r x \|_{L^p(0,T;E)}  = \left( \bigintsss_{r}^{T} \left\|  \frac{\sin((t-r) \sqrt{\nA})}{\sqrt{\nA}}x  \right\|_{E}^p \, dt   \right)^{\frac{1}{p}} \leq C_T ~\| x\|_{H},
	\end{align}	
	with the RHS being independent of $r$.  Consequently, we have
	\begin{align}\label{lem-AuxRes1-2}
	\|\Lambda_r \|_{\mathcal{L}(H,L^p(0,T;E))} &  = \sup_{ \substack{x \in H \\ \| x\|_H \leq 1}  } \| \Lambda_r x \|_{L^p(0,T;E)}  \leq C_T.
	\end{align}
	To move further with the proof let us set the following useful notation:
	\begin{equation}\label{lem-AuxRes1-3}
	f_r := \mathds{1}_{[r,T]}, r \in [0,T] \qquad \textrm{ and } \qquad g(t,x) := \frac{\sin(t \sqrt{\nA})}{\sqrt{\nA}}x, (t,x) \in [0,T] \times H.
	\end{equation}
	We have the following result.
	\begin{lemma}\label{lem-Lambda_cont_wrt_r}
		If $x \in H$, then the function $$ [0,T] \ni r \mapsto \Lambda_r x \in L^p(0,T;E),$$
		is continuous.
	\end{lemma}
	\begin{proof}[Proof of Lemma \ref{lem-Lambda_cont_wrt_r}]
		
	We first assume that  $x \in \Dom(A^{k-\frac{1}{2}})$, where 	 $k \in \mathbb{N}$ is such that $\Dom(A^k) \hookrightarrow E$.  Since $\Dom(A^{k-\frac{1}{2}})$  is dense in $H$, it is sufficient to prove the lemma for $x \in \Dom(A^{k-\frac{1}{2}})$.  Note that, by property \eqref{ineq-L^infty-regular} for $x \in \Dom(A^{k-\frac{1}{2}})$ we have
		\begin{equation*}
		g(\cdot,x) \in \mathrm{C}([0,T]; \Dom(A^k)) \subset \mathrm{C}([0,T]; E).
		\end{equation*}
		Since $\Vert f_r - f_{r_0} \Vert_{L^p(0,T;\mathbb{R})}^p = r-r_0$, we infer that the function
		\begin{equation*}
		[0,T] \ni r \mapsto f_r \in L^p(0,T;\mathbb{R}),
		\end{equation*}
		is continuous.
		Next we claim  that if $g(\cdot,x) \in \mathrm{C}([0,T];E)$, then the map
		\begin{equation*}
		h: [0,T] \ni r \mapsto f_r g(\cdot-r,x) =: h_r \in L^p(0,T;E),
		\end{equation*}
		is continuous. For this we note that  for $0 \leq r < s \leq T$ we have
		\begin{equation}\label{claim2-0}
		\Vert h_r-h_s  \Vert _{L^p(0,T;E)}^p = \int_{r}^{s} \|h_r(t) \|_E^p \, dt + \int_{s}^{T} \|h_r(t) - h_s(t) \|_E^p \, dt,
		\end{equation}
because
		\begin{equation*}
			(h_r-h_s)(t) = \threepartdef {0}{t \in [0,r]}   {h_r(t)}{t \in [r,s]}    {h_r(t) - h_s(t)}{t \geq s}.
		\end{equation*}
		Concerning  the first integral on the RHS of \eqref{claim2-0} we have
		\begin{equation*}
			 \int_{r}^{s} \|h_r(t) \|_E^p \, dt  =  \int_{0}^{s-r} \|g(\tau,x) \|_E^p \, d \tau=\int_{0}^{\delta} \|g(\tau,x) \|_E^p \, d \tau \leq  (s-r) \Vert g(\cdot,x) \Vert_{C([0,T];E)}^p.
		\end{equation*}
		For the second integral in the RHS of \eqref{claim2-0}, we have
		\begin{align*}
			\int_{s}^{T} \|h_r(t) - h_s(t) \|_E^p \, dt = \int_{0}^{T-s} \|g(s-r+\tau,x) - g(\tau,x) \|_E^p \, dt \leq T \Delta(s-r),
		\end{align*}
		where, by the  uniform continuity of function $g(\cdot,x): [0,T] \to E$,

		$$ \Delta(\delta):= \sup_{\substack{s,t\in [0,T] \\ |t-s| \leq \delta }  } \|g(t,x) - g(s,x)\|_E \to 0 \mbox{ as } \delta \to 0. $$
		Hence the continuity of function $h$ follows and we are done with the proof of Lemma~ \ref{lem-Lambda_cont_wrt_r}.
	\end{proof}
	Thus, by coupling the Lemma \ref{lem-Lambda_cont_wrt_r} with Proposition \ref{prop-Neidhard Lemma 44} and \cite[Corollary 1.1.29]{HVVW2016B}, we infer that the process $\Xi$  is progressively measurable.

	It only remains to prove the inequality \eqref{AuxRes1-eqn1}. For this aim let us fix  $r \in [0,T]$. By invoking the inhomogeneous Strichartz estimate from Theorem  \ref{thm-InHomStrichartz} and  \eqref{L2-LpConvergenceRel}, followed by  \eqref{lem-AuxRes1-2}, we obtain
	\begin{align}
	& \|	\Lambda_r \circ \xi \|_{\gamma(K,L^p(0,T;E))}  \simeq \bigl( \tE \left[\bigg\Vert \sum\limits_{j \in \mathbb{N}} \beta_j ~ \Lambda_r(\xi(e_j)) \bigg\Vert_{L^p(0,T;E)}^p \right] \bigr)^{\frac{1}{p}} \nonumber\\
	& = \bigl( \tE \left[  \bigg\Vert \Lambda_r \bigl( \sum\limits_{j \in \mathbb{N}} \beta_j ~ \xi(e_j) \bigr) \bigg\Vert_{L^p(0,T;E)}^p \right] \bigr)^{\frac{1}{p}} \leq \| \Lambda_r\|_{\mathcal{L}(H,L^p(0,T;E))} \bigl( \tE \left[  \bigg\Vert  \sum\limits_{j \in \mathbb{N}} \beta_j ~ \xi(e_j)  \bigg\Vert_{H}^p \right] \bigr)^{\frac{1}{p}} \nonumber\\
	&  \leq C(p,H) ~C_T ~ \| \xi\|_{\gamma(K,H)}. \nonumber
	\end{align}
	Hence the proof of Lemma \ref{lem-AuxRes1} is complete.
\end{proof}

The following main result of this section is one of the most important ingredient in the proof of the local existence theorem in Section \ref{sec-main}.

\begin{theorem}[Stochastic Strichartz Estimates]\label{thm-StocStrichartzEstimate}
	Let us assume that $T>0$ and $p \in [2,\infty)$. Then there exist  constants\footnote[1]{The constant $K$ depends on $T$ only in the case of Neumann boundary conditions.} $\tilde{K}_T := K(p,T,H)>0$ and  $\tilde{C}_T := \tilde{C}(p,T,E,H)>0$ such that if  a process $\xi$  belongs to
$\mathcal{M}^{p}(\mathcal{L}^{2}([0,T],\gamma(K,H))$, then the following assertions hold.
\begin{itemize}
\item[(I)]  There exists an $\mathbb{F}$-adapted,  and $\rV\times H$-valued\footnote[2]{Let us recall that $\rV={\Dom(\nA^{\frac12 })}$. } continuous   process  $(\tilde{u},\tilde{v})$
such that
\begin{equation}\label{StocStrichartzEstimate1}
	\tE \left[\sup_{t \in [0,T]}  \Vert \tilde  u(t) \Vert_{\rV}^p   + \sup_{t \in [0,T]}  \Vert \tilde  v(t) \Vert_{H}^p \right] \leq \tilde{K}_T~ \tE \left[\int_{0}^{T}\Vert \xi(t)\Vert_{\gamma(K,H)}^2 \, dt\right]^{\frac{p}{2}},
\end{equation}
\begin{equation}\label{eqn-modification}
\tilde{u}(t)=\int_{0}^{t} \frac{\sin((t-s) \sqrt{\nA})}{\sqrt{\nA}} \xi(s) ~ dW(s), \;\; \tP\mbox{-a.s.}, \;\; \mbox{for every } t \in  [0,T],
\end{equation}
\begin{equation}\label{eqn-u_t=v}
\tilde{u}(t)=\tilde{u}(0)+\int_{0}^{t} \tilde{v}(s) ~ ds \mbox{ in }H, \;\; \tP\mbox{-a.s.}, \;\; \mbox{for every } t \in  [0,T].
\end{equation}

\item[(II)] There exists an $E$-valued  progressively measurable process $\tilde{\tilde{u}}$ such that
\begin{equation}\label{StocStrichartzEstimate2}
	\tE \left[ \int_{0}^{T} \Vert \tilde{\tilde{u}}(t) \Vert_{E}^p \, dt \right] \leq \tilde{C}_T~ \tE \left[\int_{0}^{T}  ~ \Vert \xi(t) \Vert_{\gamma(K,H)}^2 \, dt \right]^{\frac{p}{2}},
\end{equation}
\begin{equation}\label{eqn-equality almost sure}
j (\tilde{\tilde{u}}(t,\omega))  = i (\tilde{u}(t,\omega)) \;\; \mbox{ for $\Leb\otimes \tP$-almost all } (t,\omega)\in [0,T]\times \Omega,
\end{equation} where $i: \rV \hookrightarrow H$ and $j: E \hookrightarrow H$ are the natural embeddings. \\%
\item[(III)] Moreover, if processes $\xi_1$ and $\xi_2$ are equivalent, then so are the corresponding processes $\tilde{\tilde{u_1}}$ and $\tilde{\tilde{u_2}}$.
In particular, the map
$$\mathcal{M}^{p}(\mathcal{L}^{2}([0,T],\gamma(K,H)) \ni \xi \mapsto  \tilde{\tilde{u}} \in \mathcal{L}^p (0,T;E), $$
extends in a unique way to the following bounded and linear map
\begin{equation}\label{JxiDefn}
\mathbb{M}^{p}({L}^{2}([0,T],\gamma(K,H)) \ni [\xi] \mapsto  [\tilde{\tilde{u}}] \in L^p (0,T;E).
\end{equation}

\end{itemize}
\end{theorem}

\begin{remark}\label{rem-stopping time versions of Burkholder inequality}
Suppose that $\tau$ is a stopping time such that $\tau \leq T$ for some $T>0$ and $\xi(s)$, $s\in [0,\tau)$ is an $\gamma(K,E)$-valued  progressively measurable process. Since the process $\mathds{1}_{[0,\tau)}(s)$, $s\in [0,\infty)$ is well-measurable, see \cite[Proposition 4.2]{Metivier_1982} and, see  \cite[Theorem 1.6]{Metivier_1982}, the $\sigma$-field of well-measurable sets is smaller than the $\sigma$-field of progressively measurable sets, it follows that the process $\mathds{1}_{[0,\tau)}(s)$, $s\in [0,\infty)$ is progressively measurable.
In particular, the process $\{\mathds{1}_{[0,\tau)}(s)\xi(s), s\in [0,\infty)\}$ is progressively measurable. Hence,
by applying inequalities \eqref{StocStrichartzEstimate1} and \eqref{StocStrichartzEstimate2} to the process $\mathds{1}_{[0,\tau)}\xi$ we infer the following stopped versions of those inequalities. \\
\begin{itemize}
\item[(i)]  There exists an $\mathbb{F}$-adapted,  $\rV\times H$-valued continuous   local process  $(\tilde{u},\tilde{v})=\bigl((\tilde{u}(t),\tilde{v}(t)), t\in [0,T]\bigr)$,
such that
\begin{equation}\label{eqn-StocStrichartzEstimate1}
	\tE \left[\sup_{t \in [0,T]}  \Vert \tilde  u(t\wedge \tau  ) \Vert_{\rV}^p  + \sup_{t \in [0,T]}  \Vert \tilde  v(t\wedge \tau  ) \Vert_{H}^p \right] \leq \tilde{K}_T~ \tE \left[\int_{0}^{\tau}\Vert \xi(t)\Vert_{\gamma(K,H)}^2 \, dt\right]^{\frac{p}{2}},
\end{equation}%
\begin{equation}\label{eqn-modification-stopping time}
\tilde{u}(t  )=\int_{0}^{t} \mathds{1}_{[0,t \wedge \tau) }(s) \frac{\sin((t-s) \sqrt{\nA})}{\sqrt{\nA}} \xi(s) ~ dW(s), \;\; \tP\mbox{-a.s.}, \;\; \mbox{for every } t \in  [0,T],
\end{equation}
\begin{equation}\label{eqn-u_t=v-stopping time}
\tilde{u}(t)=\tilde{u}(0)+\int_{0}^{t} \tilde{v}(s) ~ ds\mbox{ in }H, \;\; \tP\mbox{-a.s.}, \;\; \mbox{for every } t \in  [0,\tau).
\end{equation}
\item[(ii)] There exists an $E$-valued  progressively measurable local process $\tilde{\tilde{u}}(t)$, $t\in [0,\tau]$, such that
\begin{equation}\label{eqn-StocStrichartzEstimate2}
	\tE \left[ \int_{0}^{\tau } \Vert \tilde{\tilde{u}}(t) \Vert_{E}^p \, dt \right] \leq \tilde{C}_T~ \tE \left[\int_{0}^{\tau}  ~ \Vert \xi(t) \Vert_{\gamma(K,H)}^2 \, dt \right]^{\frac{p}{2}},
\end{equation}
\begin{equation}\label{eqn-equality almost sure1}
j (\tilde{\tilde{u}}(t,\omega))  = i (\tilde{u}(t,\omega)) \;\; \mbox{ for $\Leb\otimes \tP$-almost all } (t,\omega)\in [0,\tau]\times \Omega.
\end{equation}%
\end{itemize}
\end{remark}
\begin{remark}\label{rem-StocStrichartzEstimate} It follows from the proof that
$\tilde{K}_T\leq M_1  e^{mT}  B_{p}(\mathcal{H})$ for some constants $m \geq 0$ and $M_1 \geq 1$ and $\tilde{C}_T \leq  C_T ~C(p,H)~ \hat{B}_p(E)$.
\end{remark}

\begin{proof}[\textbf{Proof of part (I) of Theorem \ref{thm-StocStrichartzEstimate}}] In what follows we fix either the Dirichlet or the Neumann boundary conditions. Let us fix $p \in [2,\infty)$.
To prove the first assertion, let us consider the Hilbert space $\mathcal{H}:= \rV \times H=\Dom(\nA^{\frac12 }) \times H$ and a linear  operator $\mathcal{A}$ on $\mathcal{H}$ defined by
\begin{equation} \label{eqn-group generator}
\begin{split}%{align}
 \Dom(\mathcal{A}) &= \Dom(\nA)\times \Dom(\nA^{\frac12 }),\\
 \mathcal{A}(x,y)&=(y, -Ax), \;\; (x,y) \in  \Dom(\mathcal{A}).
 \end{split}%{align}
\end{equation}
It is well known that  	since $\nA$ is non-negative and self-adjoint in $L^2({\cO})$, one may prove that $\mathcal{A}$ generates a $C_0$-group on $\mathcal{H}$, denoted by $\{S(t)\}_{t \geq 0}$. Moreover, for $t \geq 0$,
	\begin{equation}\label{eqn-group S(t)}
	S(t)(x,y) =
	\bigl( \cos(t\sqrt{\nA})x+ \sin(t \sqrt{\nA})/\sqrt{\nA}y,
	- \sqrt{\nA} \sin(t \sqrt{\nA})x +  \cos(t\sqrt{\nA}) y\bigr), (x,y) \in \mathcal{H}.
	\end{equation}
i.e., using a matrix notation,

	\begin{equation*}
	S(t) = \begin{pmatrix}
	\cos(t\sqrt{\nA}) & \sin(t \sqrt{\nA})/\sqrt{\nA} \\
	- \sqrt{\nA} \sin(t \sqrt{\nA}) &  \cos(t\sqrt{\nA})
	\end{pmatrix},
	\end{equation*}

	It follows, with being $\pi_1:\mathcal{H} \to \rV$  the natural projection,  that for $ t\geq 0$,
	\begin{equation}\label{eqn-pi_1S(t)}
\pi_1 \bigl(S(t)(0,y)\bigr)= \sin(t \sqrt{\nA})/\sqrt{\nA} y, \;\; y\in H.
\end{equation}
Let us  introduce the following auxiliary  $\gamma(K,\mathcal{H})$-valued  process $\tilde{\xi}$
	\begin{align}  \tilde{\xi}(t)[k] = \bigl(
	0,
	\xi(t)[k] \bigr), \;\; k \in K,\;\; t \in [0,T] .
	\end{align}
 Our argument now is based on  \cite{Hausenblas+Seidler_2001}.
  We begin by observing, see e.g. \cite[Theorem~12.2]{Wentzel_1981},  that there exists an $\mathcal{H}$-valued continuous  process
  $\tilde{\mathfrak{q}}= \{\tilde{\mathfrak{q}}(t),t\in [0,T]\}$   such that
   \[\tilde{\mathfrak{q}}(t)=\int_{0}^{t} S(-s) ~ \tilde{\xi}(s)~ dW(s),\;\; \tP\mbox{-a.s., for every } t\in [0,T].
  \]
  Since $\{S(t)\}_{t \geq 0}$  is a $C_0$-group, we infer that  the process $\tilde{\mathfrak{u}}$
 defined  by
	\begin{align}
	\tilde{\mathfrak{u}}(t)  &:=S(t)\tilde{\mathfrak{q}}(t) , \;\;\; t\in [0,T], \nonumber
	\end{align}
is a continuous $\mathcal{H}$-valued and 	that
\begin{equation}\label{eqn-u tilde}
		\tilde{\mathfrak{u}}(t) =\int_{0}^{t} S(t-s) ~ \tilde{\xi}(s)~ dW(s), \;\; \tP\mbox{-a.s., for every } t \in [0,T].
	\end{equation}
 With $\pi_2:\mathcal{H} \to H$ being the natural projection, we define  continuous $\rV$ and $H$-valued   processes $\tilde{u}$ and $\tilde{v}$, respectively,  by
	\begin{equation}\label{eqn-u tilde-2}
		\tilde{u}(t):=\pi_1(\tilde{\mathfrak{u}}(t)), \qquad \textrm{ and } \qquad  \tilde{v}(t):=\pi_2(\tilde{\mathfrak{u}}(t)),  \;\;\; t \in [0,T].
	\end{equation}
	 Applying  identity \eqref{eqn-pi_1S(t)} to  the two previous equalities \eqref{eqn-u tilde}-\eqref{eqn-u tilde-2} we infer that
\[
\tilde{u}(t)=\int_{0}^{t} \frac{\sin((t-s) \sqrt{\nA})}{\sqrt{\nA}} \xi(s) ~ dW(s), \;\;\;\tP\mbox{ -a.s. for all } t \in [0,T].
\]
Thus, we proved that   $\tilde{u}$ is  ${H}_A$-valued continuous and $\mathbb{F}$-adapted process satisfying equality \eqref{eqn-modification}.
	Moreover, using the Burkholder inequality \eqref{BurkholderIneq} and the bound property of $C_0$-group,   we get the following train of inequalities:
	\begin{align}
	& \tE \left[\sup_{t \in [0,T]}  \Vert \tilde u(t)\Vert_{\Dom(\nA^{\frac12 })}^p \right]  \leq \tE \left[\sup_{t \in [0,T]}  \Vert \tilde{\mathfrak{u}}(t)\Vert_{\mathcal{H}}^p \right]  = \tE\left[\sup_{t \in [0,T]} \bigg\Vert S(t)  \int_{0}^{t} S(-s) \tilde{\xi}(s)~ dW(s) \bigg\Vert_{\mathcal{H}}^p \right]\nonumber\\
	&\quad  \leq K_T' B_p(\mathcal{H}) ~ \tE\left[ \int_{0}^{T} \Vert  S(-s) \tilde{\xi}(s) \Vert_{\gamma(K,\mathcal{H})}^2 ds  \right]^{\frac{p}{2}} \nonumber\\
	& \quad \leq K_T B_p(\mathcal{H}) ~ \tE\left[ \int_{0}^{T} \Vert  \tilde{\xi}(s) \Vert_{\gamma(K,\mathcal{H})}^2 ds  \right]^{\frac{p}{2}}
	= K_T B_p(\mathcal{H}) ~ \tE\left[ \int_{0}^{T} \Vert  \xi(s) \Vert_{\gamma(K,H)}^2 ds  \right]^{\frac{p}{2}}, \nonumber
	\end{align}
	where $K_T \leq M_1 e^{mT}$  for some constants $m \geq 0$ and $M_1 \geq 1$. This  yields  inequality \eqref {StocStrichartzEstimate1} and, in particular,  assertion of part  (i).\\
Equality \eqref{eqn-u_t=v-stopping time} follows from Proposition \ref{prop-equivalence} and equalities \eqref{eqn-u tilde-2}.
\end{proof}

\begin{proof}[\textbf{Proof of part (II) of Theorem \ref{thm-StocStrichartzEstimate}}]

We split the proof  into two steps. First we will prove the theorem for  more regular processes.  Then we will transfer the results to the right class of processes  by employing a suitable approximation.

\textbf{Step 1:}
We begin by observing that by the classical Sobolev embedding theorem there exists
  natural number $k$ such  that the Hilbert space
\begin{equation}\label{eqn-embed}
\mbox{ $\Dom(\nA^{k+\frac12})$ is continuously embedded into the Banach space $E=\Dom(\nA_q^{\frac{1-r}{2}})$}.
\end{equation}
Let us fix $p \in [2,\infty)$. Let us assume that a process  $\xi$ belongs to  $\mathcal{M}^{p}(\mathcal{L}^{2}([0,T],\gamma(K,\Dom(A^k)))$. By assertion (i), we infer that there exists an $\mathbb{F}$-adapted  ${\Dom(\nA^{k+\frac12})}$-valued continuous    process $\tilde{u}$ which satisfies condition \eqref{eqn-modification} and the following inequality
\begin{equation*}
	\tE \Bigl[~\sup_{t \in [0,T]}  \Vert \tilde  u(t) \Vert_{\Dom(\nA^{k+\frac12})}^p \Bigr] \leq K(p,T,H)~ \tE \left[\int_{0}^{T}\Vert \xi(t)\Vert_{\gamma(K,\Dom(A^k))}^2 \, dt\right]^{\frac{p}{2}}<\infty.
\end{equation*}

Also, let us note that, in view  of our additional assumption  \eqref{eqn-embed},  the process $\tilde{u}$ is an $\mathbb{F}$-adapted $\rV$-valued continuous  (hence progressively measurable), and
\begin{equation*}
	\tE \bigl[ \Vert \tilde u \Vert_{\mathrm{C}([0,T];\rV)}^p \bigr]< \infty.
\end{equation*}
Next, we define an $\gamma(K, L^p(0,T; \Dom(A^k)))$-valued process $\Xi=\{\Xi_r : r \in [0,T] \}$ by formula \eqref{eqn-Xi}.
This process, in view of Lemma \ref{lem-AuxRes1}, is  progressively measurable and, by   the Burkholder inequality \eqref{eqn-BurkholderIneqCorr} together with inequality \eqref{AuxRes1-eqn1},  it satisfies the following inequality
  \begin{align}\label{ineq-Burkholder-02}
 \tE \left[ \left\Vert  \int_{0}^{T} \Xi_r \, dW(r) \right\Vert_{L^p(0,T;E)}^p \right]
	& \leq C_T ~ C(p,H)~ \hat{B}_p(E) ~ \tE \left[ \int_{0}^{T}  \Vert \xi(r) \Vert_{\gamma(K,H)}^2 \, dr \right]^{\frac{p}{2}}.
\end{align}
Let us choose  an $\mathscr{F}_T$/$\mathscr{B}(L^p(0,T;E) )$-measurable function  $\tilde{\tilde{\mathbf{u}}}: \Omega \to L^p(0,T;E) $  such that, see Remark \ref{remar-ito integral},

 \begin{equation}\label{eqn-test 00}
\tilde{\tilde{\mathbf{u}}}=\int_{0}^{T} \Xi_r \, dW(r) \mbox{ in }L^p(0,T;E), \;\;\; \tP\mbox{-a.s.}.
\end{equation}
By the first part of Proposition \ref{prop-pointwise evaluation}, there exists
an $\mathscr{B}([0,T])\otimes \mathscr{F}_T$/$\mathscr{B}(E)$-measurable function
\[
\tilde{\tilde{u}}:[0,T] \times \Omega \to E,
\]
and a set  $\Omega^\prime \subset \Omega$ such that $\tP(\Omega^\prime)=1$ and for every $\omega \in \Omega^\prime$,
\[
\tilde{\tilde{u}}(\cdot,\omega)=\tilde{\tilde{\mathbf{u}}}(\omega) \mbox{ in } L^p(0,T;E).
\]
Later on we will  show that $\tilde{\tilde{u}}$ is progressively measurable process.
Then, by inequality \eqref{ineq-Burkholder-02} we infer that process $\tilde{\tilde{u}}$ satisfies  inequality \eqref{StocStrichartzEstimate2}.

Let us recall that $i: \rV \hookrightarrow H$ and $j: E \hookrightarrow H$ are the natural embeddings. We define corresponding Nemytski type  embeddings $I$ and $J$ by
 \begin{align}
 \label{eqn-I}
  I:\mathrm{C}([0,T];\rV) &\ni f \mapsto i \circ f \in  L^2(0,T;H),
 \\\label{eqn-J}
  J:L^p(0,T;E) &\ni f \mapsto j \circ f \in  L^2(0,T;H),
\end{align}
 and observe that both $I$ and  $J$ are continuous. Therefore we  deduce that
  $J \circ \tilde{\tilde{\mathbf{u}}}$ is an $L^2(0,T;H)$-valued random variable and
\begin{equation*}%\label{eqn-test 02}
 J \circ \tilde{\tilde{\mathbf{u}}} =\int_{0}^{T} (J \circ \Xi_r) \, dW(r) \mbox{ in } L^2(0,T;H),\;\;  \tP\mbox{-a.s.}.
\end{equation*}

 \medskip

 Since the process $\tilde{u}$ has $\tP$-almost surely  continuous $\rV$-valued trajectories, by \cite[Proposition 3.18]{DaPZ2014B} it induces, in a natural way, an $\mathscr{F}_T$/$\mathscr{B}(\mathrm{C}([0,T];\rV))$-measurable function $\tilde{\mathbf{u}}:\Omega \to \mathrm{C}([0,T];\rV)$. Because the map $I
 :\mathrm{C}([0,T];\rV) \to  L^2(0,T;H)$ is continuous,  $I (\tilde{\mathbf{u}}):\Omega \to L^2(0,T;H)$ is  $\mathscr{F}_T$/$\mathscr{B}(L^2(0,T;H))$-measurable.
We claim that
\begin{equation}\label{eqn-test 03}
I(\tilde{\mathbf{u}})= \int_0^T (J \circ \Xi_r)\,  dW(r) \mbox{ in } L^2(0,T;H),\;\;  \tP\mbox{-a.s.}.\end{equation}

\begin{proof}[Proof of equality \eqref{eqn-test 03}] Since $H$ is a separable Hilbert space, the Banach space $\mathrm{C}([0,T];H)$ is also separable.
Let us choose  a dense subset $\{h_n: n\in \mathbb{N}\}$ in $\mathrm{C}([0,T];H)$. Let us choose and fix $n\in \mathbb{N}$ and define a linear and bounded operator
\begin{equation}\label{eqn-H_n}
H_n: L^2(0,T;H)\ni x \mapsto \lb x,h_n\rb_{L^2(0,T;H)} \in \mathbb{R}.
\end{equation}
Thus, we infer that
\begin{equation}\label{eqn-test 03b}
H_n(\int_0^T (J \circ \Xi_r)\,  dW(r))= \lb \int_0^T (J \circ \Xi_r)\,  dW(r) ,h_n\rb_{L^2(0,T;H)} \mbox{ on } \Omega,
\end{equation}
and
\begin{equation*}%\label{eqn-test 04}
H_n(\int_0^T (J \circ \Xi_r)\,  dW(r))=\int_0^T (H_n \circ J \circ \Xi_r)\,  dW(r) \mbox{ in } \mathbb{R},\;\;  \tP\mbox{-a.s.}.
\end{equation*}
Let us note that by  definitions \eqref{eqn-H_n} of $H_n$, \eqref{eqn-J} of  $J$ and   \eqref{eqn-Xi} of $\Xi_r$ we have, for all $k\in K$ and  $r \in [0,T]$, the following equality
\begin{align}%\label{eqn-test 05}
& (H_n \circ J \circ \Xi_r)[k]=H_n  \bigl( (J \circ \Xi_r)[k]\bigr) =
\lb (J \circ \Xi_r),h_n\rb_{L^2(0,T;H)} \nonumber\\
&\quad =\int_0^T \lb j\bigl((\Xi_r [k] )(t)\bigr),h_n(t)\rb_{H}\, dt
 =\int_0^T \mathds{1}_{[r,T]}(t)  \lb \frac{\sin((t-r) \sqrt{\nA})}{\sqrt{\nA}} (\xi(r)[k]) ,h_n(t)\rb_{H}\, dt.
\nonumber \end{align}
Therefore, we infer that   $\tP$-almost surely
\begin{align}\nonumber
H_n(\int_0^T (J \circ \Xi_r)\,  dW(r))
&=
\int_0^T \int_0^T \mathds{1}_{[r,T]}(t)  \lb  \frac{\sin((t-r) \sqrt{\nA})}{\sqrt{\nA}} \xi(r) ,h_n(t)\rb_{H}\, dt \,  dW(r)
\\
&\lefteqn{\hbox{\hspace{-1.0truecm}}=
\int_0^T \int_0^T  \mathds{1}_{[0,t]}(r)   \lb  \frac{\sin((t-r) \sqrt{\nA})}{\sqrt{\nA}} \xi(r) ,h_n(t)\rb_{H}\, dW(r)   \,  dt,}
\label{eqn-test 06}%\nonumber
\end{align}
where the last equality is a consequence of the stochastic Fubini Theorem \cite{Brz+G+R+W_2021}, which is a generalization of  \cite[Theorem 2.4.16]{Carroll1999T} and \cite[Theorem 2.2]{Veraar_2012}.

On the other hand, by the definition \eqref{eqn-I} of the map $I$, we have
\begin{align}\label{eqn-test 07}
\lb I(\tilde{\mathbf{u}}),h_n\rb_{L^2(0,T;H)} &=\int_0^T \lb I(\tilde{\mathbf{u}})(t),h_n(t)\rb_{H}  \,dt
\\
=\int_0^T \lb i(\tilde{u}(t)),h_n(t)\rb_{H} &=\int_0^T \lb \int_{0}^{t} \frac{\sin((t-s) \sqrt{\nA})}{\sqrt{\nA}} \xi(s) ~ dW(s),h_n(t)\rb_{H} \,dt
\nonumber\\
&=\int_0^T \lb \int_{0}^{T} \mathds{1}_{[0,t]}(r)
 \frac{\sin((t-s) \sqrt{\nA})}{\sqrt{\nA}} \xi(s) ~ dW(s),h_n(t)\rb_{H} \,dt.
\nonumber
\end{align}

Thus, from \eqref{eqn-test 03b}, \eqref{eqn-test 07} and \eqref{eqn-test 06}, we infer that for every $n\in \mathbb{N}$,
\[
\lb I(\tilde{\mathbf{u}}),h_n\rb_{L^2(0,T;H)}= \lb \int_{0}^{T} (J \circ \Xi_r) \, dW(r) ,h_n\rb_{L^2(0,T;H)}, \;\;\tP\mbox{-a.s.}.
\]
By the density of the countable set $\{h_n\}$ in $L^2(0,T;H)  $ we deduce \eqref{eqn-test 03}.
\end{proof}

From the just proven equation \eqref{eqn-test 03} and equality \eqref{eqn-test 00}  we infer that

\[
I \circ \tilde{\mathbf{u}}=J \circ \tilde{\tilde{\mathbf{u}}}\mbox{ in } L^2(0,T;H),\;\;\tP\mbox{-a.s.}.
\]
Hence, by the second part of Proposition \ref{prop-pointwise evaluation}, we infer that  the $H$-valued processes   $i(\tilde{u}):[0,T] \times \Omega \to H$ and $j (\tilde{\tilde{u}}):[0,T] \times \Omega \to H$,
are $\Leb \otimes \tP$ equal. Since, the former is $H$-valued progressively measurable, by the Kuratowski Theorem, see e.g. \cite[Corollary~I.3.3]{Parthasarathy1967B} and the argument in the proof of \cite[Proposition A.1]{Brz+Haus+Mot_2013}  we infer that process
$\tilde{\tilde{u}}$ is $E$-valued progressively measurable.  This concludes the proof of \textbf{Step 1}.

\textbf{Step 2:} The result follows by applying  Step 1. Let $\xi$ be a progressively measurable process from the space $\mathcal{M}^{p}(\mathcal{L}^{2}([0,T],\gamma(K,H))$, and $k \in \mathbb{N}$ as in Step 1.

We choose a sequence $\{ \xi_n\}_{n \in \mathbb{N}}$ of processes from $\mathcal{M}^{p}(\mathcal{L}^{2}([0,T],\gamma(K,\Dom(A^k)))$ s.t.
\begin{equation}\label{eqn-converegnce}
\Vert \xi_n - \xi\Vert_{\mathcal{M}^{p}(\mathcal{L}^{2}([0,T],\gamma(K,H))} \to 0  \mbox{ sufficiently fast as } n \to \infty. \end{equation}
We denote the corresponding processes for $\xi_n$, from the previous step, by $\tilde{u}_n$ and $\tilde{\tilde{u}}_n$. By Step 1, for each $n$, the processes $\tilde{u}_n$ and $\tilde{\tilde{u}}_n$ satisfy the condition \eqref{eqn-equality almost sure}, the process $\tilde{u}_n$ satisfies inequality \eqref{StocStrichartzEstimate1} and the process $\tilde{\tilde{u}}_n$ satisfies inequality \eqref{StocStrichartzEstimate2}. Thus, both sequences are Cauchy in the appropriate Banach spaces
$\mathbb{M}^{p}({L}^{\infty}([0,T],\rV))$ and $M^{p}([0,T],E)$, respectively.
Hence, there exist unique elements in those spaces, whose representatives, respectively, we denote by
 $\tilde{u}$ and $\tilde{\tilde{u}}$. Because the convergence \eqref{eqn-converegnce} is sufficiently fast, we deduce that
 $\tP$-almost surely,
 $\tilde{\tilde{u}}_n \to \tilde{\tilde{u}}$ in $L^p(0,T;E)$ and
 $\tilde{u}_n \to \tilde{u}$ in $\mathrm{C}([0,T];\rV)$. Hence,  we infer that  $\tilde{u}$ is $\rV$-valued $\mathbb{F}$-adapted and  continuous process and $\tilde{\tilde{u}}$ is an $E$-valued progressively measurable process. Moreover, the processes $\tilde{u}$ and $\tilde{\tilde{u}}$ satisfy the condition \eqref{eqn-equality almost sure}.
Hence we are done with the proof of part (ii) of Theorem~\ref{thm-StocStrichartzEstimate}.

\end{proof}

\begin{proof}[\textbf{Proof of part (III) of Theorem \ref{thm-StocStrichartzEstimate}}]
This part  follows straightforwardly from the second part of Proposition \ref{prop-pointwise evaluation}.
\end{proof}

\section{Local well-posedness - preliminary results} \label{sec-localExistence}
The aim of this section is to formulate and prove some preliminary results which will be helpful in Section \ref{sec-main} where we show the existence and uniqueness of solutions to the stochastic wave equation \eqref{AbsWe}. Recall that we are working in the setting mentioned in the Subsection~\ref{subsec-mainAssump}.

Let us also recall the following definitions. For any $T>0$, we put
\begin{equation}\label{eqn-X_T+Z_T}
\Xa_T: = \mathrm{C}\bigl([0,T];\rV\bigr),  \; \Ya_T: = L^p(0,T;\rE), \; Q_T:= \mathrm{C}([0,T];H)    \textrm{ and }   \Za_T:= \Xa_T  \cap \Ya_T.
\end{equation}

Obviously,  $\Xa_T$, $\Ya_T$, $Q_T$   and $\Za_T$ are (separable) Banach spaces with naturally defined  norms, see e.g. \eqref{eqn-X_T-norm}-\eqref{eqn-Q_T-norm}.

By  $\mathbb{M}^{p}(\Za_T)$  we denote the Banach space of (equivalence classes) of all  $E$-valued progressively measurable processes $\{ u(t) :t \in [0,T]\}$ having a continuous $\rV$-valued modification and  satisfying
\begin{equation}
\label{eqn-MpYT-Norm}
\Vert u \Vert_{\mathbb{M}^{p}(\Za_T)}^p := \tE \bigl[ \Vert u \Vert_{\Za_T}^p \bigr] = \tE \bigl[ \Vert u \Vert_{\mathrm{C}([0,T];\rV)}^p + \Vert u \Vert_{L^p(0,T;\rE)}^p \bigr] < \infty.
\end{equation}

Analogously, by  $\mathbb{M}^{p}(Q_T)$  we denote the Banach space of  all  $H$-valued continuous and progressively measurable processes $\{ u(t) ,t \in [0,T]\}$  satisfying
\begin{equation}
\label{eqn-MpQT-Norm}
\Vert u \Vert_{\mathbb{M}^{p}(Q_T)}^p :=  \tE \bigl[ \Vert u \Vert_{\mathrm{C}([0,T];H)}^p  \bigr] < \infty.
\end{equation}

If $\tau$ is a bounded stopping time, by $\mathbb{M}^{p}(\Za_\tau)$ we mean the Banach space of (equivalence classes of) all progressively measurable processes
$	u :[0,\tau] \times \Omega \to E$
which have a continuous $\rV$-valued modification such that for each $\omega \in \Omega$, $u(\cdot,\omega) \in \Za_{\tau(\omega)}$ and
\begin{equation*}
	\int_\Omega	 \bigl[ \Vert u(\omega) \Vert_{\mathrm{C}([0,\tau(\omega)];\rV)}^p + \Vert u(\omega) \Vert_{L^p(0,\tau(\omega);\rE)}^p \bigr]\, d\tP (\omega) < \infty.
\end{equation*}

\subsection{SNLWE and assumptions}
Here we recall the stochastic nonlinear wave equation we consider here and state the assumptions on the drift and diffusion terms.  To be precise, we consider the following Cauchy problem for stochastic nonlinear wave equation with the Dirichlet or the Neumann boundary condition
\begin{numcases}{}
	u_{tt} +\nA u + F(u) = G(u) \dot{W},&  \label{eqn-Stoch wave main}\\
	u(0) = u_0, ~~  u_t(0) = u_1, & \label{eqn-Stoch wave main-IC}
\end{numcases}
where  $\nA$ is either $-\Delta_{D}$ or $-\Delta_{N}$; $(u_0,u_1) \in \rV \times L^2({\cO})$ and $W = \{W(t), t \geq 0  \}$  is a cylindrical Wiener process on some real separable Hilbert space $K$ such that $K \embed L^\infty({\cO}) $ and for some orthonormal basis $\{f_j\}_{j \in \mathbb{N}}$ of $K$,
\begin{equation}\label{KAssump}
 \sum\limits_{j \in \mathbb{N}} \| f_j \|_{L^\infty({\cO})}^2 < \infty.
\end{equation}
We assume the following hypotheses for the nonlinearity $F$ and the diffusion coefficient $G$ in equation \eqref{eqn-Stoch wave main}.
\begin{enumerate}[label=\textbf{A.\arabic*}]
 	\item \label{Fassump}
 	 Assume that a map
 	\begin{equation*}
 	F: B_{\rV}(0,1) \cap E \to  H,
 	\end{equation*}
 where for $R>0$,  $B_{\rV}(0,R)$ is the open ball in the space $\rV$, centered at the origin and of radius $R$, i.e.
 \[
 B_{\rV}(0,R):=\bigl\{u \in \rV: \Vert u \Vert_{\rV}<R \bigr\}, \]
 	is  such that $$F(0)=0,$$ and there exists a $\gamma_{F} >0$ such that for every $M \in (0,1)$  there exists a positive real number $C_F: = C_F(M)$,
   such that the following holds
  \begin{align}\label{ineq-F}
  \Vert F(u) - F(v) \Vert_H &  \le C_F(M) \left[ 1 + \| u\|_{E} + \| v \|_{E} \right]^{\gamma_{F}}~ \Vert u - v \Vert_{\rV},
  \end{align}
	 provided
  	\begin{equation}\label{eqn-FCondn}
 		u,v \in \bar{B}_{\rV}(0,M) \cap E.
 	\end{equation}

 	\item \label{Gassump}
 	 Assume  that a map
  		\[G: B_{\rV}(0,1) \cap E \to \gamma(K,H),\]
  	is  such that $$G(0)=0,$$ and  there exists a $\gamma_{G} >0$ such that for every   $M \in (0,1)$  there exists a positive real number $C_G(M): = C_G(M,\gamma_{G})$ such that if $u,v$ satisfy \eqref{eqn-FCondn}, then
	 \begin{align}\label{ineq-G}
		\Vert 	G( u ) - G( v) \Vert_{\gamma(K,H)} & \le C_G(M) \left[ 1 + \| u\|_{E} + \| v \|_{E} \right]^{\gamma_{G}}~ \Vert u - v \Vert_{\rV}.
 	\end{align}
	 \end{enumerate}

\begin{remark}\label{rem-gamma}
	Without loss of generality  we will assume that $\gamma_{F} = \gamma_{G} =: \gamma$.
\end{remark}

\begin{remark}\label{rem-extensionFandG}
	Note that, for given $M<1$, since $F$ and $G$ are generalized Lipschitz functions, i.e. satisfy \eqref{ineq-F} and \eqref{ineq-G}, respectively, on a closed  ball $\bar{B}_{\rV}(0,M)$, by Theorem \ref{thm-brz+rana_extension}, there exist maps $\tilde{F}=\tilde{F}_M$ and $\tilde{G}=\tilde{G}_M$, defined on $H_A \cap E$, taking values in spaces $H$ and $\gamma(K,H)$, respectively, such that
\begin{align}\label{eqn-extensions}
\tilde{F}_M(u)=F(u) \mbox{ and } \tilde{G}_M(u)=G(u) \mbox { if } u \in \bar{B}_{\rV}(0,M) \cap E,
\end{align}
and the inequalities \eqref{ineq-F} and \eqref{ineq-G} hold true for every $u,v \in H_A \cap E$. In particular, with Remark \ref{rem-gamma}, there exists a $\gamma >0$ such that for every $u, v \in H_A \cap E$
	\begin{align}
		& \Vert \tilde{F}_{M}(u) - \tilde{F}_{M}(v) \Vert_H   \le 3C_F(M) \left[ 1 + \| u\|_{E} + \| v \|_{E} \right]^{\gamma}~ \Vert u - v \Vert_{\rV}, \label{ineq-tildeF} \\
		& \Vert \tilde{G}_{M}( u ) - \tilde{G}_{M}( v) \Vert_{\gamma(K,H)}  \le 3C_G(M) \left[ 1 + \| u\|_{E} + \| v \|_{E} \right]^{\gamma}~ \Vert u - v \Vert_{\rV}. \label{ineq-tildeG}
	\end{align}

\end{remark}

\noindent The next two lemmata are a straightforward but important consequences of Remark \ref{rem-extensionFandG} and Assumptions~\ref{Fassump} and \ref{Gassump}.

\begin{lemma}\label{lem-F}
	Let us assume that  $M \in (0,1)$. Let us assume that the function  $F: \bar{B}_{\rV}(0,M)  \cap E \to H$  satisfies  Assumption~\ref{Fassump} with $\gamma = \gamma_{F}$ and let $T > 0$ and $p \geq \gamma$.  Let $\tilde{F}=\tilde{F}_M$ be the extension of $F$ as introduced  in Remark \ref{rem-extensionFandG} above.
	Then there exists a number $K(\gamma)>0$ such that the following inequality holds, 	provided $	u_1, u_2 \in \Za_T$,
	\begin{align}
		&	\Vert \tilde{F}_M(u_1) - \tilde{F}_M(u_2) \Vert_{L^1(0,T; H)} \leq 3K(\gamma) C_{F}(M) \left[ T + T^{1-\frac{\gamma}{p}} \bigl( \Vert u_1 \Vert_{\Ya_T}^{\gamma} +  \Vert u_2 \Vert_{\Ya_T}^{\gamma}\bigr)\right]~ \Vert u_1 - u_2 \Vert_{\Xa_T}.
		\nonumber	\end{align}

\end{lemma}
\begin{remark}\label{rem-lem-F}
Lemma \ref{lem-F} means that the function  $\tilde{F}:X_T \to L^1(0,T; H)$ is uniformly Lipschitz on the sets $\{ u \in X_T: \Vert u \Vert_{\Ya_T} \leq R\}$, $R>0$.
Analogous remark can be made for Lemma~\ref{lem-G}.
\end{remark}

\begin{proof}[\textbf{Proof of Lemma \ref{lem-F}}]
	
	Let us choose and fix $u_1,u_2 \in \Ya_T \cap \Xa_T=\Za_T$.
	Then,
	by using inequality \eqref{ineq-tildeF}, followed by the H\"older inequality, we get
	 \begin{align}
		& \Vert \tilde{F}(u_1) - \tilde{F}(u_2) \Vert_{L^1(0,T; H)}  \leq 3 C_{F}(M)   \int_{0}^{T} \bigl[ 1 + \| u_1(t)\|_{E} + \| u_2(t) \|_{E} \bigr]^\gamma~ \Vert u_1(t) - u_2(t) \Vert_{\rV}  \, dt \nonumber\\
		& \leq K(\gamma) 3 C_{F}(M)  \Vert u_1 - u_2 \Vert_{\Xa_T} \left[ T + T^{1-\frac{\gamma}{p}} \left( \int_{0}^{T} \Vert u_1(t) \Vert_E^p \, dt\right)^{\frac{\gamma}{p}}  + T^{1-\frac{\gamma}{p}}  \left( \int_{0}^{T} \Vert u_2(t) \Vert_E^p\, dt\right)^{\frac{\gamma}{p}}\right].  \nonumber
	\end{align}  Hence  Lemma \ref{lem-F} follows.
\end{proof}

\begin{lemma}\label{lem-G}
	Let us assume that  $M \in (0,1)$. Let us assume that the function  $G: \bar{B}_{\rV}(0,M)  \cap E \to \gamma(K,H)$   satisfies   Assumption~\ref{Gassump} with $\gamma = \gamma_{G}$ and let $T > 0$ and $p \geq 2\gamma$. Let $\tilde{G}$ be the extension of $G$ as in Remark \ref{rem-extensionFandG}.
	  Then there exists a real number $K(2\gamma)>0$ such that  the following inequality holds,
	provided
	$	u_1$ and  $u_2$ belong to $ \Za_T$,	%
	\begin{align}
		& \Vert \tilde{G}_M(u_1) - \tilde{G}_M(u_2) \Vert_{L^2(0,T; \gamma(K,H))}^2 \leq 9 K(2\gamma) C_{G}^2(M) \bigl[ T + T^{1-\frac{2\gamma}{p}} \bigl( \Vert u_1 \Vert_{\Ya_T}^{2\gamma}
		+  \Vert u_2 \Vert_{\Ya_T}^{2\gamma}\bigr)\bigr] \Vert u_1 - u_2 \Vert_{\Xa_T}^2, \nonumber
	\end{align}
\end{lemma}
\begin{proof}[\textbf{Proof of Lemma \ref{lem-G}}]
	Let us choose and fix $u_1,u_2 \in \Ya_T \cap \Xa_T$. Then, invoking inequality \ref{ineq-tildeG} and the H\"older inequality, we obtain
	\begin{align}
		& \Vert \tilde{G}(u_1) - \tilde{G}(u_2)\Vert_{L^2(0,T; \gamma(K,H))}^2  \leq 9 C_{G}^2(M)  \int_{0}^{T} \left[ 1 + \| u_1(t)\|_{E} + \| u_2(t) \|_{E} \right]^{2\gamma}~ \Vert u_1(t) - u_2(t) \Vert_{\rV}^2 \, dt \nonumber\\
		& \leq 9 K(2\gamma) C_{G}^2(M) \Vert u_1 - u_2 \Vert_{\Xa_T}^2 \left[ T + T^{1-\frac{2\gamma}{p}} \left( \int_{0}^{T} \Vert u_1(t) \Vert_E^p \, dt\right)^{\frac{2\gamma}{p}}  + T^{1-\frac{2\gamma}{p}} \left( \int_{0}^{T} \Vert u_2(t) \Vert_E^p\, dt\right)^{\frac{2\gamma}{p}}\right].  \nonumber
	\end{align}  Hence the proof of Lemma \ref{lem-G} is complete.
\end{proof}
We set the notation $C_{\tilde{F}_M} := 3C_{F}(M) $  and $C_{\tilde{G}_M}:= 3 C_{G}(M)$ and use them in Section \ref{sec-approximation problem}.  The Lemmata \ref{lem-F} and \ref{lem-G} have been formulated in the language of the extensions $\tilde{F}$ and $\tilde{G}$ of functions $F$ and $G$. The following corollary is written in the language of the original functions $F$ and $G$.
\begin{corollary}\label{cor-F and G}
	Let  $T > 0$ and $M \in (0,1)$.  Let us assume that the function  $F:B_{\rV}(0,1) \cap E \to H$  satisfies  Assumption~\ref{Fassump} with $\gamma = \gamma_{F}$ and $p \geq \gamma$.
Then there exists $K(\gamma)>0$ such that the following inequality holds
	\begin{align}
		&	\Vert {F}(u_1) - {F}(u_2) \Vert_{L^1(0,T; H)} \leq K(\gamma) C_{F}(M) \Bigl[ T + T^{1-\frac{\gamma}{p}} \Bigl( \frac{\Vert u_1 \Vert_{\Ya_T}^{\gamma}}{M^{\gamma}} +
\frac{\Vert u_2 \Vert_{\Ya_T}^{\gamma}}{M^{\gamma}} \Bigr)\Bigr]~ \Vert u_1 - u_2 \Vert_{\Xa_T},
		\nonumber	\end{align} 	
for all 	$	u_1, u_2 \in \Za_T$ satisfying
  	\begin{equation}\label{eqn-FCondn-2}
 		 \sup_{t \in [0,T]}\Vert u(t) \Vert_{\rV} \leq M \mbox{ and } \sup_{t \in [0,T]} \Vert v(t) \Vert_{\rV} \leq M.
 	\end{equation}
 	Let us assume that the function  $G:B_{\rV}(0,1) \cap E \to \gamma(K,H)$   satisfies   Assumption~\ref{Gassump} with $\gamma = \gamma_{G}$ and $p \geq 2\gamma$.
	Then there exists $K(2\gamma)>0$ such that  the following inequality holds
	\begin{align}
		& \Vert {G}(u_1) - {G}(u_2) \Vert_{L^2(0,T; \gamma(K,H))}^2 \leq K(2\gamma) C^2_{G}(M) \Bigl[ T + T^{1-\frac{2\gamma}{p}} \Bigl(
\frac{\Vert u_1 \Vert_{\Ya_T}^{2\gamma}}{M^{2\gamma}}
		+  \frac{\Vert u_2 \Vert_{\Ya_T}^{2\gamma}}{M^{2\gamma}} \Bigr)\Bigr] \Vert u_1 - u_2 \Vert_{\Xa_T}^2, \nonumber
	\end{align}
for all 	$	u_1, u_2 \in \Za_T$ satisfying condition \eqref{eqn-FCondn-2}.
 \end{corollary}

To prove the main result of this Section \ref{sec-localExistence} we need the following known results. The first one is from \cite{R2005}.
\begin{theorem}\label{thm-Moser Trudinger Ineq}[Moser-Trudinger Inequality]\\
	Let ${\cO} \subseteq \mathbb{R}^2$ be a domain (bounded or unbounded) and $\alpha \leq 4 \pi$. Then
	\begin{equation}\label{MTIneq}
	C(\alpha,{\cO}):= \sup_{  \substack{u \in H^{1,2}({\cO}), \\  \Vert u \Vert_{H^{1,2}({\cO})} \leq 1 }   } \int_{{\cO}} \bigl(e^{\alpha (u(x))^2} - 1 \bigr)  \, dx   < + \infty.
	\end{equation}
	Moreover, this result is sharp in the sense that if $\alpha > 4\pi$ then $C(\alpha,{\cO})=\infty$.
\end{theorem}

The next required result is the well-known Logarithmic inequality from \cite{O1995}.
\begin{theorem}\label{thm-LogIneqThm}
	Let ${\cO}$ be a domain in $\mathbb{R}^d$. Let $p,q,m \in \mathbb{R}$ satisfy $1<p< \infty, 1 \leq q < \infty$, and $m > \frac{d}{q}$. Then there exists a positive constant $L$ such that for all $u \in H^{\frac{d}{p},p}({\cO}) \cap H^{m,q}({\cO})$  the following holds,
	\begin{equation}\label{LogIneq}
	\Vert u \Vert_{L^\infty({\cO})} \leq L \Vert u \Vert_{H^{\frac{d}{p},p}({\cO})} \left[ 1+  \log\bigl( 1 + \frac{\Vert u \Vert_{H^{m,q}({\cO})}}{\Vert u \Vert_{H^{\frac{d}{p},p}({\cO})}} \bigr) \right]^{1-\frac{1}{p}}.
	\end{equation}
\end{theorem}

In the next two results we provide an example of functions $f$ and $g$ such that the corresponding $F$ and $G$ satisfy the assumptions \ref{Fassump} and \ref{Gassump}, respectively. The example below comes from \cite{IMM2006} and \cite{IJ2011} when $E$ is a suitable H\"older space.
We will prove the next result in detail because we need a slightly more general version of the Moser-Trudinger inequality and the Logarithmic estimate, respectively, see Theorem \ref{thm-Moser Trudinger Ineq} and \ref{thm-LogIneqThm}, than those used in \cite{IMM2006} and \cite{IJ2011}.
\begin{lemma}\label{lem-fandgExample} Assume that ${\cO} \subseteq \mathbb{R}^2$ is a bounded domain.
	Let $L$ be a constant from the Theorem \ref{thm-LogIneqThm}. Let $h: \mathbb{R} \to \mathbb{R}$ be a function defined by $h(x) = \pm x \bigl( e^{4 \pi x^2} - 1\bigr)$ for $x \in \mathbb{R}$.
Assume that a pair  $(q,r)$ of positive numbers satisfies
	\begin{equation}\label{pairCondn}
	0 < r + \frac{2}{q}< 1.   	\end{equation}
Then  for every $M \in (0,1)$
	there exists a constant $C_{M,L} >0$, which depending only on $M$ and $L$, such that
	\begin{equation}\label{eqn-Nemytski}
		\Vert h \circ u - h \circ v \Vert_{L^{2}({\cO})} \leq C_{M,L}  \left[ 1 + \| u\|_{H^{1-r,q}({\cO})} + \| v \|_{H^{1-r,q}({\cO})} \right]^{2\pi L^2}~ \Vert u - v \Vert_{H^{1,2}({\cO})},
	\end{equation}
	provided $u,v\in H^{1,2}({\cO}) \cap H^{1-r,q}({\cO})$ satisfy the following condition \begin{equation}\label{Mcondn}
		\|u\|_{H^{1,2}({\cO})} \leq M \quad \textrm{ and } \quad \|v \|_{H^{1,2}({\cO})} \leq M.
	\end{equation}
\end{lemma}

In the  next result, which is about a generalized Nemytskii operator $G$ associated with function $g$,  $\{f_j\}_{j \in \mathbb{N}}$ is an ONB of a Hilbert space $K$.
\begin{lemma}\label{lem-G1}
	Assume that  condition \eqref{KAssump} holds. Assume that $g(x) = x \bigl( e^{4 \pi x^2} - 1\bigr)$, $x \in \mathbb{R}$ and a pair  $(q,r)$ exists which satisfies condition \eqref{pairCondn}.
 Let $L$ is a constant from the Theorem \ref{thm-LogIneqThm} and $G$ be  defined by
	\begin{equation*}
 		G(u):=\{ K \ni k \mapsto (g \circ u) \cdot k \in L^{2}({\cO}) \}, \;\; u \in B_{H^{1,2}({\cO})}(0,1) \cap H^{1-r,q}({\cO}) .
 	\end{equation*}
	Then for every $M \in (0,1)$  the following inequality holds
	\begin{equation*}
		\Vert G(u)  - G(v) \Vert_{\gamma(K,L^{2}({\cO}))} \leq \bar{C}_G \left[ 1 + \| u\|_{H^{1-r,q}({\cO})}
+ \| v \|_{H^{1-r,q}({\cO})} \right]^{2\pi L^2}~ \Vert u - v \Vert_{H^{1,2}({\cO})},
	\end{equation*}
	where
	 \[\bar{C}_G := C_{M,L} \times \bigl[ \sum\limits_{j \in \mathbb{N}} \ \| f_j \|_{L^\infty({\cO})}^2 \bigr]^{\frac{1}{2}},\]
for all $u,v \in H^{1,2}({\cO}) \cap H^{1-r,q}({\cO})$ satisfying condition  \eqref{Mcondn}.
\end{lemma}
\begin{remark}\label{rem-Lem-F+G}
Both Lemmata \ref{lem-fandgExample} and \ref{lem-G} are applicable to spaces defined in \eqref{eqn-spaces} because
$$ \rV \subset H^{1,2}({\cO})  \mbox{ and } \rE  \subset H^{1-r,q}({\cO}).$$

\end{remark}
\begin{proof}[\textbf{Proof of Lemma \ref{lem-G}}]
	Let $u$ and $v$ belong to  $H^{1,2}({\cO}) \cap H^{1-r,q}({\cO})$ satisfying condition  \eqref{Mcondn}. By assumption \eqref{KAssump} and Lemma \ref{lem-fandgExample} (applied to $h=g$)  we infer that

 \begin{align}
	& \| G(u)- G(v)  \|_{\gamma(K,L^2({\cO}))}^2   = \sum\limits_{j \in \mathbb{N}} \| G(u)f_j - G(v) f_j\|_{L^2({\cO})}^2\nonumber\\
	& = \sum\limits_{j \in \mathbb{N}} \| (g \circ u) f_j - (g \circ v) f_j \|_{L^2({\cO})}^2 \leq \| g \circ u - g \circ v\|_{L^2({\cO})}^2 \sum\limits_{j \in \mathbb{N}} \ \| f_j \|_{L^\infty({\cO})}^2.  \nonumber
	\end{align}
	as required. Hence the result follows by applying the inequality \eqref{eqn-Nemytski}.
\end{proof}

\begin{proof}[\textbf{Proof of Lemma \ref{lem-fandgExample}}]
We only prove the result for $h(x) = x(e^{4\pi x^2}-1)$, since the proof for the function $-x(e^{4\pi x^2}-1)$ is analogous. We begin here with the following observation which is a consequence of  the Mean Value Theorem. If $u,v\in \mathbb{R}$, then the following equality holds
	\begin{equation*}
	h(u)-h(v)=	(u-v) \int_0^1  \left[ \bigl(1+ 8 \pi u_\theta^2 \bigr)e^{4 \pi u_\theta^2}  - 1\right] \; d\theta.
	\end{equation*}
Let us now fix $M \in (0,1)$ and choose $\zeta \in (0,1)$ and $\varepsilon >0$, such that
	\begin{equation}\label{epsAndZetaBound}
		(1+\varepsilon)(1 +\zeta)M^2 \leq 1.
	\end{equation}
Let us take arbitrary $u, v \in H^{1,2}({\cO}) \cap H^{1-r,q}({\cO})$ satisfying  condition \eqref{Mcondn}.
Applying the above for $u=u(x)$ and $v=v(x)$, for a fixed $x \in \cO$ we get, with $u_\theta(x):=(1-\theta)u(x) + \theta v(x)$,
	\begin{equation*}
	h(u(x))-h(v(x))= (u(x)-v(x)) \int_0^1  \left[ \bigl(1+ 8 \pi u_\theta^2(x) \bigr)e^{4 \pi u_\theta^2(x)}  - 1\right] \; d\theta. \;\;
	\end{equation*}
  Thus, we infer that
\begin{equation*}
	\Vert h \circ u - h \circ v \Vert_{L^2({\cO})} \leq \bigg\Vert (u-v) \int_0^1 \left[ \bigl(1+ 8 \pi u_\theta^2 \bigr)e^{4 \pi u_\theta^2}  - 1\right]\, d \theta  \bigg\Vert_{L^2({\cO})}.
	\end{equation*}
	Applying the Minkowski inequality gives
	\begin{align}\label{fandgExample-Eqn5}
	\Vert (u-v) & \int_0^1  \left[ \bigl(1+ 8 \pi u_\theta^2 \bigr)e^{4 \pi u_\theta^2}  - 1\right] \; d\theta\Vert_{L^2({\cO})} \nonumber\\
	& \qquad  \leq
	\int_0^1 \left\Vert  (u-v)  \left[ \bigl(1+ 8 \pi u_\theta^2 \bigr)e^{4 \pi u_\theta^2}  - 1\right] \right\Vert_{L^2({\cO})} \; d\theta .
	\end{align}
 Then due to the H\"older inequality, the  Sobolev embedding and  the following basic inequality
 \begin{equation*}
 (1+2a)e^a -1 \leq 2 \bigl(1 + \frac{1}{\varepsilon} \bigr) \bigl(e^{(1+\varepsilon)a} -1 \bigr), \quad  a,\varepsilon >0,
 \end{equation*}	
 we infer that
	\begin{align}\label{fandgExample-Eqn7}
	& \bigg\Vert   (u-v) \left[ \bigl(1+ 8 \pi u_\theta^2 \bigr)e^{4 \pi u_\theta^2}  - 1\right] \bigg\Vert_{L^2({\cO})}^2  \leq 4 \bigl(1 + \frac{1}{\varepsilon} \bigr)^2 \bigg\Vert (u-v) \bigl(e^{4 \pi (1+\varepsilon) u_\theta^2} - 1 \bigr)   \bigg\Vert_{L^2({\cO})}^2  \nonumber\\
	& \le 4\bigl(1 + \frac{1}{\varepsilon} \bigr)^2 \Vert u-v \Vert_{L^{2 + \frac{2}{\zeta}}({\cO})}^2 \bigg\Vert \bigl( e^{4 \pi (1+\varepsilon) u_\theta^2} - 1 \bigr)^2  \bigg\Vert_{L^{1+\zeta}({\cO})}  \nonumber\\
	& \le  4C_\zeta \bigl(1 + \frac{1}{\varepsilon} \bigr)^2   \Vert u-v\Vert_{\Dom(\nA^{\frac12 })}^2  ~ e^{4 \pi (1+\varepsilon) \Vert u_\theta^2\Vert_{L^\infty({\cO})} } \Vert e^{4 \pi (1+\varepsilon) u_\theta^2} - 1   \Vert_{L^{1+\zeta}({\cO})}.
	\end{align}
	Moreover, since $\theta \in (0,1)$, $\Vert  u_\theta \Vert_{H^{1,2}({\cO})}   \leq  M$ and therefore \eqref{epsAndZetaBound} holds.  Thus, the Moser-Trudinger inequality from Theorem \ref{thm-Moser Trudinger Ineq} gives
	\begin{align}\label{fandgExample-Eqn8}
	& \big\Vert  e^{4\pi (1+\varepsilon) u_\theta^2} -1\big\Vert_{L^{1 + \zeta}({\cO})}^{1+ \zeta} \leq  \big\Vert e^{4\pi (1+\varepsilon)(1+\zeta)u_\theta^2} -1\big\Vert_{L^1({\cO})}   \leq  C(4\pi,{\cO}).
	\end{align}
 	Invoking the  estimate from Theorem \ref{thm-LogIneqThm}, which is possible due to \eqref{pairCondn} and Lemma \ref{lem-EquiSpace1},
 	we obtain  \begin{align}
	e^{ 4 \pi (1+\varepsilon) \Vert u_\theta\Vert_{L^\infty({\cO})}^2} \leq  \exp\left[ 4 \pi L^2 ( 1+\varepsilon) \Vert u_\theta \Vert_{H^{1,2}({\cO})}^2 \left\{ 1+ \log\bigl( 1 + \frac{\Vert u_\theta \Vert_{H^{1-r,q}({\cO}) }}{\Vert u_\theta\Vert_{H^{1,2}({\cO})}} \bigr)  \right\}     \right]. \nonumber
	\end{align}

	 \noindent Using the fact that if $B_1 > 1, B_2 > 0$, then the function $x \mapsto x^2 \Bigl( 1 + \log \bigl( B_1 + \frac{B_2}{x}\bigr)\Bigr)$ is non-decreasing, we deduce that,
	\begin{align}\label{fandgExample-Eqn10}
	e^{ 4 \pi (1+\varepsilon) \Vert u_\theta\Vert_{L^\infty({\cO})}^2} & \leq  \left[ e \bigl(  1 + \frac{\Vert u_\theta \Vert_{H^{1-r,q}({\cO})}}{M}  \bigr) \right]^{4 \pi L^2 ( 1+\varepsilon) M^2}.
	\end{align}
	Let us put  \begin{equation}	\label{gamma}
		\gamma := 2 \pi L^2(1+\varepsilon) M^2.
	\end{equation}
Note that, since $(1+\varepsilon) M^2 <1$, from \eqref{epsAndZetaBound}  we get that $\gamma < 2 \pi L^2$.  Next,  from \eqref{fandgExample-Eqn5}, \eqref{fandgExample-Eqn7}, \eqref{fandgExample-Eqn8} and  \eqref{fandgExample-Eqn10}, we infer that
	\begin{align}
		 \Vert h \circ u - h \circ v \Vert_{L^2({\cO})}  & \le 2 \sqrt{C_{\zeta}} \bigl(1 + \frac{1}{\varepsilon} \bigr) \bigl( \frac{e}{M} \bigr)^{2\pi L^2} (C(4 \pi,{\cO}))^{\frac{1}{2(1+\zeta)}}  \Vert u-v \Vert_{\Dom(\nA^{\frac12 })} \times \nonumber\\
		& \qquad \times \int_0^1 \bigl(  1 + \| u\|_{H^{1-r,q}({\cO})}
		+ \| v \|_{H^{1-r,q}({\cO})}    \bigr)^{2\pi L^2} \, d \theta . \nonumber
	\end{align} Hence Lemma \ref{lem-fandgExample} follows.
\end{proof}

\begin{remark}\label{rem-polynomial} For polynomial functions a stronger version of previous two lemmata hold,  see e.g. \cite{Brz+Millet_2014}.\\
Let us also observe that we proved that if $h(x) = x(e^{4\pi x^2}-1)$, then
$ h \circ u  \in L^{2}({\cO})$ for every  $u\in B_{H^{1,2}({\cO})}(0,1) \cap H^{1-r,q}({\cO})$.  If $u\in H^{1,2}({\cO}) \cap H^{1-r,q}({\cO})$ and $\Vert u \Vert_{H^{1,2}({\cO})} \geq 1$, then $ h \circ u$ is a well defined function  but we not know if it belongs to the space $L^{2}({\cO})$. One can ask if the following assertions are  true: if $v\in H^{1,2}({\cO}) \cap H^{1-r,q}({\cO})$ and  $M=\Vert v \Vert_{H^{1,2}({\cO})} \in (0, 1)$, then
\begin{align} \label{question-1}
\lim_{ u\in H^{1,2}({\cO}) \cap H^{1-r,q}({\cO}), \Vert u \Vert_{H^{1,2}({\cO})} \toup 1  } & \Vert h \circ u  \Vert_{ L^{2}({\cO})}=\infty?
\\
\label{question-2}
\lim_{ r\toup M^{-1} }& \Vert h \circ (rv)  \Vert_{ L^{2}({\cO})}=\infty?
\end{align}
\end{remark}

\subsection{Definition of a local mild solution}
In this subsection we introduce the definitions of local and maximal local solutions  we adopt in this paper. They are modifications of definitions used in earlier papers, see e.g. \cite{Brz+M+S_2005}.

\begin{definition}\label{del-local solution}
Let $p \in [2,\infty)$.  Assume that
 $(u_0,u_1)$ is an $\mathcal{F}_0$-measurable  $ H_A  \times H$  random variable which   satisfies the following condition
\begin{align}\label{eqn-u_0}
	&\Vert u_0 \Vert_{H_A}  < 1, \;\;\; \mathbb{P}\mbox{-almost surely,}
\end{align}
\begin{itemize}
\item[\textbf{A.}]	A local mild solution to problem \eqref{eqn-Stoch wave main}-\eqref{eqn-Stoch wave main-IC} is a $\Dom(\nA^{\frac12 })$-valued continuous and $\mathbb{F}$-adapted process $u= \{  u(t): t \in [0,\tau) \}$ satisfying the following conditions
\begin{enumerate}
		\item $\tau$ is an accessible stopping time,		
\item  the condition \eqref{eqn-u_0} is preserved, i.e.
\begin{equation}\label{eq-norm-preservance}
		\Vert u(t) \Vert_{H_A}  < 1, \quad  \mbox{for $t \in [0,\tau)$,} \;\; \tP\textrm{-a.s.},
		\end{equation}
		\item there exists an  announcing sequence $\{ \tau_k \}_{k \geq 1}$ of the stopping times for  $\tau$, such that
		\begin{equation*}
		\mbox{
		$u$ belongs to $\mathbb{M}^{p}(\Za_{t \wedge \tau_k})$, for all $t\geq 0$ and every $k$,}		
		\end{equation*}
		and, for all $t \geq 0$ and $k \in \mathbb{N}$,
\begin{align}\label{eqn-StocMildSolnDefn}
		\quad \quad u(t \wedge \tau_k) & = \cos((t \wedge \tau_k)\sqrt{\nA})u_0 + \frac{\sin((t \wedge \tau_k) \sqrt{\nA})}{\sqrt{\nA}}u_1 \nonumber\\
		& \quad + \int_{0}^{t \wedge \tau_k} \frac{\sin((t\wedge \tau_k -s)  \sqrt{\nA})}{\sqrt{\nA}} F(u(s)) \, ds + I(\tau_k,G)(t \wedge \tau_k), \;\;\tP\textrm{-a.s.},
		\end{align}
		where $I(\tau_k,G)$ is a process  defined by
\begin{equation}\label{eqn-stoppedStocConv}
		I(\tau_k,G)(t) = \int_{0}^{t} \mathds{1}_{[0,\tau_k)}(s) \frac{\sin((t-s) \sqrt{\nA})}{\sqrt{\nA}} G(u(s)) \, dW(s), \;\; t \geq 0.
		\end{equation}
	\end{enumerate}
	\item[\textbf{B.}]	A local mild solution $u= \{  u(t):t \in [0,\tau) \}$ to problem \eqref{eqn-Stoch wave main}-\eqref{eqn-Stoch wave main-IC} is unique iff for any other local solution $\hat{u}= \{  \hat{u}(t): t \in [0,\hat{\tau}) \}$ to problem \eqref{eqn-Stoch wave main}-\eqref{eqn-Stoch wave main-IC},  the restricted processes $\restr{u}{[0,\tau \wedge \hat{\tau}) \times \Omega}$ and $\restr{\hat{u}}{[0,\tau \wedge \hat{\tau} ) \times \Omega}$ are equivalent.

\item[\textbf{C.}] A local mild solution $u= \{  u(t): t \in [0,\tau) \}$ to problem \eqref{eqn-Stoch wave main}-\eqref{eqn-Stoch wave main-IC} is not maximal iff  there exists a  local solution $\hat{u}= \{  \hat{u}(t): t \in [0,\hat{\tau}) \}$ to problem \eqref{eqn-Stoch wave main}-\eqref{eqn-Stoch wave main-IC} such that $\tP(\tau \leq \hat{\tau})=1$, $\tP(\tau < \hat{\tau})>0$ and processes $u$ and $\restr{\hat{u}}{[0,\tau ) \times \Omega}$ are equivalent.
\end{itemize}
\end{definition}

\begin{remark}\label{rem-local solution} The definition of the process $I_{\tau_n}(G)$ is explained in Lemma \ref{lem-A.1} of Appendix \ref{sec-stopped processes}. The use of processes $I_{\tau_n}(G)$  was first introduced for the SPDEs of parabolic type in \cite{Brz+G_1999} and \cite{Carroll1999T} and in \cite{Brz+M+S_2005} for the hyperbolic SPDEs. The definition we use above is only in terms of the process $u$ and thus it is different from the one used in \cite{Brz+M+S_2005} which is in terms of pair processes $(u,u_t)$. In Appendix \ref{sec-def sol} we discuss the equivalence between these two approaches. \end{remark}

\begin{remark}\label{rem-local solution-2} It can be shown, see Remark 2.22 in \cite{Brz+E_2000}, that if  an $\mathbb{F}$-adapted and $\Dom(\nA^{\frac12 })$-valued continuous process $u= \{  u(t): t \in [0,\tau) \}$ is a  local mild solution to problem \eqref{eqn-Stoch wave main}-\eqref{eqn-Stoch wave main-IC}, then equality \eqref{eqn-StocMildSolnDefn} with \eqref{eqn-stoppedStocConv} hold for any stopping time of the
form $\tau_k\wedge \sigma$, where $\sigma$ is an accessible stopping time.
\end{remark}

 It can be shown that the concept of a local maximal solution introduced in part \textbf{C.} of Definition \ref{del-local solution} is equivalent to the following set theoretical one.
A natural continuation of this new definition is the so called ``Amalgamation Lemma", see  \cite[Lemma III 6A and 6B]{Elworthy-1982B} and  Definition 3.11 in \cite{Brz+H+R_2020-liquid crystals penalised}.
	\begin{definition}\label{def-maxsol-0}
		Let us denote the set of all local solutions  $(u,\tau)$ to  the problem \eqref{eqn-Stoch wave main}-\eqref{eqn-Stoch wave main-IC} by $\mathcal{ LS}$. For any two
		elements $(u,\tau), (v,\sigma) \in \mathcal{ LS} $ we write that
		$(u,\tau)\preceq (v,\sigma)$ iff $\tau \leq \sigma$ $\tP$-a.s. and
		$v_{\vert [0,\tau)\times \Omega} \sim u$, see Definition \ref{defn-EquivSP} for the notation $\sim$. We write $(u,\tau)\prec
		(v,\sigma)$ iff $(u,\tau)\preceq (v,\sigma)$ and $(u,\tau)\not\sim
		(v,\sigma)$. It is straightforward to show that $\preceq$ is a partial order on  $\mathcal{ LS}$.
		
		We say that $(u,\tau)$ is a maximal element of $(\mathcal{ LS},\preceq)$ iff there is no $(v,\sigma) \in (\mathcal{ LS},\preceq)$ such that $(u,\tau)\prec
		(v,\sigma)$.
		Each maximal element $(u,\tau)$ in the set $(\mathcal{
			LS},\preceq)$
		is called a maximal local solution to  the problem  \eqref{eqn-Stoch wave main}-\eqref{eqn-Stoch wave main-IC}.
	\end{definition}

\section{The approximation problem}\label{sec-approximation problem}
In this section we will study the following approximated version of problem  \eqref{eqn-Stoch wave main}-\eqref{eqn-Stoch wave main-IC}, for $n\in \mathbb{N}\setminus\{0\}$,
\begin{numcases}{}
	u_{tt}(t) +\nA u(t) +  \theta_n(\Vert u \Vert_{\Za_t}) \tilde{F}_{n}(u(t)) = \theta_n(\Vert u \Vert_{\Za_t})   \tilde{G}_{n}(u(t)) \dot{W}(t), &  \label{eqn-Stoch wave approx}\\
	u(0) = u_0, ~~  u_t(0) = u_1, & \label{eqn-Stoch wave approx-IC}
\end{numcases}
where $\theta : \mathbb{R}_+ \to [0, 1]$ be an auxiliary smooth function with compact support such that
\begin{equation*}
	\inf_{x \in \mathbb{R}_+} \theta^\prime(x) \geq -2, \quad 	\mathds{1}_{[0,1]} \le \theta \le \mathds{1}_{[0,2)},
\end{equation*}
and for $n \geq 1$ set $\theta_n (\cdot) = \theta\bigl( \frac{\cdot}{n}\bigr)$.

The following lemma states the basic properties of $\theta_n$.
\begin{lemma}\label{localExist2-Lem1}
	The functions $\theta_n$ is Lipschitz and bounded and, for all $x,y\in [0,\infty)$,
	\begin{equation*}
		\vert \theta_n(x) - \theta_n(y) \vert \leq \frac{ 2}{n} \vert x-y \vert.
	\end{equation*}
	Moreover, if $h: \mathbb{R}_+ \to \mathbb{R}_+$ is a non decreasing function, then for all $x\in [0,\infty)$,
	\begin{equation*}
		\theta_n(x) h(x) \leq h(2n).
	\end{equation*}
\end{lemma}

\begin{remark}
	Let us  point out that there is a typo in the lower bound of $\inf_{x \in \mathbb{R}_+} \theta(x)$ in \cite[condition (4.10)]{Brz+Millet_2014} and the value of this lower bound should be strictly smaller than $-1$. Indeed,   we can easily show that there does not exist  a smooth function  satisfying  \cite[condition (4.10)]{Brz+Millet_2014}.  Consequently, the Lipschitz constant in \cite[Lemma ~4.3]{Brz+Millet_2014} should be strictly greater than $\frac{1}{n}$.
\end{remark}
To move ahead let us consider a numerical sequence $M=(M_n)_{n=1}^\infty=(1-\frac1{n+1})_{n=1}^\infty= (\frac n{n+1})_{n=1}^\infty$.

The approximating problem \eqref{eqn-Stoch wave approx} approximates the original problem \eqref{eqn-Stoch wave main} in a two-fold way. Firstly,  the very bad non-linearities (of possibly exponential growth)
$F$ and $G$ are  replaced   by nonlinearities $\tilde{F}_{n}$ and $\tilde{G}_{n}$ depending  on $M_n$ which are globally Lipschitz and bounded in the sense of  Remark \ref{rem-extensionFandG} and Theorem \ref{thm-brz+rana_extension}.
Secondly, these new nonlinearities are multiplied  by a suitable cut-off function   $\theta_n(\Vert u \Vert_{\Za_t})$ which depends on the whole history up to time $t$ of the solution. Note however that this cut-off function  takes values in the interval $[0,1]$.
The second modification makes the coefficients of the approximating problem \eqref{eqn-Stoch wave approx} not only time dependent but also random. This two-fold modification makes the coefficients in problem \eqref{eqn-Stoch wave approx} enough Lipschitz  so that the corresponding fixed point problem has a unique solution  and problem \eqref{eqn-Stoch wave approx} has a unique global  solution.

\begin{definition} \label{def-solution approximating equation}
	Assume that  $(u_0,u_1)$   is an $H_A  \times H$-valued $\mathcal{F}_0$-measurable  random variable. 	A mild solution to problem  \eqref{eqn-Stoch wave approx}-\eqref{eqn-Stoch wave approx-IC} on time interval $[0,T]$, where $T>0$, is a process
	$u $ belonging to space $\mathbb{M}^{p}(\Za_T)$ such that for every $t\in [0,T]$, the following equality holds $\mathbb{P}$-almost surely,
	\begin{align}\label{eqn-sol approx eqn}
		u(t)  = \cos(t \sqrt{\nA}) u_0 + \frac{\sin(t \sqrt{\nA})}{\sqrt{\nA}} u_1
		&+ \int_{0}^{t}  \theta_n(\Vert u \Vert_{\Za_s}) \, \frac{\sin((t-s) \sqrt{\nA})}{\sqrt{\nA}}   \tilde{F}_{n}(u(s)) \, ds\\
		&  + \int_{0}^{t} (\theta_n(\Vert u \Vert_{\Za_s}))\frac{\sin((t-s)\sqrt{\nA})}{\sqrt{\nA}}    \tilde{G}_{n}(u(s)) \, dW(s).
		\nonumber
	\end{align}
	A global mild solution to problem \eqref{eqn-Stoch wave approx}-\eqref{eqn-Stoch wave approx-IC} is a process $u$  such that for every $T>0$,
	$u \in \mathbb{M}^{p}(\Za_T)$  is a mild solution to problem \eqref{eqn-Stoch wave approx}-\eqref{eqn-Stoch wave approx-IC} on time interval $[0,T]$.
\end{definition}

The following result is a consequence of Proposition \ref{prop-equivalence}.
\begin{lemma}\label{lem-velocity process}
Assume that  $T>0$ and    $(u_0,u_1)$   is an $H_A  \times H$-valued $\mathcal{F}_0$-measurable  random variable
	satisfying
\begin{equation}\label{eqn-u_0+u_1}
	\mathbb{E}\bigl( (\Vert u_0\Vert_{H_A}^2+\Vert u_1\Vert_{H}^2)^{\frac{p}{2}}\bigr)= \mathbb{E}  \bigl(\Vert (u_0,u_1)\Vert^p_{H_A \times H} \bigr) <\infty.
\end{equation}
If a process $u\in \mathbb{M}^{p}(\Za_T)$ mild solution to problem \eqref{eqn-Stoch wave approx}-\eqref{eqn-Stoch wave approx-IC} on time interval $[0,T]$, then the $\mathbb{P}$-almost surely the trajectories of the process $u$ belong to the space $C([0,T];H_A)\cap C^1([0,T];H)$ and the velocity process  $u^\prime=u_t $  belongs to the class $\mathbb{M}^{p}(Q_T)$, see \eqref{eqn-MpQT-Norm}.
\end{lemma}

The main aim of this section is to prove that problem \eqref{eqn-Stoch wave approx}-\eqref{eqn-Stoch wave approx-IC} has a unique global solution, i.e. the following result.
\begin{theorem}\label{thm-global solution to auxiliary problem}
	Let $n\in \mathbb{N}\setminus\{0\}$  and $M_n \in (0,1)$ be the corresponding number.  Let us assume that a   triple $(p,q,r)$ satisfies  condition  \eqref{InHomStrichartzIneqCondn}.   Let $H,\rV$ and $E$ be Hilbert and respectively Banach spaces defined in \eqref{eqn-spaces}. Let us assume that the maps
	$F : E \cap B_{\rV}(0,1) \to H$ and $ G: E \cap B_{\rV}(0,1) \to \gamma(K,H)$, where $K$ is a separable Hilbert space, satisfy assumptions \ref{Fassump} and \ref{Gassump} with $\gamma$, independent of $M_n$, satisfying
	\begin{equation}\label{eqn-gamma}
		0 <2 \gamma < p.
	\end{equation}
Assume that $\tilde{F}_{n}$, $\tilde{G}_{n}$ and $\theta_n $ are as above.
	Assume that  $(u_0,u_1)$   is an $H_A  \times H$-valued $\mathcal{F}_0$-measurable  random variable {satisfying  condition \eqref{eqn-u_0+u_1}}. Then the problem \eqref{eqn-Stoch wave approx}-\eqref{eqn-Stoch wave approx-IC} has a unique global solution,  denoted by  $u_n$,  in the sense of Definition \ref{def-solution approximating equation}. In particular, for every $T>0$, $u_n$ belongs to the space $\mathbb{M}^{p}(\Za_T)$.
	
	Moreover, if $0<T_1<T_2\leq \infty$ and $u_i$, $i=1,2$ is a solution to problem \eqref{eqn-Stoch wave approx}-\eqref{eqn-Stoch wave approx-IC} on respectively interval $[0,T_1]$ and $[0,T_2]$ $($or $[0,\infty)$ if $T_2=\infty)$, then $u_2$ restricted to the smaller interval $[0,T_1]$ coincides with $u_1$.

\end{theorem}

For  the remainder of this section  we choose and fix a natural number $n\in \mathbb{N}\setminus\{0\}$  and a  number $M_n \in (0,1)$.  Our proof of  Theorem \ref{thm-global solution to auxiliary problem} we will use  a modified fixed point argument.
In order to prepare for it we need to introduce the following additional notation.

Now, if $S\geq0$ and $T\geq 0$,   let us  denote by  $\mathbb{M}^{p}(\Za_{[S,S+T]})$   the Banach space of (equivalence classes) of all  $E$-valued  progressively  measurable processes $\{ u(t) ,t \in [S,S+T]\}$ having a continuous $\rV$-valued modification and  satisfying
\begin{equation}
\label{eqn-M^p(Z_S,T)}
\Vert u \Vert_{\mathbb{M}^{p}(\Za_{[S,S+T]})}^p := \tE \bigl[ \Vert u \Vert_{\Za_{[S,S+T]}}^p \bigr] = \tE \bigl[ \Vert u \Vert_{\mathrm{C}([S,S+T];\rV)}^p + \Vert u \Vert_{L^p(S,S+T;\rE)}^p \bigr] < \infty.
\end{equation}
Let us make a trivial observation that the space  $\Za_{[0,0]}$ is isomorphic with the space $\rV$ and the space $\mathbb{M}^{p}(\Za_{[0,0]})$ is isomorphic with the space  $L^p(\Omega,\mathscr{F}_0;\rV)$. In case $S=0$, we will use both the notation $\mathbb{M}^{p}(\Za_{T})$ and $\mathbb{M}^{p}(\Za_{[0,T]})$ interchangeably.  We will use similar notational convention even if we have $X,Y$ or $Q$ in place of $Z$.

\begin{definition}\label{def-Spaces depending on v}
Assume that $S \geq 0$ and $T>0$. If processes $v \in \mathbb{M}^{p}(\Za_{[0,S]})$
and $u \in \mathbb{M}^{p}(\Za_{[S,S+T]})$ agree at time  $S$, i.e. satisfy
the following compatibility condition
\begin{align}\label{eqn-compatibility}
v(S)=u(S), \;\; \mathbb{P}\mbox{-almost surely},
\end{align}
then we  define the concatenation process $\bar{u}$ by  \[\bar{u}=v \cup u, \]
and note that $\bar{u} \in \mathbb{M}^{p}(\Za_{[0,S+T]})$. The set of all processes $\bar{u}$ obtained in such a way will be denoted by   $\mathbb{M}^{p}(\Za_{[0,S+T]}^{v})$. This set  is a closed affine subspace of $\mathbb{M}^{p}(\Za_{[0,S+T]})$ and we consider the distance on  $\mathbb{M}^{p}(\Za_{[0,S+T]}^{v})$ induced by the norm on $\mathbb{M}^{p}(\Za_{[0,S+T]})$.\\
In particular, if $\mathbf{u}_0$   is an $H_A $-valued $\mathcal{F}_S$-measurable  random variable
satisfying   the following condition
\begin{equation}\label{eqn-u_0-2}
 \mathbb{E}  \bigl(\Vert \mathbf{u}_0\Vert^p_{H_A} \bigr) <\infty,
\end{equation}
then by $\mathbb{M}^{p}(\Za_{[S,S+T]}^{\mathbf{u}_0})$ we denote the closed subspace of $\mathbb{M}^{p}(\Za_{[S,S+T]})$ consisting of all processes  $u \in \mathbb{M}^{p}(\Za_{[S,S+T]})$ which satisfy
\begin{align}\label{eqn-M^p u_0}
u(S)=\mathbf{u}_0, \;\; \mathbb{P}\mbox{-almost surely}.
\end{align}
\end{definition}

\begin{definition}\label{def-PsiMap}
Assume that $S \geq 0$ and  $T>0$ and  that a  process $v\in \mathbb{M}^{p}(\Za_{[0,S]})$,  is such that $\mathbb{P}$-almost surely it's trajectories  belong to the space $C([0,S];H_A)\cap C^1([0,S];H)$ and the velocity process  $v_t $  belongs to the class $\mathbb{M}^{p}(Q_{[0,S]})$.

	A  mild solution to equation  \eqref{eqn-Stoch wave approx} on time interval $[S,S+T]$  with the \textit{history} process $v$ is a process $u $ belonging to space $\mathbb{M}^{p}(\Za_{[S,S+T]})$ such that for every $t\in [S,S+T]$, the following equality holds $\mathbb{P}$-almost surely,
	\begin{align}\nonumber
		u(t)  = \cos((t-S) \sqrt{\nA})v(S) &+ \frac{\sin((t-S) \sqrt{\nA})}{\sqrt{\nA}}v_t(S)
		+ \int_{S}^{t}  \frac{\sin((t-s) \sqrt{\nA})}{\sqrt{\nA}}  (\theta_n(\Vert \bar{u} \Vert_{\Za_{[0,s]}}) \tilde{F}_{n}(u(s))) \, ds\\
		&  + \int_{S}^{t} \frac{\sin((t-s)\sqrt{\nA})}{\sqrt{\nA}} (\theta_n(\Vert \bar{u} \Vert_{\Za_{[0,s]}})   \tilde{G}_{n}(u(s))) \, dW(s),
		\label{eqn-sol approx eqn-S}
	\end{align}
	where $\bar{u}=v \cup u$ is the concatenation process of processes $v$ and $u$ and
$\Vert \cdot \Vert_{\Za_{[0,s]}}$ is the norm in the space $\Za_{[0,s]}$.
\end{definition}

In the framework of Definition \ref{def-PsiMap} let us define a map
	\begin{equation}\label{eqn-PsiMap}
	\Psi_{[S,S+T]}^n : \mathbb{M}^{p}(\Za_{[S,S+T]}^{v(S)}) \ni u \mapsto \Psi_{[S,S+T]}^n(u) \in \mathbb{M}^{p}(\Za_{[S,S+T]}^{v(S)}),
	\end{equation}
	by  the following formula, for every $u \in \mathbb{M}^{p}(\Za_{[S,S+T]}^{v(S)})$ and all $t  \in [S,S+T]$,
	\begin{align}\label{eqn-PsiMapImage}
	[\Psi_{[S,S+T]}^n(u)](t) & = \cos\bigl((t-S) \sqrt{\nA}\bigr)v(S) + \frac{\sin\bigl((t-S) \sqrt{\nA}\bigr)}{\sqrt{\nA}}v_t(S)
	 \\
& \quad + \int_{S}^{t}  \frac{\sin((t-s) \sqrt{\nA})}{\sqrt{\nA}}  (\theta_n(\Vert \bar{u} \Vert_{\Za_{[0,s]}}) \tilde{F}_{n}(u(s))) \, ds \nonumber\\
	&  \quad + \int_{S}^{t} \frac{\sin((t-s)\sqrt{\nA})}{\sqrt{\nA}} (\theta_n(\Vert \bar{u} \Vert_{\Za_{[0,s]}})   \tilde{G}_{n}(u(s))) \, dW(s), ~\tP\textrm{-a.s..}
\nonumber
	\end{align}

We will show in the following pages that the map $\Psi_{[S,S+T]}^n$ is well defined. Let us observe that
if $T>0$,   a  process $v\in \mathbb{M}^{p}(\Za_{[0,S]})$,  is such that $\mathbb{P}$-almost surely it's trajectories  belong to the space $C([0,S];H_A)\cap C^1([0,S];H)$ and the velocity process  $v_t $  belongs to the class $\mathbb{M}^{p}(Q_{[0,S]})$ and
 $u \in \mathbb{M}^{p}(\Za_{[S,S+T]}^{v(S)})$, then  the concatenation process $\bar{u}$ is well defined and belongs to
$\mathbb{M}^{p}(\Za_{[0,S+T]}^v)$. Moreover, it follows from Proposition \ref{prop-local solution to auxilairy problem} that  if  a process  $u$ belongs to  $\mathbb{M}^{p}(\Za_{[S,S+T]}^{v(S)})$,
$\Psi_{[S,S+T]}^n(u)$  belongs to the same subspace $\mathbb{M}^{p}(\Za_{[S,S+T]}^{v(S)})$.

We continue with the following a'priori estimates and the uniqueness results for the solutions of the approximating problem \eqref{eqn-Stoch wave approx}. The proof of this result is based on some estimates proven in the following assertions.

\begin{theorem}\label{thm-apriori for approximating equation}
Let  $T>0$ and assume that the conditions in the statement of the Theorem \ref{thm-global solution to auxiliary problem} hold. Let $u\in  \mathbb{M}^{p}(\Za_T) $ be a 	mild solution to problem \eqref{eqn-Stoch wave approx}-\eqref{eqn-Stoch wave approx-IC} on the time interval $[0,T]$ in the sense of Definition \ref{def-solution approximating equation}. Then, there exists two constants $C_1 (T,p)$ and $C_2(T,n,p)$ such that
\begin{equation}\label{eqn-apriori for approximating equation}
\Vert u \Vert_{\mathbb{M}^{p}(\Za_{[0,T]})}^p = \tE \bigl[ \Vert u  \Vert_{\Za_{[0,T]}}^p \bigr]
\leq C_1 \tE \left[ \Vert u_0 \Vert_{H_A}^p  + \Vert u_1 \Vert_{H}^p  \right] +C_2(T,n,p).
  \end{equation}
The constants $C_1 (T,p)$ and $C_2(T,n,p)$ can be chosen in such a way that they are increasing functions w.r.t. $T$.
\end{theorem}

\begin{theorem}\label{thm-uniqueness approx eqn}
Assume that $S\geq 0$ and $T>0$.
%	Assume that  $(\mathbf{u}_0,\mathbf{u}_1)$   is an $H_A  \times H$-valued $\mathcal{F}_S$-measurable  random variable 	satisfying  condition \eqref{eqn-u_0+u_1}.
 Let $v\in \mathbb{M}^{p}(\Za_{[0,S]})$ be a process  such that $\mathbb{P}$-almost surely it's trajectories  belong to the space $C([0,S];H_A)\cap C^1([0,S];H)$ and the velocity process  $v_t $  belongs to the class $\mathbb{M}^{p}(Q_{[0,S]})$.
Then there exists at most one mild solution to equation  \eqref{eqn-Stoch wave approx} on time interval $[S,S+T]$  with the \textit{history} process $v$, in the sense of Definition \ref{def-PsiMap}.
%	to Problem  \eqref{eqn-Stoch wave approx}-\eqref{eqn-Stoch wave approx-IC-2} on time interval $[S,S+T]$.
\end{theorem}

\begin{proof}[\textbf{Proof of Theorem \ref{thm-apriori for approximating equation}}]
Recall that $n \in \mathbb{N}$ is fixed. Let us choose and fix $T>0$  and $u$   a  mild	solution to problem \eqref{eqn-Stoch wave approx}-\eqref{eqn-Stoch wave approx-IC} on the time interval $[0,T]$. Then, by definition, $u$
satisfies \eqref{eqn-sol approx eqn} and therefore,
 \begin{align}\label{eqn-norm-PsiMapImage}
	\Vert u \Vert_{\mathbb{M}^{p}(\Za_T)}^p & \leq 3^p \Vert t \mapsto \cos(t \sqrt{\nA})u_0 + \frac{\sin(t \sqrt{\nA})}{\sqrt{\nA}}u_1 \Vert_{\mathbb{M}^{p}(\Za_T)}^p  \nonumber\\
	& \quad +3^p   \Vert t \mapsto  \int_{0}^{t}  \frac{\sin((t-s) \sqrt{\nA})}{\sqrt{\nA}}  (\theta_n(\Vert u \Vert_{\Za_s}) \tilde{F}_{n}(u(s))) \, ds \Vert_{\mathbb{M}^{p}(\Za_T)}^p \nonumber\\
	& \quad + 3^p  \Vert t \mapsto  \int_{0}^{t} \frac{\sin((t-s)\sqrt{\nA})}{\sqrt{\nA}} (\theta_n(\Vert u \Vert_{\Za_s})   \tilde{G}_{n}(u(s))) \, dW(s)\Vert_{\mathbb{M}^{p}(\Za_T)}^p.
\end{align}
The first term in the right hand side is equal to $\Vert \mathscr{I}_1^n(u_0,u_1) \Vert_{\mathbb{M}^{p}(\Za_T)}^p$ as in Lemma \ref{lem_I1WellDef}, with $S=0$, and thus we have
\begin{align}
	\Vert t \mapsto \cos(t \sqrt{\nA})u_0 &+ \frac{\sin(t \sqrt{\nA})}{\sqrt{\nA}}u_1 \Vert_{\mathbb{M}^{p}(\Za_T)}^p
	\leq \bigl(C_T^p+\bar{C}_T^p)2^{p-1} \left[ \tE\Vert u_0 \Vert_{H_A}^p  + \tE\Vert u_1 \Vert_{H}^p  \right] \nonumber\\
	& =: C_1(T,p) \left[ \tE\Vert u_0 \Vert_{H_A}^p  + \tE\Vert u_1 \Vert_{H}^p  \right].
\end{align}
The second and third terms in the right hand side in \eqref{eqn-norm-PsiMapImage} can be estimated by applying Lemmata \ref{I2WellDef} and \ref{I3WellDef}, respectively, with $S=0$. In particular we obtain
\begin{align}
	& \Vert t \mapsto   \int_{0}^{t}  \frac{\sin((t-s) \sqrt{\nA})}{\sqrt{\nA}}  (\theta_n(\Vert u \Vert_{\Za_s}) \tilde{F}_{n}(u(s))) \, ds \Vert_{\mathbb{M}^{p}(\Za_T)}^p  \nonumber\\
	& \qquad  \leq (2n)^p (C_{\tilde{F}_{n}}^\prime)^p ~ (C_T^p+K_T^p)~  \bigl( T + T^{1-\frac{\gamma}{p}} (2n)^\gamma \bigr)^p, \nonumber\\
	& \Vert t \mapsto \int_{0}^{t} \frac{\sin((t-s)\sqrt{\nA})}{\sqrt{\nA}} (\theta_n(\Vert u \Vert_{\Za_s})   \tilde{G}_{n}(u(s))) \, dW(s)\Vert_{\mathbb{M}^{p}(\Za_T)}^p \nonumber\\
	& \qquad 	\leq (2n)^p ~ (C_{\tilde{G}_{n}}^\prime)^p ~ (\tilde{C}_T +\tilde{K}_T) ~ \left[  T + T^{1-\frac{2\gamma}{p}}  (2n)^{2 \gamma} \right]^{\frac{p}{2}}, \nonumber
\end{align}
where $\gamma$ satisfies \eqref{eqn-gamma} as in the statement of Theorem \ref{thm-global solution to auxiliary problem}.
By setting the sum of right hand sides of last two estimates as $C_2(T,n,p)$ we complete the proof of Theorem \ref{thm-apriori for approximating equation}.
\end{proof}

We skip the proof of Theorem \ref{thm-uniqueness approx eqn} because it can be done in a similar (in fact much simpler) way to the proof of Theorem \ref{thm-Uniqueness}.

	\begin{remark}
		It is important to note that the cut-off function $\theta_n$ plays essential role in the fixed point argument we display here because of the quasi-Lipschitz properties of $\tilde{F}_n$ and $\tilde{G}_n$, see  \eqref{ineq-tildeF} and \eqref{ineq-tildeG}, respectively.
	\end{remark}
Let us recall that the $Y_t$ norm has been defined in \eqref{eqn-Y_T-norm}.
We will show that there exists $T_n >0$ such that $\Psi_{[S,S+T_n]}^n$ is a strict contraction. We divide our argument in a couple of lemmata.
	\begin{lemma}\label{lem_I1WellDef}
		 Let $(\mathbf{u}_0,\mathbf{u}_1)$  be  an $H_A  \times H$-valued $\mathcal{F}_S$-measurable  random variable such that \eqref{eqn-u_0-2} holds true.
		If $S\geq 0$ and  $T>0$, then the map
		\begin{equation}
					\mathscr{I}_1^n : H_A \times H \ni (\mathbf{u}_0,\mathbf{u}_1) \mapsto \Big\{ [S,S+T] \ni t \mapsto  \cos( (t-S) \sqrt{\nA})\mathbf{u}_0 +
		\frac{\sin((t-S) \sqrt{\nA})}{\sqrt{\nA}}\mathbf{u}_1 \Big\} \in \mathbb{M}^{p}(\Za_{[S,S+T]}^{\mathbf{u}_0}), \nonumber
				\end{equation} is well-defined and
		\begin{align}\label{eqn-5.29}
		\Vert \mathscr{I}_1^n(\mathbf{u}_0,\mathbf{u}_1) \Vert_{\mathbb{M}^{p}(\Za_{[S,S+T]})}^p  &
		\leq \bigl(C_T^p+\bar{C}_T^p)2^{p-1} \left[ \tE\Vert \mathbf{u}_0 \Vert_{H_A}^p  + \tE\Vert \mathbf{u}_1 \Vert_{H}^p  \right].
		\end{align}
	\end{lemma}
	\begin{proof}[\textbf{Proof of Lemma \ref{lem_I1WellDef}}]
				It is obvious to see that $[\mathscr{I}_1^n(\mathbf{u}_0,\mathbf{u}_1)](S)=\mathbf{u}_0, \mathbb{P}$-a.s.. Let us choose and fix $T>0$. It is known that, see \cite{Ball_1977}, $w:= \mathscr{I}_1^n(\mathbf{u}_0,\mathbf{u}_1)$ is the unique solution, for each $\omega \in \Omega$,  of the following homogeneous wave equation, with the Dirichlet or the Neumann boundary condition,
		\begin{equation*}\left\{
		\begin{array}{ll}
		\partial_{tt}w - \Delta w = 0, \;\; t\geq S  \\
		w(S,\cdot) = \mathbf{u}_0(\cdot), ~~ \partial_t w(S,\cdot) = \mathbf{u}_1(\cdot).
		\end{array}
		\right.
		\end{equation*}
		Moreover, for each $\omega \in \Omega$,  see Remark \ref{rem-L^infty},  $w$ belongs to $\mathrm{C}\bigl([S,S+T]; \rV\bigr) = \Xa_{[S,S+T]}$ with
\begin{equation*}
		\Vert w \Vert_{\Xa_{[S,S+T]}} \leq \bar{C}_T \left[\Vert \mathbf{u}_0 \Vert_{H_A}  + \Vert \mathbf{u}_1 \Vert_{H}  \right],
		\end{equation*}
		and, by  Theorem \ref{thm-InHomStrichartz}, $w$ belongs to $\Ya_{[S,S+T]}$ and satisfy
\begin{equation*}
		\Vert w \Vert_{\Ya_{[S,S+T]}} \leq C_T \left[\Vert \mathbf{u}_0 \Vert_{H_A}  + \Vert \mathbf{u}_1 \Vert_{H}  \right].
		\end{equation*}
		So, for every $\omega \in \Omega$, $\mathscr{I}_1^n(\mathbf{u}_0,\mathbf{u}_1) \in \Ya_{[S,S+T]} \cap \Xa_{[S,S+T]}$ and, in view of \eqref{eqn-MpYT-Norm}, we have
		\begin{align*}
			\Vert \mathscr{I}_1^n(\mathbf{u}_0,\mathbf{u}_1) \Vert_{\mathbb{M}^{p}(\Za_{[S,S+T]})}^p  \leq \bigl(C_T^p+\bar{C}_T^p)2^{p-1} \left[ \tE\Vert \mathbf{u}_0 \Vert_{H_A}^p  + \tE\Vert \mathbf{u}_1 \Vert_{H}^p  \right].
		\end{align*}
		 Furthermore, since the process $\mathscr{I}_1^n(\mathbf{u}_0,\mathbf{u}_1)$ is $\mathbb{F}$-adapted and continuous, it is progressively measurable and, hence Lemma \ref{lem_I1WellDef} follows.
\end{proof}

	In the remaining part of this section, let us fix $v\in \mathbb{M}^{p}(\Za_{[0,S]})$ a process  such that $\mathbb{P}$-almost surely it's trajectories  belong to the space $C([0,S];H_A)\cap C^1([0,S];H)$ and the velocity process  $v_t $  belongs to the class $\mathbb{M}^{p}(Q_{[0,S]})$. Moreover, we will write $\bar{u}= v \cup u$ for $u \in \mathbb{M}^{p}(\Za_{[S,S+T]})$, see Definition \ref{def-Spaces depending on v}. Furthermore,  the assumptions and notation of the Theorem \ref{thm-global solution to auxiliary problem} hold true in the rest of the section  unless otherwise stated.
	\begin{lemma}\label{I2WellDef}
		If $S\geq 0$ and $T>0$,
		then the   map
\begin{align}\label{eqn-I_2^n}
		\mathscr{I}_2^n : \mathbb{M}^{p}(\Za_{[S,S+T]}^v)& \ni u \mapsto
		 \tilde{u}  =\Big\{ [S,S+T] \ni t \mapsto \\&  \int_{S}^{t} \frac{\sin((t-s) \sqrt{\nA})}{\sqrt{\nA}} (\theta_n(\Vert \bar{u} \Vert_{\Za_{[0,s]}})\tilde{F}_{n}(u(s))) \, ds \Big\} \in \mathbb{M}^{p}(\Za_{[S,S+T]}), \nonumber
		\end{align}
		is well-defined and satisfies
		\begin{equation}
\label{eqn--I_2^n-estimates}
		\tE \Bigl[ \Vert  \tilde{u} \Vert_{\Za_{[S,S+T]}}^p \Bigr] \leq  (2n)^p (C_{\tilde{F}_{n}}^\prime)^p ~ (C_T^p+K_T^p)~  \bigl( T + T^{1-\frac{\gamma}{p}} (2n)^\gamma \bigr)^p.
		\end{equation}
	\end{lemma}

	\begin{proof}[\textbf{Proof of Lemma \ref{I2WellDef}}] 				
	Let us choose and fix $T>0$. The computation is pathwise and thus we will not write $\omega$ explicitly unless any confusion arises. Take an arbitrary  $u \in \mathbb{M}^{p}(\Za_{[S,S+T]}^v)$ and $\tilde{u} = \mathscr{I}_2^n (u)$. Then, by definition \eqref{eqn-I_2^n} of the map $\mathscr{I}_2^n$,   we have
		\begin{align}\label{I2WellDef-t2}
		\Vert  \tilde{u} \Vert_{\Xa_{[S,S+T]}}   & \leq   \sup_{t \in [S,S+T]} \Bigl\|
\int_{S}^{t}   \frac{\sin((t-s) \sqrt{\nA}) }{\sqrt{\nA}} ( \theta_n(\Vert \bar{u} \Vert_{\Za_{[0,s]}})  \tilde{F}_{n}(u(s)) ) \, ds
\Bigr\|_{\rV}
\\
		& \leq K_T  \int_{S}^{S+T} \theta_n(\Vert \bar{u} \Vert_{\Za_{[0,s]}})  \|\tilde{F}_{n}(u(s))\|_{H} \, ds  .
\nonumber		\end{align}
		Note that, in the last step we used the following bound, which is a consequence of \cite[Lemma~2.2]{Pazy1983B} applied to $C_0$-group associated to the wave operator $\partial_{tt}-\Delta$,
		\begin{equation}\label{groupBound}
			\left\|  \frac{\sin((t-S) \sqrt{\nA})}{\sqrt{\nA}} \right\|_{\mathcal{L}(H,H_A)} \leq K_T, ~~ t \in [S,S+T],
		\end{equation} where $K_T := \tilde{M} e^{mT}$ for some constants $m \geq 0$ and $\tilde{M} \geq 1$.
		Let $T^\ast$ be stopping time defined by
		\begin{align}\label{I2WellDef-t3}
				T^\ast  & : = \inf\{ t \in [0,T]: \Vert \bar{u} \Vert_{\Za_{[0,S+t]}} \geq 2n \}.
		\end{align}

		If the set in the definition of $T^\ast$ is empty, then we set $T^\ast = T$. 	
		Returning to  \eqref{I2WellDef-t2}, by applying \eqref{I2WellDef-t3} we get
		\begin{align}\label{I2WellDef-t5}
		\int_{S}^{S+T} & \theta_n(\Vert \bar{u} \Vert_{\Za_{[0,s]}})~ \Vert \tilde{F}_{n}(u(s))\Vert_{H} \, ds  \leq    \Vert \tilde{F}_{n}(u)\Vert_{L^1(S,S+T^\ast;H)}.
		\end{align}
		Moreover, since in view of \eqref{I2WellDef-t3},   $\Vert u \Vert_{\Xa_{[S,S+T^\ast]}}, \Vert u \Vert_{\Ya_{[S,S+T^\ast]}} \leq 2n$ and  since  $\tilde{F}_{n}(0)=0$, by Lemma \ref{lem-F}  we infer 	the following inequality
		\begin{align}\label{I2WellDef-t6}
		\Vert \tilde{F}_{n}(u)\Vert_{L^1(S,S+T^\ast;H)}  & \leq C_{\tilde{F}_{n}}'~   \bigl( T^\ast + (T^\ast)^{1-\frac{\gamma}{p}} \Vert u \Vert_{\Ya_{[S,S+T^\ast]}}^\gamma \bigr) \Vert u \Vert_{\Xa_{[S,S+T^\ast]}}   \nonumber\\
		& \leq  (2n) C_{\tilde{F}_{n}}'~  \bigl( T + T^{1-\frac{\gamma}{p}} (2n)^\gamma \bigr)  ,
		\end{align}
		where $C_{\tilde{F}_{n}}' := K(\gamma) C_{\tilde{F}_{n}}$.
		Combining \eqref{I2WellDef-t2}, \eqref{I2WellDef-t5} and \eqref{I2WellDef-t6} followed by taking the expectation give \begin{equation}\label{I2WellDef-t7}
			\tE \left[\Vert \tilde{u}\Vert_{\Xa_{[S,S+T]}}^p \right] \leq (2n)^p (C_{\tilde{F}_{n}}^\prime)^p ~K_T^p ~    \bigl( T + T^{1-\frac{\gamma}{p}}  (2n)^\gamma \bigr)^p.
		\end{equation}
By definition \eqref{eqn-I_2^n} of the map $\mathscr{I}_2^n$, invoking the inhomogeneous Strichartz estimates from Theorem \ref{thm-InHomStrichartz} followed by \eqref{I2WellDef-t6} we get		
		\begin{align}
		\Vert  \tilde{u} \Vert_{\Ya_{[S,S+T]}}  &  \leq C_T  \Vert \tilde{F}_{n}(u)\Vert_{L^1(S,S+T^\ast;H)}   \leq (2n) C_F' ~ C_T~  \bigl( T + T^{1-\frac{\gamma}{p}}(2n)^\gamma \bigr). \nonumber
		\end{align}
		Which consequently, after taking the expectation,  gives,
		\begin{align}\label{I2WellDef-t8}
		\tE \left[ \Vert  \tilde{u} \Vert_{\Ya_{[S,S+T]}}^p \right] \leq (2n)^p (C_{\tilde{F}_{n}}^\prime)^p ~ C_T^p~  \bigl( T + T^{1-\frac{\gamma}{p}} (2n)^\gamma \bigr)^p.
		\end{align}
		Hence, by estimates \eqref{I2WellDef-t7} and \eqref{I2WellDef-t8} we get the Lemma \ref{I2WellDef}.					
	\end{proof}
	
	The next result establishes the Lipschitz property of $\mathscr{I}_2^n$ as a map acting on $\mathbb{M}^{p}(\Za_{[S,S+T]})$.	
	\begin{lemma}\label{lem-I2Contraction}
		Let $S\geq 0$ be given. For every  $T>0$, then 	there exists a constant $L_2^n(T) >0$ such that the following assertions are true:
\begin{trivlist}
			\item[(i)]  for every $n \in \mathbb{N}$, the function $L_2^n(\cdot)$ is non decreasing;
			\item[(ii)]  for every $n \in \mathbb{N}$, $\lim\limits_{T \to 0}L_2^n(T) = 0$;
			\item[(iii)]  for all $u_1, u_2 \in \mathbb{M}^{p}(\Za_{[S,S+T]})$, the following inequality holds,
			\begin{equation*}
			\Vert \mathscr{I}_2^n(u_1) - \mathscr{I}_2^n(u_2) \Vert_{\mathbb{M}^{p}(\Za_{[S,S+T]})} \leq  L_2^n(T) \Vert u_1-u_2 \Vert_{\mathbb{M}^{p}(\Za_{[S,S+T]})}.
			\end{equation*}
		\end{trivlist}
	\end{lemma}
	\begin{proof}[\textbf{Proof of Lemma \ref{lem-I2Contraction}}] 						
		Let us choose and fix two arbitrary elements $u_1,u_2 \in \mathbb{M}^{p}(\Za_{[S,S+T]}^v)$. As in before, we will not write $\omega$ explicitly unless any confusion arises.  Since $\mathscr{I}_2^n$ is well-defined, we denote $\tilde{u}_1 := \mathscr{I}_2^n(u_1), \tilde{u}_2 := \mathscr{I}_2^n(u_2) \in \mathbb{M}^{p}(\Za_{[S,S+T]})$.
		As in the proof of Lemma \ref{I2WellDef}, we define the following stopping times $T_1$, $T_2$ and  $T^\ast$
		%
		%
%		\begin{align}
%			& S+T_i : = \inf\{ t \in [S,S+T]: \Vert \bar{u}_i \Vert_{\Za_{[S,t]}} \geq 2n \},  ~~ i=1,2,  \quad \textrm{ and } \quad T^\ast := \max\{ T_1, T_2 \}. \nonumber
%		\end{align} 	
		\begin{align}
			& T_i : = \inf\{ t \in [0,T]: \Vert \bar{u}_i \Vert_{\Za_{[0,S+t]}} \geq 2n \},  ~~ i=1,2,  \quad \textrm{ and } \quad T^\ast := \max\{ T_1, T_2 \}. \nonumber
		\end{align} 	
Let us observe that for $ t \in [S,S+T]$, $\mathbb{P}$-a.s.,
\begin{align*}
		& \widetilde{u_2}(t)-\widetilde{u_1}(t)=  \int_{S}^{t} \frac{\sin((t-s) \sqrt{\nA})}{\sqrt{\nA}} \bigl[ \theta_n(\Vert \bar{u}_2 \Vert_{\Za_{[0,s]}}) (\tilde{F}_{n}(u_2(s)))-
\theta_n(\Vert \bar{u}_1 \Vert_{\Za_{[0,s]}}) (\tilde{F}_{n}(u_1(s)))  \bigr]\, ds.
		\end{align*}
Moreover,  similarly to inequality \eqref{I2WellDef-t6} we get the following inequality
	\begin{align*}
  \int_{S}^{S+T_2} \Vert \tilde{F}_{n}(u_2(t)) \Vert_H \, dt  &
  	 \leq C_{\tilde{F}_{n}}'~   \bigl( T_2 + T_2^{1-\frac{\gamma}{p}} \Vert u_2 \Vert_{\Ya_{[S,S+T_2]}}^\gamma \bigr) \Vert u_2 \Vert_{\Xa_{[S,S+T_2]}}  \\
   &\leq 2n C_{\tilde{F}_{n}}'~  \bigl( T + T^{1-\frac{\gamma}{p}} (2n)^\gamma \bigr) ,
	\end{align*}

		Therefore, by invoking the inhomogeneous Strichartz estimates from Theorem \ref{thm-InHomStrichartz}, we get the following inequality
		\begin{align}\label{eqn-I2Contraction-t10}
		& \tE \left[ \Vert \tilde{u}_1 - \tilde{u}_2 \Vert_{\Ya_{[S,S+T]}}^p \right]
=		\tE\left[ \int_{S}^{S+T} \vert \tilde{u}_1(t) - \tilde{u}_2(t) \vert_{\rE}^p \, dt \right] \nonumber\\
& \leq C_T^p \tE \left[ \int_{S}^{S+T}  \vert \theta_n(\Vert \bar{u}_2 \Vert_{\Za_{[0,s]}}) (\tilde{F}_{n}(u_2(s)))-
\theta_n(\Vert \bar{u}_1 \Vert_{\Za_{[0,s]}}) (\tilde{F}_{n}(u_1(s))) \vert_{H} \, ds \right]
\nonumber\\
&= C_T^p \tE \left[ \int_{S}^{S+T^\ast}  \vert \theta_n(\Vert \bar{u}_2 \Vert_{\Za_{[0,s]}}) (\tilde{F}_{n}(u_2(s)))-
\theta_n(\Vert \bar{u}_1 \Vert_{\Za_{[0,s]}}) (\tilde{F}_{n}(u_1(s))) \vert_{H} \, ds \right]
\nonumber\\
&  \leq \bar{C}_T^p  \tE\left[ \1_{\{T_1 \leq T_2\}} \int_{S}^{S+T^\ast} \vert \theta_n(\Vert \bar{u}_2 \Vert_{\Za_{[0,s]}}) (\tilde{F}_{n}(u_2(s)))-
\theta_n(\Vert \bar{u}_1 \Vert_{\Za_{[S,s]}}) (\tilde{F}_{n}(u_1(s))) \vert_{H}  \, ds \right]^p
\nonumber \\
& \quad + \bar{C}_T^p\tE\left[ \1_{\{T_2 \leq T_1\}} \int_{S}^{S+T^\ast} \vert \theta_n(\Vert \bar{u}_2 \Vert_{\Za_{[0,s]}}) (\tilde{F}_{n}(u_2(s)))-
\theta_n(\Vert \bar{u}_1 \Vert_{\Za_{[0,s]}}) (\tilde{F}_{n}(u_1(s))) \vert_{H}  \, ds \right]^p,
\end{align}
where $\bar{C}_T^p := 2^{p-1} C_T^p$.

Next, since $ \theta_n(\Vert \bar{u}_1 \Vert_{\Za_{[0,s]}})=0$ for $s \geq S+T_1$ and $\Vert \bar{u}_2\Vert_{\Za_{[0,t]}} \leq 2n$ for $t \in [S,S+T_2]$, $\mathbb{ P}$-a.s., by using the Lemmata \ref{lem-F} and  \ref{localExist2-Lem1}, we estimate the first of the two integrals in the right hand side above as
\begin{align}\label{I2Contraction-t11}
&   \tE\left[ \1_{\{T_1 \leq T_2\}} \int_{S}^{S+T^\ast} \vert \theta_n(\Vert \bar{u}_2 \Vert_{\Za_{[0,s]}}) (\tilde{F}_{n}(u_2(s)))-
\theta_n(\Vert\bar{u}_1 \Vert_{\Za_{[0,s]}}) (\tilde{F}_{n}(u_1(s))) \vert_{H}  \, ds \right]^p
\nonumber \\
&  \quad \leq 2^{p-1}   \tE\left[ \1_{\{T_1 \leq T_2\}} \int_{S}^{T_2} \theta_n(\Vert \bar{u}_1 \Vert_{\Za_{[0,s]}}) \Vert \tilde{F}_{n}(u_2(s)) -  \tilde{F}_{n}(u_1(s)) \Vert_H \, ds \right]^p
\nonumber \\
		& \quad \quad +   2^{p-1}  \tE\left[ \1_{\{T_1 \leq T_2\}} \int_{S}^{T_2}  \vert \theta_n(\Vert \bar{u}_2 \Vert_{\Za_{[0,s]}}) - \theta_n(\Vert \bar{u}_1 \Vert_{\Za_{[0,s]}})  \vert~ \Vert \tilde{F}_{n}(u_2(s)) \Vert_H \, ds \right]^p \nonumber\\
		&  \quad =   2^{p-1}   \tE\left[ \1_{\{T_1 \leq T_2\}} \int_{S}^{S+T_1} \theta_n(\Vert \bar{u}_1\Vert_{\Za_{[0,s]}}) \Vert \tilde{F}_{n}(u_2(s)) -  \tilde{F}_{n}(u_1(s)) \Vert_H \, ds \right]^p
		\nonumber \\
		& \quad \quad +   2^{p-1}  \tE\left[ \1_{\{T_1 \leq T_2\}} \int_{S}^{T_2}  \vert \theta_n(\Vert \bar{u}_2 \Vert_{\Za_{[0,s]}}) - \theta_n(\Vert \bar{u}_1 \Vert_{\Za_{[0,s]}})  \vert~ \Vert \tilde{F}_{n}(u_2(s)) \Vert_H \, ds \right]^p \nonumber\\
		& \quad \leq 2^{p-1}   (C_{\tilde{F}_{n}}^\prime)^p  ~ \tE \left[ \1_{\{T_1 \leq T_2\}}\| u_1 - u_2 \|_{\Xa_{[S,S+T_1]}} \bigl( T_1 + T_1^{1-\frac{\gamma}{p}} \| u_1 \|_{\Ya_{[S,S+T_1]}}^\gamma + T_1^{1-\frac{\gamma}{p}}  \| u_2 \|_{\Ya_{[S,S+T_1]}}^\gamma \bigr) \right]^p  \nonumber\\
		&   \quad \quad + \frac{2^{p-1}  }{n^p} ~\tE\left[ \1_{\{T_1 \leq T_2\}} \Vert \bar{u}_1 -\bar{u}_2\Vert_{\Za_{[0,S+T_2]}} \int_{S}^{T_2} \Vert \tilde{F}_{n}(u_2(t)) \Vert_H \, dt \right]^p \nonumber\\
		& \quad = 2^{p-1}   (C_{\tilde{F}_{n}}^\prime)^p  ~ \tE \left[ \1_{\{T_1 \leq T_2\}}\| u_1 - u_2 \|_{\Xa_{[S,S+T_1]}} \bigl( T_1 + T_1^{1-\frac{\gamma}{p}} \| u_1 \|_{\Ya_{[S,S+T_1]}}^\gamma + T_1^{1-\frac{\gamma}{p}}  \| u_2 \|_{\Ya_{[S,S+T_1]}}^\gamma \bigr) \right]^p  \nonumber\\
		&   \quad \quad + \frac{2^{p-1}  }{n^p} ~\tE\left[ \1_{\{T_1 \leq T_2\}} \Vert u_1 -u_2\Vert_{\Za_{[S,S+T_2]}} \int_{S}^{T_2} \Vert \tilde{F}_{n}(u_2(t)) \Vert_H \, dt \right]^p \nonumber\\
		&  \quad \leq 2^{p-1}  (C_{\tilde{F}_{n}}^\prime)^p  ~\bigl( T_1 + 2 T_1^{1-\frac{\gamma}{p}} (2n)^\gamma  \bigr)^p  ~ \tE \left[\| u_1 - u_2 \|_{\Xa_{[S,S+T_1]}}^p \right]  \nonumber\\
		& \quad \quad + \frac{2^{p-1} }{n^p} (C_{\tilde{F}_{n}}^\prime)^p ~(2n)^p ~\bigl( T_2 +T_2^{1-\frac{\gamma}{p}} (2n)^\gamma   \bigr)^p  \tE\left[ \Vert u_1 -u_2\Vert_{\Za_{[S,S+T_2]}}^p   \right]  \nonumber\\
		&  \quad \leq  2^{2p}   (C_{\tilde{F}_{n}}^\prime)^p ~\| u_1-u_2 \|_{\mathbb{M}^{p}(\Za_{[S,S+T]})}^p ~\bigl( T + 2 T^{1-\frac{\gamma}{p}} (2n)^\gamma \bigr)^p.
		\end{align}
	Swapping between $u_2$ and $u_1$ we can analogously show the following estimate the second of the two integrals in  the RHS of \eqref{eqn-I2Contraction-t10}, i.e.
	\begin{align}\label{I2Contraction-t12}
	&\tE\left[ \1_{\{T_2 \leq T_1\}} \int_{S}^{S+T^\ast} \vert \theta_n(\Vert \bar{u}_2 \Vert_{\Za_{[0,s]}}) (\tilde{F}_{n}(u_2(s)))-
	\theta_n(\Vert \bar{u}_1 \Vert_{\Za_{[0,s]}}) (\tilde{F}_{n}(u_1(s))) \vert_{H}  \, ds \right]^p
	\nonumber \\
	& \qquad \leq 2^{2p}(C_{\tilde{F}_{n}}^\prime)^p ~  ~\| u_1-u_2 \|_{\mathbb{M}^{p}(\Za_{[S,S+T]})}^p ~\bigl( T + 2T^{1-\frac{\gamma}{p}} (2n)^\gamma \bigr)^p.
	\end{align}
Thus, by combining the computation from \eqref{eqn-I2Contraction-t10}-\eqref{I2Contraction-t12} we obtain
\begin{align*}
	& \tE \left[ \Vert  \tilde{u}_1 - \tilde{u}_2 \Vert_{\Ya_{[S,S+T]}}^p  \right]  \leq 2^{2p+1} C_T^p  ~(C_{\tilde{F}_{n}}^\prime)^p ~ \| u_1-u_2 \|_{\mathbb{M}^{p}(\Za_{[S,S+T]})}^p ~\bigl( T + 2T^{1-\frac{\gamma}{p}} (2n)^\gamma \bigr)^p.
\end{align*}
		Next, using the inequality \eqref{groupBound}, followed by repeating  calculations as in \eqref{eqn-I2Contraction-t10} and \eqref{I2Contraction-t11}, we obtain
		\begin{align*}
		& \tE \left[ \Vert  \tilde{u}_1 - \tilde{u}_2 \Vert_{\Xa_{[S,S+T]}}^p  \right]  \leq 2^{2p+1} K_T^p  ~(C_{\tilde{F}_{n}}^\prime)^p ~ \| u_1-u_2 \|_{\mathbb{M}^{p}(\Za_{[S,S+T]})}^p ~\bigl( T + 2T^{1-\frac{\gamma}{p}} (2n)^\gamma \bigr)^p.
		\end{align*}
		Hence, in combination with the estimate \eqref{eqn-I2Contraction-t10} we get
		\begin{align}
		\Vert  \tilde{u}_1 - \tilde{u}_2 \Vert_{\mathbb{M}^{p}(\Za_{[S,S+T]})}^p & \leq 2^{2p+1}  (C_T^p+K_T^p) ~ (C_{\tilde{F}_{n}}^\prime)^p  \| u_1-u_2 \|_{\mathbb{M}^{p}(\Za_{[S,S+T]})}^p \bigl( T + T^{1-\frac{\gamma}{p}} (2n)^\gamma \bigr)^p \nonumber\\
		&  =: (L_2^n(T))^p ~~ \Vert u_1 - u_2 \Vert_{\mathbb{M}^{p}(\Za_{[S,S+T]})}^p. \nonumber
		\end{align}
		Since $\gamma < p$, by definition of $L_2^n(T)$, it is clear that, for each $n \in \mathbb{N}$, $\lim\limits_{T \to 0}L_2^n(T) = 0$. Thus, the proof  of  Lemma \ref{lem-I2Contraction} is complete.		
	\end{proof}

	Next, we set $$\xi_u^n(t):= \theta_n(\Vert \bar{u} \Vert_{\Za_{[0,t]}}) \tilde{G}_{n}(u(t)), \quad t \in [S,S+T].$$ Then  we can write
  	\begin{align}\label{JXi-t1}
  	\int_{S}^{t}&  \frac{\sin((t-r)\sqrt{\nA})}{\sqrt{\nA}} (\theta_n(\Vert \bar{u} \Vert_{\Za_{[0,r]}})  \tilde{G}_{n}(u(r))) \, dW(r) \nonumber\\
  	& =: \int_{S}^{t} \frac{\sin((t-r)\sqrt{\nA})}{\sqrt{\nA}} \xi_u^n(r) \, dW(r)  = : [J\xi_u^n](t),\;\;\; t \in [S,S+T].
  	\end{align}
  	In the next result, we show that $\mathscr{I}_3^n$ maps $\mathbb{M}^{p}(\Za_{[S,S+T]}^v)$ into $\mathbb{M}^{p}(\Za_{[S,S+T]})$.
	\begin{lemma}\label{I3WellDef}
		For any $S\geq 0$ and $T>0$, the map \begin{equation}\label{eqn-I3Map}
		\mathscr{I}_3^n : \mathbb{M}^{p}(\Za_{[S,S+T]}^v) \ni u \mapsto J\xi_u^n \in \mathbb{M}^{p}(\Za_{[S,S+T]}),
		\end{equation} where $J\xi^n$ is as \eqref{JXi-t1}, is well-defined and satisfies
	\begin{align}
		\| J\xi_u^n \|_{\mathbb{M}^{p}(\Za_{[S,S+T]})}^p \leq (2n)^p ~ (C_{\tilde{G}_{n}}^\prime)
^p ~ (\tilde{C}_T +\tilde{K}_T) ~ \left[  T + T^{1-\frac{2\gamma}{p}}  (2n)^{2 \gamma} \right]^{\frac{p}{2}}.
	\end{align}
	\end{lemma}
	\begin{proof}[\textbf{Proof of Lemma \ref{I3WellDef}}]
		Take any $u \in \mathbb{M}^{p}(\Za_{[S,S+T]}^v)$ and set $\tilde{u} := J\xi_u^n $. Observe that from \eqref{StocStrichartzEstimate2}, we have
		\begin{align}\label{I3WellDef-contractionEqn1}
		\tE\left[ \Vert [J\xi_u^n] \Vert_{L^p(S,S+T;E)}^p \right] & \leq \tilde{C}_T~ \tE \left[\int_{S}^{S+T}  ~ \Vert \xi_u^n(t) \Vert_{\gamma(K,H)}^2 \, dt \right]^{\frac{p}{2}}.
		\end{align}
		Let us define \begin{equation}\nonumber
			T^\ast := \inf\{t \in[0,T] : \Vert \bar{u} \Vert_{\Za_{[0,S+t]}} \geq  2n \} ,
		\end{equation}
	 	as a stopping time.
		Since $v \in \mathbb{M}^{p}(\Za_{[0,S]}) $ and $u \in \mathbb{M}^{p}(\Za_{[S,S+T]}^v) $, the map $t \mapsto \Vert \bar{u} \Vert_{\Za_{[0,t]}}$ is non-decreasing and continuous. Consequently, we infer that $\theta_n(\Vert \bar{u} \Vert_{\Za_{[0,t]}})  = 0$ for all $t \in (S+T^\ast,S+T]$. Moreover, since $\| u \|_{\Za_{[S,t]}} \leq \| \bar{u} \|_{\Za_{[0,t]}} $ for every $t  \in [S,S+T]$, we deduce that
		\[ \Vert u \Vert_{\Za_{[S,S+T^\ast]}} \leq 2n.
\]
		Thus, invoking Lemma \ref{lem-G}, followed by the H\"older inequality give
		\begin{align}\label{I3WellDef-contractionEqn3}
		\int_{S}^{S+T}   ~ \Vert \xi_u^n(t) \Vert_{\gamma(K,H)}^2 \, dt   & \leq \int_{S}^{S+T^\ast}  \Vert   \tilde{G}_{n}(u(t)) \Vert_{\gamma(K,H)}^2 \, dt \nonumber \\
		& \leq  (C_{\tilde{G}_{n}}^\prime)^2~   \Vert u \Vert_{\Xa_{[S,S+T^\ast]}}^2 \left[ T^\ast + (T^\ast)^{1-\frac{2\gamma}{p}}  \Vert u \Vert_{\Ya_{[S,S+T^\ast]}}^{2\gamma} \right]  \nonumber\\		
		 & \leq  (2n)^2 ~(C_{\tilde{G}_{n}}^\prime)^2 ~\left[ T + T^{1-\frac{2\gamma}{p}} (2n)^{2 \gamma}\right] .
		\end{align}
		Here $C_{\tilde{G}_{n}}':= K(\gamma) C_{\tilde{G}_{n}}$.  Consequently, by using \eqref{I3WellDef-contractionEqn3} in \eqref{I3WellDef-contractionEqn1} we obtain
		\begin{equation}\label{I3WellDef-contractionEqn4}
		\tE \left[\int_{S}^{S+T} \Vert [J\xi_u^n](t) \Vert_E^p \, dt\right]  ~\leq~ (2n)^p ~ (C_{\tilde{G}_{n}}^\prime)^p ~ \tilde{C}_T ~ \left[  T + T^{1-\frac{2\gamma}{p}}  (2n)^{2 \gamma} \right]^{\frac{p}{2}}.
		\end{equation}
		Next, to estimate $\tE\left[\Vert J\xi_u^n \Vert_{\Xa_{[S,S+T]}}^p \right]$, using the stochastic Strichartz estimates from Theorem \ref{thm-StocStrichartzEstimate}, followed by \eqref{I3WellDef-contractionEqn3}, we get
		\begin{align}\label{I3WellDef-contractionEqn5}
		\tE \left[\sup_{t \in [S,S+T]}\Vert [J\xi_u^n](t) \Vert_{\rV}^p \right] &  \leq~ (2n)^p~  (C_{\tilde{G}_{n}}^\prime)^p ~ \tilde{K}_T ~  \left[ T + T^{1-\frac{2\gamma}{p}}  (2n)^{2 \gamma}\right]^{\frac{p}{2}}.
		\end{align}
		Combining \eqref{I3WellDef-contractionEqn4} and \eqref{I3WellDef-contractionEqn5} completes the proof of Lemma \ref{I3WellDef}.		
	\end{proof}

	The next result establishes the Lipschitz property of the map $\mathscr{I}_3^n$ defined in \eqref{eqn-I3Map}.
	\begin{lemma}\label{lem_I3Contraction}
Let $S\geq 0$ be given. Assume that a Hilbert space $K$  satisfies Assumption \eqref{KAssump}. If $T>0$, then there exists a constant $L_3^n(T) >0$ such that the following assertions are true:
\begin{itemize}
			\item  $L_3^n(\cdot)$ is non decreasing;
			\item  for every $n \in \mathbb{N}$, $\lim\limits_{T \to 0}L_3^n(T) = 0$;
			\item  for $u_1, u_2 \in \mathbb{M}^{p}(\Za_{[S,S+T]})$, $\mathscr{I}_3^n$ satisfy,
			\begin{equation}\label{I3ContractionIneq}
			\Vert \mathscr{I}_3^n(u_1) - \mathscr{I}_3^n(u_2) \Vert_{\mathbb{M}^{p}(\Za_{[S,S+T]})} \lesssim_{p } L_3^n(T) \Vert u_1- u_2 \Vert_{\mathbb{M}^{p}(\Za_{[S,S+T]})}.
			\end{equation}
		\end{itemize}
	\end{lemma}
	\begin{proof}[\textbf{Proof of Lemma \ref{lem_I3Contraction}}]	
		Let us choose and fix $T>0$. We set
		$$\xi_i^n(t) = \theta_n(\Vert \bar{u}_i \Vert_{\Za_{[0,t]}})  G(u_i(t)), \quad i=1,2 \textrm{ and } t \in [S,S+T].$$
		Since $\mathscr{I}_3^n$ is well-defined, we denote $J\xi_1^n   := J\xi_{u_1}^n , J\xi_2^n   := J\xi_{u_2}^n \in \mathbb{M}^{p}(\Za_{[S,S+T]})$.
		Then, note that since for $ t \in [S,S+T]$, $\mathbb{P}$-a.s.,
		\begin{align*}
			& J\xi_1^n  (t)-J\xi_2^n  (t)=  \int_{S}^{t} \frac{\sin((t-s) \sqrt{\nA})}{\sqrt{\nA}} \bigl[ \theta_n(\Vert \bar{u}_1 \Vert_{\Za_{[0,s]}}) (\tilde{G}_{n}(u_1(s)))-
			\theta_n(\Vert \bar{u}_2 \Vert_{\Za_{[0,s]}}) (\tilde{G}_{n}(u_2(s)))  \bigr]\, dW(s).
		\end{align*}
		Thus, applying \eqref{StocStrichartzEstimate2} from Theorem \ref{thm-StocStrichartzEstimate} gives
		\begin{align}\label{I3ContractionEqn1}
		\tE & \left[ \Vert J\xi_1^n  - J\xi_2^n \Vert_{L^p(0,T;E)}^p \right]  \leq \tilde{C}_T~ \tE \left[\int_{S}^{S+T}  ~ \Vert \xi_1^n(t) - \xi_2^n(t) \Vert_{\gamma(K,H)}^2 \, dt \right]^{\frac{p}{2}}.
		\end{align}
		Now, we define the following stopping times
		\begin{align}
			& T_i : = \inf\{ t \in [0,T]: \Vert \bar{u}_i \Vert_{\Za_{[0,S+t]}} \geq 2n \},  ~~ i=1,2, \quad \textrm{ and } \quad T^\ast := \max\{ T_1, T_2 \}. \nonumber
		\end{align}
		By applying the stochastic Strichartz estimate \eqref{StocStrichartzEstimate1} from Theorem \ref{thm-StocStrichartzEstimate}, we get
		\begin{align}\label{2DSWE-I3ContractionEqn2}
		\tE& \left[\sup_{t \in [S,S+T]}\Vert [J\xi_1^n](t)  - [J\xi_2^n](t)  \Vert_{\rV}^p \right]  \leq ~\tilde{K}_T~ \tE \left[\int_{S}^{S+T}  ~ \Vert \xi_1^n(t) - \xi_2^n(t) \Vert_{\gamma(K,H)}^2 \, dt \right]^{\frac{p}{2}}.
		\end{align}	
		Next,  observe that computation similar to \eqref{eqn-I2Contraction-t10} gives
		\begin{align}\label{I3ContractionEqn3a}
			\tE &   \left[\int_{S}^{S+T}  ~ \Vert \xi_1^n(t) - \xi_2^n(t) \Vert_{\gamma(K,H)}^2 \, dt \right]^{\frac{p}{2}}    \leq 2^{p-1} \tE  \left[ \1_{\{T_1 \leq T_2\}}  \int_{S}^{S+T}  ~ \Vert \xi_1^n(t) - \xi_2^n(t) \Vert_{\gamma(K,H)}^2 \, dt \right]^{\frac{p}{2}}  \nonumber\\
			& \hspace{3cm}+ 2^{p-1} \tE  \left[\1_{\{T_2 \leq T_1\}}  \int_{S}^{S+T}  ~ \Vert \xi_1^n(t) - \xi_2^n(t) \Vert_{\gamma(K,H)}^2 \, dt \right]^{\frac{p}{2}}.
		\end{align}
		As in the proof of Lemma \ref{lem-I2Contraction}, we will only estimate the first of the two integrals in the right hand side above.  Since
$ \theta_n(\Vert \bar{u}_1 \Vert_{\Za_{[0,t]}})=0$ for $t\geq S+T_1$ and $\Vert \bar{u}_2\Vert_{\Za_{[0,t]}} \leq 2n$ for $t \in [S,S+T_2]$, $\mathbb{ P}$-a.s., by following the computation of \eqref{I2Contraction-t11}  and using the Lemmata \ref{lem-G} and  \ref{localExist2-Lem1}, we estimate the first integral in the right hand side of \eqref{I3ContractionEqn3a} as
		\begin{align}\label{I3ContractionEqn3b}
			& \tE  \left[ \1_{\{T_1 \leq T_2\}}  \int_{S}^{S+T}  ~ \Vert \xi_1^n(t) - \xi_2^n(t) \Vert_{\gamma(K,H)}^2 \, dt \right]^{\frac{p}{2}} \nonumber\\
			& \quad =   \tE\left[ \1_{\{T_1 \leq T_2\}} \int_{S}^{S+T^\ast} \Vert \theta_n(\Vert \bar{u}_2 \Vert_{\Za_{[0,s]}}) (\tilde{G}_{n}(u_2(s))) -
			\theta_n(\Vert \bar{u}_1 \Vert_{\Za_{[0,s]}}) (\tilde{G}_{n}(u_1(s))) \Vert_{\gamma(K,H)}^2  \, ds \right]^{\frac{p}{2}}
			\nonumber \\
			&  \quad \leq  2^{p-1}   \tE\left[ \1_{\{T_1 \leq T_2\}} \int_{S}^{S+T_1} (\theta_n(\Vert \bar{u}_1 \Vert_{\Za_{[0,s]}}) )^2 \Vert \tilde{G}_{n}(u_2(s)) -  \tilde{G}_{n}(u_1(s)) \Vert_{\gamma(K,H)}^2 \, ds \right]^{\frac{p}{2}}
			\nonumber \\
			& \quad \quad +   2^{p-1}   \tE\left[ \1_{\{T_1 \leq T_2\}} \int_{S}^{T_2}  \vert \theta_n(\Vert \bar{u}_2 \Vert_{\Za_{[0,s]}}) - \theta_n(\Vert \bar{u}_1 \Vert_{\Za_{[0,s]}})  \vert^2~ \Vert \tilde{G}_{n}(u_2(s))\Vert_{\gamma(K,H)}^2 \, ds \right]^{\frac{p}{2}} \nonumber\\
			& \quad \leq 2^{p-1}   (C_{\tilde{G}_{n}}^\prime)^p  ~ \tE \left[ \1_{\{T_1 \leq T_2\}}\| u_1 - u_2 \|_{\Xa_{[S,S+T_1]}}^2 \bigl( T_1 + T_1^{1-\frac{2\gamma}{p}} \| u_1 \|_{\Ya_{[S,S+T_1]}}^{2\gamma} + T_1^{1-\frac{2\gamma}{p}}  \| u_2 \|_{\Ya_{[S,S+T_1]}}^{2\gamma} \bigr) \right]^{\frac{p}{2}}  \nonumber\\
			&   \quad \quad + \frac{2^{p-1}  }{n^p} ~\tE\left[ \1_{\{T_1 \leq T_2\}} \Vert \bar{u}_1 -\bar{u}_2\Vert_{\Za_{[0,S+T_2]}}^2 \int_{S}^{S+T_2} \Vert \tilde{G}_{n}(u_2(t)) \Vert_{\gamma(K,H)}^2 \, dt \right]^{\frac{p}{2}} \nonumber\\
			&  \quad \leq 2^{p-1}    (C_{\tilde{G}_{n}}^\prime)^p  ~\bigl( T_1 + 2 T_1^{1-\frac{2\gamma}{p}} (2n)^{2\gamma}  \bigr)^{\frac{p}{2}}  ~ \tE \left[\| u_1 - u_2 \|_{\Xa_{[S,S+T_1]}}^p \right]  \nonumber\\
			& \quad \quad + \frac{2^{p-1}  }{n^p} (C_{\tilde{G}_{n}}^\prime)^p ~(2n)^p~\bigl( T_2 +T_2^{1-\frac{2\gamma}{p}} (2n)^\gamma   \bigr)^{\frac{p}{2}}  \tE\left[ \Vert u_1 -u_2\Vert_{\Za_{[S,S+T_2]}}^p  \right]  \nonumber\\
			&  \quad \leq  2^{2p}   (C_{\tilde{G}_{n}}^\prime)^p ~\| u_1-u_2 \|_{\mathbb{M}^{p}(\Za_{[S,S+T]})}^p ~\bigl( T + 2 T^{1-\frac{\gamma}{p}} (2n)^\gamma \bigr)^{\frac{p}{2}}.
		\end{align}
		Swapping between $u_2$ and $u_1$ we can analogously show the following estimate the second of the two integrals in  the RHS of \eqref{I3ContractionEqn3a}, i.e.
		\begin{align}\label{I3ContractionEqn3c}
			&\tE  \left[ \1_{\{T_2 \leq T_1\}}  \int_{S}^{S+T}  ~ \Vert \xi_1^n(t) - \xi_2^n(t) \Vert_{\gamma(K,H)}^2 \, dt \right]^{\frac{p}{2}}  \nonumber\\
			& \quad = \tE\left[ \1_{\{T_2 \leq T_1\}} \int_{S}^{S+T} \Vert \theta_n(\Vert \bar{u}_2 \Vert_{\Za_{[0,s]}}) (\tilde{G}_{n}(u_2(s)))-
			\theta_n(\Vert \bar{u}_1 \Vert_{\Za_{[0,s]}}) (\tilde{G}_{n}(u_1(s))) \Vert_{\gamma(K,H)}^2 \, ds \right]^{\frac{p}{2}}
			\nonumber \\
			& \quad \leq 2^{2p}   (C_{\tilde{G}_{n}}^\prime)^p ~\| u_1-u_2 \|_{\mathbb{M}^{p}(\Za_{[S,S+T]})}^p ~\bigl( T + 2 T^{1-\frac{\gamma}{p}} (2n)^\gamma \bigr)^{\frac{p}{2}}.
		\end{align}
		Thus, substituting \eqref{I3ContractionEqn3a}-\eqref{I3ContractionEqn3c} into \eqref{2DSWE-I3ContractionEqn2} and and \eqref{2DSWE-I3ContractionEqn2} yield, respectively,
		\begin{align*}%\label{I3ContractionEqn4}
			 \tE \left[\sup_{t \in [S,S+T]}\Vert [J\xi_1^n](t)  - [J\xi_2^n](t)  \Vert_{\rV}^p \right]  &\leq 2^{3p-1} (C_{\tilde{G}_{n}}^\prime)^p   \tilde{K}_T \bigl( T + 2T^{1-\frac{2\gamma}{p}}  (2n)^{2\gamma}  \bigr)^{\frac{p}{2}}  \Vert u_1  -  u_2 \Vert_{\mathbb{M}^{p}(\Za_{[S,S+T]})}^p, \nonumber\\
		\end{align*}	
and 		%
		\begin{align*}
			 \tE \left[ \Vert J\xi_1^n  - J\xi_2^n \Vert_{L^p(0,T;E)}^p \right]  &\leq 2^{3p-1} (C_{\tilde{G}_{n}}^\prime)^p   \tilde{C}_T \bigl( T + 2T^{1-\frac{2\gamma}{p}}  (2n)^{2\gamma}  \bigr)^{\frac{p}{2}}  \Vert u_1  -  u_2 \Vert_{\mathbb{M}^{p}(\Za_{[S,S+T]})}^p.
		\end{align*}	
		Hence, we get
		\begin{align}
		\Vert   J\xi_1^n  - J\xi_2^n \Vert_{\mathbb{M}^{p}(\Za_{[S,S+T]})}^p & \leq 2^{3p-1} (C_{\tilde{G}_{n}}^\prime)^p   (\tilde{K}_T +\tilde{C}_T) \bigl( T + 2T^{1-\frac{2\gamma}{p}}  (2n)^{2\gamma}  \bigr)^{\frac{p}{2}}  \Vert u_1  -  u_2 \Vert_{\mathbb{M}^{p}(\Za_{[S,S+T]})}^p \nonumber\\
		\qquad &  =: (L_3^n(T))^p ~~ \Vert u_1 - u_2 \Vert_{\mathbb{M}^{p}(\Za_{[S,S+T]})}^p.\nonumber
		\end{align} 	
		Since $ 2 \gamma < p$, by definition of $L_3^n(T)$, it is clear that $\lim\limits_{T \to 0}L_3^n(T) = 0$ for every $n$. Thus, we have proved  \eqref{I3ContractionIneq} and hence the proof of  Lemma \ref{lem_I3Contraction} is complete.
	\end{proof}

 Let us recall that $n \in \mathbb{N}\setminus\{0\}$ and the positive number $M_n=\frac{n}{n+1}$ are fixed.  Next we prove the following auxiliary result which gives the existence of $T_{n} >0$ such that the map $\Psi_{[S,S+T_n]}^n$  is  contractive.
\begin{proposition}\label{prop-local solution to auxilairy problem}
There exists a positive number $T_n>0$ such that the following assertion holds.
If $S\geq 0$ and   a process $v\in \mathbb{M}^{p}(\Za_{[0,S]})$  is such that $\mathbb{P}$-almost surely it's trajectories  belong to the space $C([0,S],H_A)\cap C^1([0,S],H)$ and the velocity process  $v_t $  belongs to the class $\mathbb{M}^{p}(Q_{[0,S]})$, then
 there exists a unique $u\in \mathbb{M}^{p}(\Za_{[S,S+T_n]}^{v(S)})$
 which is   a unique mild solution   to equation  \eqref{eqn-Stoch wave approx} on time interval $[S,S+T_n]$  with the \textit{history} process $v$, in the sense of Definition \ref{def-PsiMap}.
  \end{proposition}

\begin{proof}[\textbf{Proof of Proposition \ref{prop-local solution to auxilairy problem}}]	
Let us choose and fix $S \geq 0$ and a process $v\in \mathbb{M}^{p}(\Za_{[0,S]})$  is such that $\mathbb{P}$-almost surely it's trajectories  belong to the space $C([0,S],H_A)\cap C^1([0,S],H)$ and the velocity process  $v_t $  belongs to the class $\mathbb{M}^{p}(Q_{[0,S]})$.
 We will show that there exits $T_n>0$, independent of $v$, such that the map $\Psi_{[S,S+T_n]}^n$, defined in \eqref{eqn-PsiMap}-\eqref{eqn-PsiMapImage}, is a $\frac12$-contraction in the space $\mathbb{M}^{p}(\Za_{[S,S+T_n]}^{v(S)})$.

Let us  define, for $T>0$, $L_n(T):= L_2^n(T)+L_3^n(T)$. By Lemmata \ref{lem-I2Contraction} and \ref{lem_I3Contraction} we infer that
the function $T \mapsto L_n(T)$ is nondecreasing  and  that  $\lim\limits_{T \to 0}L_n(T) = 0$.  Hence we can find $T_n>0$ such that $L_n(T_n) \leq \frac12$.

Since also by  Lemmata \ref{lem_I1WellDef}-\ref{lem_I3Contraction}, we infer that  the map $\Psi_{[S,S+T]}^n$ is well-defined on $\mathbb{M}^{p}(\Za_{[S,S+T]}^{v(S)})$ and,
for all $u_1,u_2 \in \mathbb{M}^{p}(\Za_{[S,S+T]}^{v(S)})$,
	\begin{align}
	\Vert \Psi_{[S,S+T]}^n(u_1) - \Psi_{[S,S+T]}^n(u_2) \Vert_{\mathbb{M}^{p}(\Za_{[S,S+T]})} & \leq   L_n(T) \Vert u_1 - u_2 \Vert_{\mathbb{M}^{p}(\Za_{[S,S+T]})}, \nonumber
	\end{align}
 we infer that
 $\Psi_{[S,S+T_n]}^n$ is a  $\frac12$-contraction. Thus, by the Banach Fixed Point Theorem there exists a unique fixed point $u_{n} \in \mathbb{M}^{p}(\Za_{[S,S+T_{n}]}^{v(S)})$ of the map $\Psi_{[S,S+T_n]}^n$.  Hence the proof of Proposition \ref{prop-local solution to auxilairy problem} is complete by observing that a fixed point  of the map $\Psi_{[S,S+T_n]}^n$
 is a  mild solution of problem \eqref{eqn-Stoch wave approx}  on $[S,S+T_n]$ with the \textit{history} process $v$.
 \end{proof}

Finally we are ready to embark on the proof of  Theorem \ref{thm-global solution to auxiliary problem}.
\begin{proof}[\textbf{Proof of Theorem \ref{thm-global solution to auxiliary problem}}]
Let us choose and fix $n\in \mathbb{N}\setminus\{0\}$ and  an $H_A  \times H$-valued $\mathcal{F}_0$-measurable  random variable $(u_0,u_1)$   which
satisfies  condition \eqref{eqn-u_0+u_1}.

By  Proposition \ref{prop-local solution to auxilairy problem}, there exists a positive number $T_n>0$ such that for every $S \geq 0$ and  for every
 process $v\in \mathbb{M}^{p}(\Za_{[0,S]})$, which has $\mathbb{P}$-almost surely it's trajectories  belong to the space $C([0,S],H_A)\cap C^1([0,S],H)$ and the velocity process  $v_t $  belongs to the class $\mathbb{M}^{p}(Q_{[0,S]})$, there exists a unique $u\in \mathbb{M}^{p}(\Za_{[S,S+T_n]}^v)$
 which is  a unique mild solution of problem \eqref{eqn-Stoch wave approx} on $[S,S+T_n]$. i.e., \eqref{eqn-sol approx eqn-S} satisfies.

We will construct by induction a sequence $(u^k)_{k=1}^\infty$ such that
$u^k\in \mathbb{M}^{p}(\Za_{[0,kT_n]})$  such that
$\mathbb{P}$-almost surely its trajectories  belong to the space $C^1([0,kT_n],H)$ and the velocity process  $u^k_t $  belongs to the class $\mathbb{M}^{p}(Q_{[0,kT_n]})$,  and $u^k$ is a solution, in the sense of Definition \ref{def-solution approximating equation}, of problem \eqref{eqn-Stoch wave approx} on the interval $[0,kT_n]$ with initial condition \eqref{eqn-Stoch wave approx-IC}.

  For $k=1$ we take as $u^1$ the unique process mentioned above with $S=0$ and the initial data $(\mathbf{u}_0,\mathbf{u}_1)=(u_0,u_1)$. This process satisfies the listed conditions because of
  Proposition \ref{prop-local solution to auxilairy problem} and Lemma \ref{lem-velocity process}.

  Suppose next that $k\in \mathbb{N}\setminus\{0\}$ and we  have constructed a process $u^k \in \mathbb{M}^{p}(\Za_{[0,kT_n]})$  with prescribed conditions.

  Then by  Proposition \ref{prop-local solution to auxilairy problem}  with $S=kT_n$ and the history process  $v=u^k$,     there exists a unique process $u \in \mathbb{M}^{p}(\Za_{[kT_n,(k+1)T_n]}^{u^k(kT_n)})$ which solves  problem \eqref{eqn-Stoch wave approx} on the interval $[kT_n,(k+1)T_n]$ and satisfying initial condition $(u^k(kT_n),u_t^k(kT_n))$. Then we observe that the concatenation process $u^k \cup u$ makes sense and it belongs to $\mathbb{M}^{p}(\Za_{[0,(k+1)T_n]})$.
    It is sufficient to prove that $u^{k+1} :=u^k \cup u$ is a solution of  problem \eqref{eqn-Stoch wave approx} on the interval $[0,(k+1)T_n]$ with initial condition
\eqref{eqn-Stoch wave approx-IC}.  In view of Definition \ref{def-solution approximating equation}, because $u^{k}$ solves equation \eqref{eqn-Stoch wave approx} on the interval $[0,kT_n]$ with initial condition \eqref{eqn-Stoch wave approx-IC}, it is sufficient to show that $u^{k+1}$ satisfies identity \eqref{eqn-sol approx eqn} for every $t\in (kT_n,(k+1)T_n]$, i.e.
\begin{align}\label{eqn-sol approx eqn-20}
	u^{k+1}(t) & = \cos(t \sqrt{\nA})u_0  + \frac{\sin(t \sqrt{\nA})}{\sqrt{\nA}}u_1 \\
	& \quad + \int_{0}^{t}  \frac{\sin((t-s) \sqrt{\nA})}{\sqrt{\nA}}  (\theta_n(\Vert u^{k+1} \Vert_{\Za_s}) \tilde{F}_{n}(u^{k+1}(s))) \, ds \nonumber \\
	& \quad + \int_{0}^{t} \frac{\sin((t-s)\sqrt{\nA})}{\sqrt{\nA}} (\theta_n(\Vert u^{k+1} \Vert_{\Za_s})   \tilde{G}_{n}(u^{k+1}(s))) \, dW(s) ~\tP\textrm{-a.s..} \nonumber
\end{align}

Recall that $\Vert \cdot \Vert_{\Za_s}=\Vert \cdot \Vert_{\Za_{[0,s]}}$.  Let us choose and fix $t\in (kT_n,(k+1)T_n]$.
Then, on one hand by Definition \ref{def-PsiMap} of a solution on time interval $[kT_n,(k+1)T_n]$ with \textit{history} process $u^k$ and,
since by definition of the concatenation $u^{k+1}:=u^k \cup u$, the process $u=u^{k+1}$ on the interval  $ [kT_n,(k+1)T_n]$,  we have  that
	\begin{align}\label{eqn-100}
		u^{k+1}(t)
& = \cos\bigl((t-kT_n) \sqrt{\nA}\bigr)u^k(kT_n) + \frac{\sin\bigl((t-kT_n) \sqrt{\nA}\bigr)}{\sqrt{\nA}}u^k_t(kT_n)
	 \\
&\quad + \int_{kT_n}^{t}  \frac{\sin((t-s) \sqrt{\nA})}{\sqrt{\nA}}  (\theta_n(\Vert u^{k+1} \Vert_{\Za_{[0,s]}}) \tilde{F}_{n}(u^{k+1}(s))) \, ds \nonumber\\
	&  \quad+ \int_{kT_n}^{t} \frac{\sin((t-s)\sqrt{\nA})}{\sqrt{\nA}} (\theta_n(\Vert u^{k+1} \Vert_{\Za_{[0,s]}})   \tilde{G}_{n}(u^{k+1}(s))) \, dW(s), ~\tP\textrm{-a.s..}
\nonumber
	\end{align}
On the other hand,  since $u^k$ and $u^{k+1}$ are equal on $[0,kT_n]$ we have
\begin{align}\label{eqn-101}
  u^k(kT_n)&=
  \cos(kT_n \sqrt{\nA})u_0 + \frac{\sin(kT_n \sqrt{\nA})}{\sqrt{\nA}} u_1
  \\
		& \quad+ \int_{0}^{kT_n}  \frac{\sin((kT_n-s) \sqrt{\nA})}{\sqrt{\nA}}  (\theta_n(\Vert u^{k+1} \Vert_{\Za_s}) \tilde{F}_{n}(u^{k+1}(s))) \, ds
\nonumber\\
		& \quad + \int_{0}^{kT_n} \frac{\sin((kT_n-s)\sqrt{\nA})}{\sqrt{\nA}} (\theta_n(\Vert u^{k+1} \Vert_{\Za_s})   \tilde{G}_{n}(u^{k+1}(s))) \, dW(s),~\tP\textrm{-a.s.,}		\nonumber
\end{align}
and
\begin{align}\label{eqn-102}
  u^k_t(kT_n)&=
  -\sqrt{\nA} \sin(kT_n \sqrt{\nA})u_0 + \cos(kT_n \sqrt{\nA})u_1
		\\ &\quad + \int_{0}^{kT_n}  \cos((kT_n-s) \sqrt{\nA})  (\theta_n(\Vert u^{k+1} \Vert_{\Za_s}) \tilde{F}_{n}(u^{k+1}(s))) \, ds \nonumber\\
		& \quad + \int_{0}^{kT_n} \cos((kT_n-s)\sqrt{\nA}) (\theta_n(\Vert u^{k+1} \Vert_{\Za_s})   \tilde{G}_{n}(u^{k+1}(s))) \, dW(s),~\tP\textrm{-a.s..}		
\nonumber
\end{align}
Thus, in light of \eqref{eqn-100}-\eqref{eqn-101}, in order to prove \eqref{eqn-sol approx eqn-20} it is sufficient to prove  equality \eqref{eqn-sol approx eqn-3} below, i.e.,

 	\begin{align}\nonumber
		u^{k+1}_t(t)  = -\sqrt{\nA} \sin(t \sqrt{\nA})u_0 + \cos(t \sqrt{\nA}) u_1
		&+ \int_{0}^{t}  \cos((t-s) \sqrt{\nA})  (\theta_n(\Vert u^{k+1} \Vert_{\Za_s}) \tilde{F}_{n}(u^{k+1}(s))) \, ds\\
		& \hspace{-4truecm} + \int_{0}^{t} \cos((t-s)\sqrt{\nA}) (\theta_n(\Vert u^{k+1} \Vert_{\Za_s})   \tilde{G}_{n}(u^{k+1}(s))) \, dW(s),~\tP\textrm{-a.s..}			
\label{eqn-sol approx eqn-3}
	\end{align}
 However, using notation from Appendix \ref{sec-def sol},  \eqref{eqn-sol approx eqn-3}  is equivalent to  the following equality
\begin{align}\nonumber
	\mathfrak{u}^{k+1}(t) &= e^{\mathcal{A}t} \mathfrak{u}(0) + \int_{0}^{t} e^{\mathcal{A}(t-s)} (\theta_n(\Vert u^{k+1} \Vert_{\Za_s}) \mathfrak{\tilde{F}_{n}}[\mathfrak{u}^{k+1}(s) ]\, ds
	\\ &+ \int_{0}^{t} e^{\mathcal{A}(t-s)} (\theta_n(\Vert u^{k+1} \Vert_{\Za_s}) \mathfrak{\tilde{G}_{n}}[\mathfrak{u}^{k+1}(s) ]\, dW(s),~\tP\textrm{-a.s..}		
	\label{eqn-103}
\end{align}
On the other hand, identities \eqref{eqn-101}-\eqref{eqn-102} are equivalent to the following equality
\begin{align}\nonumber
	\mathfrak{u}^{k}(kT_n) &= e^{\mathcal{A}kT_n} \mathfrak{u}(0) + \int_{0}^{kT_n} e^{\mathcal{A}(kT_n-s)} (\theta_n(\Vert u^{k+1} \Vert_{\Za_s}) \mathfrak{\tilde{F}_{n}}[\mathfrak{u}^{k+1}(s) ]\, ds
	\\ &+ \int_{0}^{kT_n} e^{\mathcal{A}(kT_n-s)} (\theta_n(\Vert u^{k+1} \Vert_{\Za_s}) \mathfrak{\tilde{G}_{n}}[\mathfrak{u}^{k+1}(s) ]\, dW(s),~\tP\textrm{-a.s..}	
	\label{eqn-104}
\end{align}
Finally, we notice that equality \eqref{eqn-100} implies
\begin{align}
\nonumber
	\mathfrak{u}^{k+1}(t) &=  e^{\mathcal{A}(t-kT_n)}\mathfrak{u}^{k}(kT_n)
+ \int_{kT_n}^t e^{\mathcal{A}(t-s)} (\theta_n(\Vert u^{k+1} \Vert_{\Za_s}) \mathfrak{\tilde{F}_{n}}[\mathfrak{u}^{k+1}(s) ]\, ds
\\ &+ \int_{kT_n}^t e^{\mathcal{A}(t-s)} (\theta_n(\Vert u^{k+1} \Vert_{\Za_s}) \mathfrak{\tilde{G}_{n}}[\mathfrak{u}^{k+1}(s) ]\, dW(s).
\label{eqn-105}
\end{align}

Now, by using the semigroup property of  the $C_0$-group  $\bigl( e^{t\mathcal{A}}\bigr)_t$     we can easily show that  \eqref{eqn-105} together with
\eqref{eqn-104} imply \eqref{eqn-103}.
Thus we proved that the process $u^{k+1}$ satisfies all the required conditions.

Hence we constructed the sequence $(u^{k})_{k=1}^\infty$. Let  us also note that
 $u^k$ and $u^{k+1}$ are equal on $[0,kT_n]$, for every $k\geq 1$. This implies that the infinite concatenation process
 \[
 u_n:=\cup_{k=1}^\infty u^{k}
 \]
 is a global mild solution to the problem  \eqref{eqn-Stoch wave approx}  satisfying the initial conditions
  \eqref{eqn-Stoch wave approx-IC}. The proof is complete.
\end{proof}

%\vfill \newpage

\begin{proof}[\textbf{Alternative proof of Theorem \ref{thm-global solution to auxiliary problem}}]
Let us choose and fix $n\in \mathbb{N}\setminus\{0\}$ and  an $H_A  \times H$-valued $\mathcal{F}_0$-measurable  random variable $(u_0,u_1)$   which
satisfies  condition \eqref{eqn-u_0+u_1}. Then by applying Proposition \ref{prop-local solution to auxilairy problem} with $S=0$ we find a time $T_\ast>0$ and process
$u_\ast \in \mathbb{M}^{p}(\Za_{[0,T_\ast]})$ which is
%the unique fixed of the map  $\Psi_{[0,T_\ast]}^n$
 the (unique) solution to our auxiliary problem \eqref{eqn-Stoch wave approx}-\eqref{eqn-Stoch wave approx-IC} on $[0,T_\ast]$. Let $\mathcal{U}$ be a family of pairs $(T,u)$, where $T>0$ and $u\in \mathbb{M}^{p}(\Za_{[0,T]})$ is a (unique) solution to our auxiliary problem \eqref{eqn-Stoch wave approx}-\eqref{eqn-Stoch wave approx-IC} on $[0,T_\ast]$. By what we have just proven this family is non-empty, because $(T_\ast,u_\ast)$ belongs to it.

Suppose first that the set $\bigl\{ T: (T,u) \in \mathcal{U}\bigr\}$ is unbounded from above. Then it is easy to show that it is equal to the open interval $(0,\infty)$ and therefore we can find an $\mathcal{U}$-valued   sequence  $\{(k,u_k)\}_{k=1}^\infty$.  By the uniqueness part of Theorem \ref{thm-global solution to auxiliary problem} we infer that
the restriction of  the process $u_{k+1}$  to the interval $[0,k]$ coincides with the process $u_k$. Therefore by a standard construction there exists a unique process $u(t)$, $t\in [0,\infty)$ such that for every $k\in \mathbb{N}$,  that the restriction of  the process $u$  to the interval $[0,k]$ coincides with the process $u_k$. Finally, because $u_k$ is a solution
to our auxiliary problem \eqref{eqn-Stoch wave approx}-\eqref{eqn-Stoch wave approx-IC} on $[0,k]$ we easily can show that the process $u$ is a global solution problem \eqref{eqn-Stoch wave approx}-\eqref{eqn-Stoch wave approx-IC}.

We will prove that the set $\bigl\{ T: (T,u) \in \mathcal{U}\bigr\}$ is unbounded from above.
Suppose by contradiction  that the set $\bigl\{ T: (T,u) \in \mathcal{U}\bigr\}$ is bounded from above. Put
\begin{equation}\label{eqn-max U}
\mathbf{T}:= \sup\bigl\{ T: (T,u) \in \mathcal{U}\bigr\}<\infty.
\end{equation}
Then  it is easy to show that it is equal to the open interval $(0,\mathbf{T})$ is contained in   $\bigl\{ T: (T,u) \in \mathcal{U}\bigr\}$ and
therefore we can find an $\mathcal{U}$-valued   sequence  $\{(\mathbf{T}-\frac{1}{k},u_k)\}_{k=1}^\infty$. Then, as in the first part above,
there exists a unique process $u(t)$, $t\in [0,\mathbf{T})$ such that for every $k\in \mathbb{N}$,  that the restriction of  the process $u$  to the interval $[0,\mathbf{T}-\frac{1}{k}]$ coincides with the process $u_k$. Moreover, because $u_k$ is a solution
to our auxiliary problem \eqref{eqn-Stoch wave approx}-\eqref{eqn-Stoch wave approx-IC} on $[0,\mathbf{T}-\frac{1}{k}]$, we infer that $u$ is a
solution to  problem \eqref{eqn-Stoch wave approx}-\eqref{eqn-Stoch wave approx-IC} on $[0,\mathbf{T})$. Moreover, by Theorem
\ref{thm-apriori for approximating equation}
\begin{equation}\label{eqn-a'priori on 0,T}
C_5:=\sup_{k } \Vert u \Vert_{\mathbb{M}^{p}(\Za_{[0,\mathbf{T}-\frac{1}{k}]})}^p <\infty.
\end{equation}
Let now take $T_n=T_n( \Vert (u_0,v_1)\Vert^p_{L^p(\Omega,H_A \times H)} ,C_5  )>0$ from Proposition \ref{prop-local solution to auxilairy problem} and $k\in \mathbb{K}$ such that
$\frac{2}{k}< T_n$. Then by Proposition \ref{prop-local solution to auxilairy problem} with $S=-\mathbf{T}+\frac{1}{k}$ we can find a unique solution $v$ to the \eqref{eqn-Stoch wave approx} on the time interval
$[\mathbf{T}-\frac{1}{k},\mathbf{T}-\frac{1}{k}+T_n] \supset [\mathbf{T}-\frac{1}{k},\mathbf{T}+\frac{1}{k}]$ such that
\begin{align*}
v(\mathbf{T}-\frac{1}{k}) &= u(\mathbf{T}-\frac{1}{k}), ~~  v_t(\mathbf{T}-\frac{1}{k}) = u_t(\mathbf{T}-\frac{1}{k}).
\end{align*}
Moreover, by the uniqueness part of Theorem \ref{thm-global solution to auxiliary problem}, the union of $u$ and $v$ is a well defined process which belongs
$\mathbb{M}^{p}(\Za_{[0,\mathbf{T}+\frac{1}{k}]})$ and is as well a (unique) solution to  problem \eqref{eqn-Stoch wave approx}-\eqref{eqn-Stoch wave approx-IC} on $[0,\mathbf{T}+\frac{1}{k}]$.
Therefore $\mathbf{T}$ does not satisfy \eqref{eqn-max U}.
This contradicts the assumption that the set $\bigl\{ T: (T,u) \in \mathcal{U}\bigr\}$ is bounded from above.
Theorem \ref{thm-global solution to auxiliary problem} is thus proved.
\end{proof}

\section{The main results}\label{sec-main}

%\subsection{Statements of the results on the existence and the uniqueness of solutions}\label{subsec-statement}
%
The main result of the present paper, i.e. the existence of an unique local maximal solution  to the problem \eqref{eqn-Stoch wave main}-\eqref{eqn-Stoch wave main-IC}, will be proved in this section. Let $H,\rV$ and $E$ be Hilbert and respectively Banach spaces defined in \eqref{eqn-spaces}.
\begin{theorem}\label{thm-local existence}
Let us assume that a  triple $(p,q,r)$ satisfies  condition  \eqref{InHomStrichartzIneqCondn}.  %
Let us assume that the maps
$F : E \cap B_{\rV}(0,1) \to H$ and $ G: E \cap B_{\rV}(0,1) \to \gamma(K,H)$, where $K$ is a separable Hilbert space, satisfy assumptions \ref{Fassump} and \ref{Gassump} with $\gamma$ independent of $M\in (0,1)$ satisfying
\begin{equation}%\label{eqn-gamma}
			0 <2 \gamma < p.
		\end{equation}
Then,  for every   $H_A  \times H$-valued $\mathcal{F}_0$-measurable  random variable $(u_0,u_1)$   satisfying condition \eqref{eqn-u_0},
there exists a unique local maximal mild solution $u= \{ u(t) : t \in [0,\tau) \}$, an  announcing sequence $\{ \tau_k \}_{k \geq 1}$ of the stopping time  $\tau$, to the problem \eqref{eqn-Stoch wave main}-\eqref{eqn-Stoch wave main-IC}, in the sense of Definition \ref{del-local solution}.\\
Moreover,
if $\tP(\{\tau < \infty\})>0$, then $\tP $-almost surely on the set $\{\tau < \infty\}$, the following explosion condition is satisfied
\begin{equation}\label{eqn-explosion-n3}
\mbox{if } \lim_{t \to \tau-} \Vert u\Vert_{Z_{t}} < \infty \mbox{  then }\limsup_{t \to \tau-}  \Vert u(t)\Vert_{H_A}    =1.
\end{equation}
\end{theorem}

\begin{proof}[\textbf{Proof of Theorem \ref{thm-local existence}}]	
	The proof is divided into 3 parts. We prove the existence part in the first and the uniqueness result in the second part. We introduce the  concept of maximality with its proof in our setting in the last part. For the purpose of the  proof, as before we consider a numerical sequence $(M_n)_{n=2}^\infty=(\frac n{n+1})_{n=1}^\infty$.

	Let us choose and fix a  triple $(p,q,r)$ satisfies  condition  \eqref{InHomStrichartzIneqCondn} and  a number $\gamma \in (0,\frac{p}{2})$.
	Let  $H,\rV$  be Hilbert spaces   and let  $E$ be a Banach space  defined in \eqref{eqn-spaces}.  Let us consider maps
	$F : E \cap B_{\rV}(0,1) \to H$ and $ G: E \cap B_{\rV}(0,1) \to \gamma(K,H)$, where $K$ is a separable Hilbert space, which satisfy assumptions \ref{Fassump} and \ref{Gassump}.
	Let us recall that by  Remark \ref{rem-extensionFandG}, for every $n\in \mathbb{N}\setminus\{0,1\}$,    by Theorem \ref{thm-brz+rana_extension}, there exist maps $\tilde{F}_{n}$ and $\tilde{G}_{n}$ (we do  write the explicit dependency on $n$ for accuracy), corresponding to $M=M_n$,  defined on the whole $H_A \cap E$ and  taking values in $H$ and $\gamma(K,H)$, respectively, such that the inequalities \eqref{ineq-tildeF} and \eqref{ineq-tildeG} hold true	and satisfying the equality \eqref{eqn-extensions} with $M=M_n$.

	\subsection{A proof of the existence of a local solution}\label{subsec:existence}

	The proof of the existence of a local mild solution to problem \eqref{eqn-Stoch wave main}-\eqref{eqn-Stoch wave main-IC}, in the sense of Definition \ref{del-local solution}, is carried out here.

To prove the existence of a local mild solution to problem \eqref{eqn-Stoch wave main}-\eqref{eqn-Stoch wave main-IC} in the sense of Definition \ref{del-local solution}.
Let us choose and fix the initial data $(u_0,u_1) \in H_A  \times H$  satisfying condition \eqref{eqn-u_0}.

Let $ \{u_{0,n}, u_{1,n}\}_{n\in \mathbb{N}}$ be a sequence of initial data such that (1) $u_{0,n}=  u_0$ for every $n$; (2) $u_{1,n} \to u_1$ as $n \to \infty$ in $H$ and,  for each $n \in \mathbb{N}$ and $\omega \in  \Omega$, $u_{1,n}(\omega)$ is equal to $u_1(\omega)$ if $\vert u_1(\omega)\vert_H\leq n$ and  is equal to $0$, otherwise.

Let a process $u_{n}$ be  the unique solution  of the approximating equation \eqref{eqn-Stoch wave approx} with initial data $(u_{0,n}, u_{1,n})$ whose existence is guaranteed by Theorem \ref{thm-global solution to auxiliary problem}. Let us recall that for every $T>0$,  $u_{n} \in \mathbb{M}^{p}(\Za_{T})$.
 Let $\tau_n$ be a  stopping time defined by the   following formula
	\begin{equation}\label{eqn-tau_n}
\inf \; \bigl\{ t \in[0,\infty) : \max\{ \Vert u_{n} \Vert_{\Xa_t},  \psi\bigl(\Vert u_{n} \Vert_{\Za_{t}}\bigr)\} \geq \frac{n}{n+1} \bigr\},
%		\tau_{n} := \inf \; \bigl\{ t \in[0,\infty) : \max\{ \phi(\Vert u_{n} \Vert_{\Xa_t}),  \Vert u_{n} \Vert_{\Za_{t}}\} \geq n \bigr\},
\end{equation}
where $\psi$ is an increasing homeomorphism  defined  by  the following formula
\begin{equation}\label{eqn-phi} 
\psi:[0,\infty)\ni x \mapsto \frac{x}{1+x} \in [0,1).
\end{equation}
In view of the definition \eqref{eqn-X_T-norm} of the $\Vert  \cdot \Vert_{\Xa_t}$-norm we deduce the following equivalent definition  of the stopping time $\tau_n$,
	\begin{equation}\label{eqn-tau_n-2}
		\tau_{n} := \inf\{ t \in[0,\infty) : \max\{ \Vert u_{n}(t) \Vert_{H_A},  \psi(\Vert {u}_{n} \Vert_{\Za_{t}})\} \geq \frac{n}{n+1} \}.
\end{equation}
Let us observe that $\tau_n=0$ can be equal to $0$, for instance when $\Vert u_{0} \Vert_{H_A} \geq  \frac{n}{n+1}$.

Moreover, since we use the convention that $\inf (J) =\infty$ if $J=\emptyset$ and the function $[0,\infty) \ni t \mapsto \max\{ \Vert u_{n} \Vert_{\Xa_t},  \psi(\Vert u_{n} \Vert_{\Za_{t}})\} \in [0,\infty)$ is continuous we deduce that
\begin{align} \label{tau_n=infty}
\tau_n=\infty &\mbox{ iff } \max\{ \Vert u_{n} \Vert_{\Xa_t},  \psi(\Vert u_{n} \Vert_{\Za_{t}})\} < n \mbox{ for every } t \in [0,\infty),
\\
\label{tau_n<infty}
\mbox{ if } \tau_n< \infty &\mbox{ then  } \max\{ \Vert u_{n} \Vert_{\Xa_t},   \psi(\Vert u_{n} \Vert_{\Za_{t}})\} <  \frac{n}{n+1} \mbox{ for every } t\in[ 0,\tau_n), \\
&\mbox{ and   } \max\{ \Vert u_{n} \Vert_{\Xa_t},  \psi(\Vert u_{n} \Vert_{\Za_{t}})\} =  \frac{n}{n+1}.
\nonumber
\end{align}

 Since our   filtered probability  space $(\Omega, \mathscr{F},\tP, \mathbb{F})$ satisfies the usual hypothesis,  the notions of accessible and predictable stopping times are equivalent, see
\cite[Theorem~I.6.6]{Metivier_1982} and therefore by  \cite[Proposition I.4.14]{Metivier_1982}, the stopping time $\tau_n$ is  accessible.
Thus, we can find an announcing sequence $\{\tau_{n_k} \}_{k \in \mathbb{N}}$ for $\tau_n$.

Arguing as in the proofs of \cite[Lemma 5.2]{Brz+Millet_2014} and  \cite[Theorem 4.10]{Brz_1997},
 we  can prove  that for every  $n \in \mathbb{N}$,
\begin{equation}
	\label{eqn-increasing stopping times}
	\tau_n \leq  \tau_{n+1} \;\;\;\; \mathbb{P}\text{-a.s.}
\end{equation}
and
\begin{equation}
	\label{eqn-u_n=u_n+1}
	u_n(t)=u_{n+1}(t) \mbox{ if } t\in [0,\tau_n) \mbox{ and  } n \in \mathbb{N},\;\; \mathbb{P}\text{-a.s.}.
\end{equation}
Moreover, one can easily show, by  using once again the continuity of functions $t \mapsto \| u_{n}\|_{\Za_t}$ and $t \mapsto \| u_{n}\|_{\Xa_t}$,  that  the inequality \eqref{eqn-increasing stopping times} is strict when $\tau_n \in (0,\infty)$, i.e.
\begin{equation}
	\label{eqn-increasing stopping times-2}
	\tau_n  < \tau_{n+1} \;\;\;\; \mathbb{P}\text{-a.s. on } \{ \tau_n \in (0,\infty) \}.
\end{equation}

By taking appropriate modifications we can assume that \eqref{eqn-increasing stopping times}  is satisfied on the whole space $\Omega$.
Hence, the following limit exists
\begin{equation}\label{eqn-tau_infty=lim}
	\tau_\infty(\omega)=\lim_{n\rightarrow \infty} \tau_n(\omega), \,\, \omega \in \Omega.
\end{equation}
Note that the above definition implies that if for some $n\in\mathbb{N}$,  $\tau_n(\omega)=\infty$, then $\tau_\infty(\omega)=\infty$. Thus we infer that if
there exists $n\in\mathbb{N}$ such that  $\sup \bigl\{ \max\{ \Vert u_{n}(t) \Vert_{H_A},  \psi(\Vert {u}_{n} \Vert_{\Za_{t}}): t\in [0,\infty) \}  \bigr\}< 1$, then $\tau_\infty(\omega)=\infty$. However,
$\tau_\infty(\omega)$ can be equal to $\infty$ even if  for every $n\in\mathbb{N}$  the previous supremum is equal to $\infty$.
Finally, let us also observe that
\begin{equation}\label{eqn-tau_infty-2}
	\tau_n  < \tau_{\infty} \;\;\;\; \mathbb{P}\text{-a.s. on } \{ \tau_n \in (0,\infty) \}.
\end{equation}
On the other hand, the first assertion we observed above is equivalent to the following. If $\tau_\infty(\omega)<\infty$ then
for every $n\in\mathbb{N}$,  $\tau_n(\omega)<\infty$.

At this juncture it is important to observe that, since $\| u_n(0)\|_{\rV} < 1$, $\mathbb{P}$-almost surely, by condition  \eqref{eqn-u_0},  we infer  that there exists $n\in \mathbb{N}$
 such that
$\| u_n(0)\|_{\rV} < M_n$, $\mathbb{P}$-almost surely. Moreover,  $\Vert u_{n} \Vert_{\Za_0}=\| u_n(0)\|_{\rV}<1<n$. Therefore
since  the functions $t \mapsto \| u_{n}\|_{\Za_t}$ and $t \mapsto \| u_{n}\|_{\Xa_t}$ are continuous, we deduce that the stopping time $\tau_{n}>0$  and hence $\tau_\infty>0$.
Hence we deduce that  $\tau_\infty$ is strictly positive $\tP$-almost surely.

Arguing as before, since our probability basis  satisfies the usual hypothesis, $\tau_\infty$ is an {accessible} stopping time, see \cite[Proposition 2.3 and Lemma 2.11]{Kar-Shr-96} with $\{\tau_{n}\wedge n \}_{n \in \mathbb{N}}$ being an  announcing sequence for $\tau_\infty$. By using the claims \eqref{eqn-u_n=u_n+1} and  \eqref{eqn-increasing stopping times-2} we define a   local process  $(u,\tau_\infty)$ in  the following way
\begin{equation}
	\label{eqn-local process u}
	u(t,\omega)=u_n(t,\omega) \text{ if } t<\tau_n(\omega),\;\; \omega \in \Omega.
\end{equation}
We will show later on that  $(u,\tau_\infty)$ is a local solution to problem  \eqref{eqn-Stoch wave main}-\eqref{eqn-Stoch wave main-IC}. Before embarking on this objective let us point out that
because of the definitions \eqref{eqn-tau_n} and \eqref{eqn-tau_infty=lim} of respectively stopping times $\tau_n$ and  $\tau_\infty$  and  the process $u$ we deduce that

\begin{equation}
	\label{eqn-tau_n=hitting time}
		\tau_{n} = \inf \; \bigl\{ t \in[0,\tau_\infty) : \max\{ \Vert u \Vert_{\Xa_t},  \psi(\Vert u \Vert_{\Za_{t}})\}  \geq \frac{n}{n+1} \bigr\}.
\end{equation}
In view of the definition \eqref{eqn-X_T-norm} of the $\Vert  \cdot \Vert_{\Xa_t}$-norm we deduce the following equivalent form of identity \eqref{eqn-tau_n=hitting time}
	\begin{equation}\label{eqn-tau_n=hitting time-2}
		\tau_{n} = \inf\{ t \in[0,\tau_\infty) : \max\{ \Vert u(t) \Vert_{H_A},  \psi(\Vert u \Vert_{\Za_{t}})\}  \geq \frac{n}{n+1} \}.
\end{equation}

The above implies the following result.
\begin{corollary}\label{cor-explosion}
If  $(u,\tau_\infty)$ is a local process defined above, then, $\mathbb{P}$-almost surely the following assertion holds true. \\
If $\tau_\infty(\omega)<\infty$ and
\begin{equation}\label{eqn-Z_t sup is finite}
\sup  \bigl\{   \Vert {u} \Vert_{\Za_{t}}: t\in [0,\tau_\infty) \bigr\} <\infty,
\end{equation}
then
\begin{equation}\label{eqn-explosion}
\limsup_{t \toup \tau_\infty(\omega)} \Vert u(t,\omega) \Vert_{H_A}=1.
\end{equation}
\end{corollary}

\begin{proof}[Proof of Corollary \ref{cor-explosion}]
Indeed, $\tau_\infty(\omega)<\infty$ and condition \eqref{eqn-Z_t sup is finite} is satisfied,
then there exists $n_0\in \mathbb{N}$, such that (i) $\Vert {u} \Vert_{\Za_{t}}< n_0$ for every  $t\in [0,\tau_\infty)$, and (ii) $0<\tau_{n_0}(\omega)<\tau_\infty(\omega)$ and the sequence $\bigl(\tau_{k}(\omega)\bigr)_{k=n_0}^\infty$ is strictly increasing and convergent to $  \tau_\infty(\omega)$. Therefore, since by \eqref{eqn-tau_n=hitting time-2},
$\max\{ \Vert u(\tau_{k}(\omega)) \Vert_{H_A}, \psi(\Vert {u} \Vert_{\Za_{\tau_{k}(\omega)}})\} =\frac{k}{k+1}$ for every $k\geq n_0$, we infer that
$\Vert u(\tau_{k}(\omega)) \Vert_{H_A}=\frac{k}{k+1}$,  $k\geq n_0$. This implies that
\begin{equation}\label{cor-explosion-2}
\lim_{k\to \infty} \Vert u(\tau_{k}(\omega) \Vert_{H_A})=1.
\end{equation}
This concludes the proof of Corollary \ref{cor-explosion}. Note however, that we proved a more precise assertion \eqref{cor-explosion-2} than the claimed assertion \eqref{cor-explosion}.
\end{proof}

Let us also notice that from the above Corollary \ref{cor-explosion} we infer that the local process  $(u,\tau_\infty)$ constructed above  satisfies that last part of our main Theorem \ref{thm-local existence}.

Now we present the promised proof of the assertion  that $(u,\tau_\infty)$ is a local solution to problem  \eqref{eqn-Stoch wave main}-\eqref{eqn-Stoch wave main-IC}.  Since the accessibility of $\tau_{\infty}$ has already been established, we need to prove properties (2) and (3) of Definition \ref{del-local solution} for $\{u(t):\; t \in [0,\tau_{\infty})\}$.
The property  $\| u(t)\|_{H_A} <1$ for $t \in [0,\tau_\infty), \mathbb{P}$-a.s. directly follows be definition of $\tau_{\infty}$ and the map $\psi$. Next, since $u_n \in \mathbb{M}^{p}(\Za_{T})$ for every $T>0$, due to \eqref{eqn-local process u} we infer that $u$ belongs to $\mathbb{M}^{p}(\Za_{t \wedge \tau_{n}})$, for all $t\geq 0$.

Since, from Theorem \ref{thm-apriori for approximating equation}, $u_n \in \mathbb{M}^{p}(\Za_{T})$ for every $T>0$,  expression \eqref{eqn-sol approx eqn} gives that the following satisfies, for every $t \geq 0$,
\begin{align}\label{localExistThmEqn5}
	u_n(t) & - \cos(t \sqrt{\nA})u_0 -\frac{\sin(t \sqrt{\nA})}{\sqrt{\nA}}u_1 - \int_{0}^{t} \frac{\sin((t-s) \sqrt{\nA})}{\sqrt{\nA}} ( \theta_n(\Vert u_n\Vert_{\Za_{s}})  \tilde{F}_{n}(u_n(s)) )\, ds \nonumber\\
	& \quad \quad = \int_{0}^{t} \frac{\sin((t-s) \sqrt{\nA})}{\sqrt{\nA}} (\theta_n(\Vert u_n \Vert_{\Za_{s}})   \tilde{G}_{n}(u_n(s)) )\, dW(s), ~~ \tP\textrm{-a.s.}.
\end{align}
Observe that, since from the definition of $\mathbb{M}^{p}(\Za_{T})$, the processes on both sides of equality \eqref{localExistThmEqn5} are $\Dom(\nA^{\frac12 })$-valued continuous, and by \cite[Theorem I.2]{Protter2004B}, any two modified stochastic processes are indistinguishable, we conclude that the equality even holds when the fixed deterministic time is replaced by the random one, in particular, \eqref{localExistThmEqn5} holds for $t \wedge \tau_{n}$. Moreover, from definition of $\tau_n$ we infer that $\| u_n \|_{X_{t \wedge \tau_{n}}} \leq \frac{n}{n+1}$ for all $n \in \mathbb{N}$. Therefore, $\Vert u_n \Vert_{\Xa_{t \wedge \tau_{n}}} \leq 1- \frac{1}{n} =M_n$ for all $n \in \mathbb{N}$. Thus,  the equality \eqref{eqn-local process u} gives that
\begin{equation*}
	F(u(s)) = \tilde{F}_{n}(u_n(s)) \textrm{ and } G(u(s)) = \tilde{G}_{n}(u_n(s)), \quad \forall s \in [0,t \wedge \tau_n) \textrm{ and } n \in \mathbb{N}.
\end{equation*}
Consequently, since by the definition of $\theta_n$ and $\tau_{n}$ we have $\theta_n(\Vert u_n\Vert_{\Za_{t \wedge \tau_n}}) = 1$  for all $n,k \in \mathbb{N}$, we deduce that for every $t\geq 0$,
\begin{align}
	\int_{0}^{t\wedge \tau_{n}} & \frac{\sin((t\wedge \tau_{n}-s) \sqrt{\nA})}{\sqrt{\nA}} ( \theta_n(\Vert u_n\Vert_{\Za_{s}})  \tilde{F}_{n}(u_n(s)) ) \, ds \nonumber \\
	& \quad = \int_{0}^{t\wedge \tau_{n}} \frac{\sin((t\wedge \tau_n-s) \sqrt{\nA})}{\sqrt{\nA}} ( F(u(s)) ) \, ds, ~~ \tP\textrm{-a.s..} \nonumber
\end{align}	
Moreover, by invoking Lemma \ref{lem-A.1}, which is a generalization of \cite[Lemma A.1]{Brz+M+S_2005}, we deduce that  for every $t\geq 0$,
\begin{align}
	I(t \wedge \tau_{n}) & =  \int_{0}^{t} \frac{\sin((t-s) \sqrt{\nA})}{\sqrt{\nA}} \bigl( \theta_n(\Vert u_n \Vert_{\Za_{s\wedge \tau_{n}}})
	\mathds{1}_{[0,\tau_{n})}(s) \tilde{G}_{n}(u_n(s)) \bigr) \, dW(s) \nonumber\\
	& \quad  = I_{\tau_{n}}(G(u))(t \wedge \tau_{n}), \quad \tP\textrm{-a.s.}, \nonumber
\end{align}
where $I_{\tau_{n}}(G(u))(t)$ is defined in \eqref{eqn-stoppedStocConv}. This proves that the equality \eqref{eqn-StocMildSolnDefn}
is satisfied by the process $u$. Hence, we have completed the existence of a local solution part of Theorem \ref{thm-local existence}.

	\subsection{A proof of the uniqueness of a local solution}\label{subsec:uniqueness}
	The proof below is based on the proof of \cite[Theorem 5.14]{Brz+H+R_2020-liquid crystals penalised}.
	\begin{theorem}\label{thm-Uniqueness} Assume that $(u_0,u_1)\in \rV \times H$ and that  condition \eqref{eqn-u_0} is satisfied.
		Assume that  $(u,\tau) $ and $(v,\sigma)$  are  two local solutions of problem \eqref{eqn-Stoch wave main}-\eqref{eqn-Stoch wave main-IC},  with the same  initial data $(u_0,u_1)$.
		Then,
		\begin{equation*}
		(u_{\lvert_{[0,\sigma \wedge \tau)\times \Omega } }, \sigma \wedge \tau)\sim (v_{\lvert_{[0,\sigma \wedge \tau)\times \Omega}},\sigma \wedge \tau).
		\end{equation*}
	\end{theorem}
	\begin{proof}[Proof of Theorem \ref{thm-Uniqueness}]
Let us choose and fix   two local solutions $(u,\tau) $ and $(v,\sigma)$   of problem \eqref{eqn-Stoch wave main}-\eqref{eqn-Stoch wave main-IC}, respectively, with the same  fixed initial data $(u_0,u_1)\in \rV \times H$  satisfying condition \eqref{eqn-u_0}.
Without loss of generality, we can assume that $\tau \vee \sigma \leq T$ for some $T>0$ which we fix for the remaining of the proof.

		Let $\{\tau_n\}_{n \in \mathbb{N}}$ and $\{\sigma_n\}_{n \in \mathbb{N}}$ be the announcing sequences of $\tau$ and $\sigma$, respectively. Since the considered probability space satisfies the usual hypothesis, by  \cite[Propositions 4.3 and 4.11, and Theorem 6.6]{Metivier_1982} the stopping time $\varrho:=\tau \wedge \sigma$ is accessible and it is easy to show that $\{\varrho_n\}_{n \in \mathbb{N}}:= \{\tau_n \wedge \sigma_n\}_{n \in \mathbb{N}}$ is an announcing sequence of $\varrho$.
		
		Let us fix $n\in \mathbb{N}$. Since $(v,\sigma)$ is a local solution to the problem \eqref{eqn-Stoch wave main}-\eqref{eqn-Stoch wave main-IC} we have that, for all $t\in [0,T]$,
		\begin{align}
		v(t \wedge \sigma_n) & = \cos((t \wedge \sigma_n)\sqrt{\nA})u_0 + \frac{\sin((t \wedge \sigma_n) \sqrt{\nA})}{\sqrt{\nA}}u_1 \nonumber\\
		& \quad + \int_{0}^{t \wedge \sigma_n} \frac{\sin((t\wedge \sigma_n -s)  \sqrt{\nA})}{\sqrt{\nA}} F(v(s)) \, ds + I_{\sigma_n}(G(v))(t \wedge \sigma_n), \tP\textrm{-a.s.},  \nonumber
		\end{align}
		where
		\begin{equation*}
		I_{\sigma_n}(G(v))(t) = \int_{0}^{t} \mathds{1}_{[0,\sigma_n)}(s) \frac{\sin((t-s) \sqrt{\nA})}{\sqrt{\nA}} G(v(s)) \, dW(s), \;\; t \in [0,T].
		\end{equation*}
		Since the above holds true for for all $t\in [0,T]$, by replacing $t$ by $ t \wedge \varrho_n$ we get, $\tP$-a.s.,
		\begin{align}
		v( t \wedge \varrho_n\wedge \sigma_n) & = \cos((t \wedge \varrho_n \wedge \sigma_n)\sqrt{\nA})u_0 + \frac{\sin((t \wedge \varrho_n \wedge \sigma_n) \sqrt{\nA})}{\sqrt{\nA}}u_1 \nonumber\\
		& \quad + \int_{0}^{t \wedge \varrho_n \wedge \sigma_n} \frac{\sin((t \wedge \varrho_n \wedge \sigma_n -s)  \sqrt{\nA})}{\sqrt{\nA}} F(v(s)) \, ds + I_{\sigma_n}(G(v))(t \wedge \varrho_n \wedge \sigma_n). \nonumber
		\end{align}
		Consequently, since $\varrho_n \le \sigma_n$, Lemma \eqref{lem-A.2} yields, for every $t \in [0,T]$, $\tP$-almost surely
		\begin{align}\label{eqn-Uniqueness-2}
		v( t \wedge \varrho_n)  & = \cos((t \wedge \varrho_n)\sqrt{\nA})u_0 + \frac{\sin((t \wedge \varrho_n) \sqrt{\nA})}{\sqrt{\nA}}u_1 \nonumber\\
		& \quad + \int_{0}^{t \wedge \varrho_n} \frac{\sin((t \wedge \varrho_n -s)  \sqrt{\nA})}{\sqrt{\nA}} F(v(s)) \, ds + I_{ \varrho_n}(G(v))(t \wedge \varrho_n).
		\end{align}

This proves that $(v,\varrho)$ is a local solution to problem \eqref{eqn-Stoch wave main}-\eqref{eqn-Stoch wave main-IC}. In a similar way, we can also prove that $(u, \varrho)$ is a local solution  to problem \eqref{eqn-Stoch wave main}-\eqref{eqn-Stoch wave main-IC}.

Let us choose and fix  an arbitrary $M \in (0,1)$  such that $\Vert u_0 \Vert_{H_A}    < M$ and fix $k \in \mathbb{N}$.  We define the following six stopping times
		\begin{align*}
		& \tilde{\tau}_k :=\inf\{t\in [0,T]: \Vert u \Vert_{\Za_t} \ge k  \} \wedge \tau , \hspace{1.4truecm} 	 \tilde{\sigma}_k :=\inf\{t\in [0,T]: \Vert v \Vert_{\Za_t} \ge k  \} \wedge \sigma,  \\
		 & \tilde{\tau}_M  :=\inf\{t\in [0,T]: \Vert u \Vert_{\Xa_t} \ge M  \} \wedge \tau , \hspace{1truecm} 	 \tilde{\sigma}_M :=\inf\{t\in [0,T]: \Vert v \Vert_{\Xa_t} \ge M  \} \wedge \sigma,  \\
&  \varrho_{k,M} :=\tilde{\tau}_k \wedge \tilde{\sigma}_k \wedge \tilde{\tau}_M \wedge \tilde{\sigma}_M , \hspace{3truecm}
\varrho_{n,k,M} : =\tau_n \wedge \sigma_n\wedge  \varrho_{k,M}.
		\end{align*}
		 Arguing as in the proof of \eqref{eqn-Uniqueness-2} we can show that for all $t \in [0,T]$, $\tP$-a.s.
		\begin{align}
		v( t \wedge \varrho_{n,k,M} ) & = \cos((t \wedge \varrho_{n,k,M} )\sqrt{\nA})u_0 + \frac{\sin((t \wedge \varrho_{n,k,M} ) \sqrt{\nA})}{\sqrt{\nA}}u_1 \nonumber\\
		& \quad + \int_{0}^{t \wedge \varrho_{n,k,M} } \frac{\sin((t \wedge \varrho_{n,k,M}  -s)  \sqrt{\nA})}{\sqrt{\nA}} F(v(s)) \, ds + I_{ \varrho_{n,k,M} }(G(v)(t \wedge \varrho_{n,k,M} ). \nonumber
		\end{align}
		Similarly, we can prove that  the above identity holds with $v$ replaced by $u$. Hence, by setting $w=u-v$ we obtain, for all $t\in [0,T]$, $\tP$-a.s.,
		\begin{align}
		w( t \wedge \varrho_{n,k,M} ) & =  \int_{0}^{t \wedge \varrho_{n,k,M} } \frac{\sin((t \wedge \varrho_{n,k,M}  -s)  \sqrt{\nA})}{\sqrt{\nA}} [ F(u(s))- F(v(s))] \, ds \nonumber
\\ & \quad + \tilde{I}_{ \varrho_{n,k,M} }(G(u)-G(v) )(t \wedge \varrho_{n,k,M} ), \nonumber
		\end{align}
		where \begin{equation}
		\tilde{I}_{ \varrho_{n,k,M} }(G(u)-G(v))(t ) :=  \int_{0}^{t} \mathds{1}_{[0,\varrho_{n,k,M}) }(s)  \frac{\sin((t  -s)  \sqrt{\nA})}{\sqrt{\nA}} [ G(u(s))- G(v(s))] \, dW(s). \nonumber
		\end{equation}
		In order to move forward, we set $$ w_0 := u(0) - v(0), \qquad  w_1 := u_t(0) - v_t(0). $$
		It is obvious that $w_0=w_1=0$. Observe that
	\begin{equation*}
		\sup_{t \in [0,T]} \Vert u(t \wedge \varrho_{n,k,M}) \Vert_{\rV} \leq M \quad \mbox{ and } \quad \sup_{t \in [0,T]} \Vert v(t\wedge \varrho_{n,k,M}) \Vert_{\rV} \leq M.
	\end{equation*}
	 Invoking the inequality \eqref{ineq-L^inftyH} and the stochastic Strichartz estimates \eqref{eqn-StocStrichartzEstimate1}, with $p=2$, and observing that
the inequalities  \eqref{ineq-L^inftyH} and \eqref{eqn-StocStrichartzEstimate1} hold for derivative also,
 followed by the Assumptions \ref{Fassump} and \ref{Gassump}  we infer that,
		\begin{align}\label{ineq-5.60}
		&  \tE \left[ \Vert w \Vert_{\mathrm{C}([0,t \wedge \varrho_{n,k,M}];H_A)}^2+ \Vert w_t \Vert_{\mathrm{C}([0,t \wedge \varrho_{n,k,M}];H)}^2 \right]  \lesssim   \tE \left[  \bar{C}_{T} \bigl( \Vert w_0 \Vert_{\rV}  + \Vert w_1 \Vert_{H}  \bigr) \right]^2  \\
		& \hspace{5truecm} + \tE \left[   \bar{C}_{T} \int_{0}^{t \wedge \varrho_{n,k,M} } \| [ F(u(r))- F(v(r))]\|_H  \, dr   \right]^2  \nonumber\\
		& \hspace{5truecm} +   \tE \left[ \tilde{K}_{T}  \int_{0}^{t \wedge \varrho_{n,k,M} } \| [ G(u(r))- G(v(r))]\|_{\gamma(K,H)}^2  \, dr   \right] \nonumber\\
		& \lesssim  \tE \left[  \bar{C}_{T} \bigl( \Vert w_0 \Vert_{\rV}  + \Vert w_1 \Vert_{H}  \bigr) \right]^2   + \tE \Bigl[ \bigl[ \bar{C}_{T}^2   (C_F)^2 + \tilde{K}_{T} (C_G)^2 \bigr] \nonumber\\
&\hspace{1.5truecm}   \times \int_{0}^{t \wedge \varrho_{n,k,M} }   \left\{ \sup_{s \in [0,r]} \| w(s)\|_{H_A}^2  +  \sup_{s \in [0,r]} \| w_t(s)\|_{H}^2 \right\} \bigl( 1 + \| u(r)\|_E +  \| v(r)\|_E  \bigr)^{2\gamma} \, dr   \Bigr].
\nonumber
		\end{align}
Now we can apply Lemma \ref{lem-Gronwall} with the following choice of processes.
\begin{align*}
X(t)&:= \Vert w \Vert_{\mathrm{C}([0,t \wedge \varrho_{n,k,M}];H_A)}^2+ \Vert w_t \Vert_{\mathrm{C}([0,t \wedge \varrho_{n,k,M}];H)}^2, \quad t \ge 0, \\
%Y(s)&=Z(s):=0, \\
R(s)&:=  1 + \| u(s)\|_E^{2\gamma} +  \| v(s)\|_E^{2\gamma}, \quad s \ge 0, \\
C_0&:=\bar{C}_{T}^2   (C_F)^2 + \tilde{K}_{T} (C_G)^2.
\end{align*}
Note that by the definition of the stopping time $\varrho_{n,k,M}$ the following inequality holds and therefore Assumption \eqref{ineq-E2} from Lemma \ref{lem-Gronwall} is satisfied.
\begin{align*}
\int_0^{\varrho_{n,k,M}} R(s)\, ds-\varrho_{n,k,M} &= \int_0^{\varrho_{n,k,M}} \bigl( \| u(r)\|_E^{2\gamma}  +\| v(r)\|_E^{2\gamma}  \bigr) \, ds \leq 2 k^{2\gamma}.
\end{align*}

		Finally, inequality \eqref{ineq-5.60} holds for all possible pairs $(\tau_a,\tau_b)$ of accessible stopping times such that $0\leq \tau_a \leq \tau_b \leq \varrho_{n,k,M}$. Thus,  since $w(0)=0$,
by Lemma \ref{lem-Gronwall} we infer that
	$$\tE \left[  \sup_{s \in [0, \varrho_{n,k,M}]}  \Vert w(s)  \Vert_{\rV}^2  \right]   = 0. $$
  This implies that  there exists $\Omega_{n,k,M} \in \mathscr{F}$ such that $\mathbb{ P}(\Omega_{n,k,M})=1$ and
  \begin{equation*}
	\sup_{s \in [0, \varrho_{n,k,M}(\omega)]}  \Vert w(s)  \Vert_{\rV}^2  =0,  \quad \mbox{ if } \omega \in \Omega_{n,k,M}.
	\end{equation*}
	Put $\Omega_{k,M}=\cup_{n} \Omega_{n,k,M}$. Note that $\mathbb{ P}(\Omega_{k,M})=1$.  Since $\tau_n \wedge \sigma_n \toup \tau \wedge \sigma$ we infer that
\begin{equation*}
	\sup_{s \in [0, \varrho_{k,M}(\omega))}  \Vert w(s)  \Vert_{\rV}^2  =0,  \quad \mbox{ if } \omega \in \Omega_{k,M}.
	\end{equation*}
Put $\Omega_0 := \cup_{M\in (0,1)\cap \mathbb{Q}} \cup_{k} \Omega_{k,M}$.  Note that $\mathbb{ P}(\Omega_0)=1$.
Since $\lim_{M\toup 1} \lim_{k\to \infty} \varrho_{k,M}=\tau \wedge \sigma$, $\tP$-almost surely, we infer that
there exists  $\Omega_1 \in \mathscr{F}$ such that $\mathbb{ P}(\Omega_1)=1$ and $\varrho_{k,M} \toup \tau(\omega) \wedge \sigma(\omega)$ for each $\omega \in \Omega_1$.
Hence we infer that
\begin{equation*}
	\sup_{s \in [0,\tau \wedge \sigma (\omega))}  \Vert w(s)  \Vert_{\rV}^2  =0,  \quad  \omega \in \Omega_0 \cap \Omega_1.
	\end{equation*}
Because $\tP(\Omega_0 \cap \Omega_1)=1$, this completes the proof of uniqueness.

\subsection{Constructed solution $(u,\tau_\infty)$ is a local maximal solution}
Let $(u,\tau_\infty)$ is the local  solution to problem \eqref{eqn-Stoch wave main}-\eqref{eqn-Stoch wave main} constructed in the earlier parts of this section.

We will prove  that $(u,\tau_\infty)$ is a local maximal solution to problem \eqref{eqn-Stoch wave main}-\eqref{eqn-Stoch wave main}.
For this purpose let us  consider another local solution $\hat{u}= \{  \hat{u}(t): t \in [0,\hat{\tau}) \}$ to problem \eqref{eqn-Stoch wave main}-\eqref{eqn-Stoch wave main-IC} such that $\tP(\tau_\infty \leq \hat{\tau})=1$.  By the uniqueness of local solutions, see Theorem \ref{thm-Uniqueness}, we  deduce that the  processes $u$ and $\hat{u}_{\lvert_{[0,\tau_\infty)\times \Omega}}$ are equivalent. Hence, it is sufficient  to show that $\tP(\tau_\infty = \hat{\tau})=1$.

 Firstly  we observe that
$\tau_\infty = \hat{\tau}$,  on the set $\{ \tau_\infty=\infty\}$, $\mathbb{P}$-almost surely. So it is enough to prove that
 $\mathbb{P}(\{\tau_\infty < \hat{\tau}\})=0$.

Let us introduce two auxiliary events $\Omega_1$ and $\Omega_2=\Omega \setminus \Omega_1$ defined by

\begin{equation}\label{eqn-Omega_1}
\Omega_1:=\Bigl\{\omega \in \Omega:   \sup_{t\in [0,\tau_\infty(\omega ))} \Vert u(\cdot,\omega) \Vert_{\Za_{t}} <\infty \Bigr\}.
\end{equation}
Suppose that $\mathbb{P}(\Omega_1 \cap \{\tau_\infty < \hat{\tau}\})>0$. Let $\bigl(\hat{\tau}_n\bigr)$ be the announcing sequence of $\hat{\tau}$. Since $\hat{\tau}_n \toup \hat{\tau}$ $\mathbb{P}$-a.s., we infer that
there exists  $n \in \mathbb{N}$ such that $\mathbb{P} \bigl( \Omega_1 \cap \{ \tau_\infty< \hat{\tau}_n< \hat{\tau}   \}\bigr)>0$.
Let us observe that if $\omega \in \Omega_1 \cap \{ \tau_\infty< \hat{\tau}_n< \hat{\tau}\}  $, then
\[
 \sup_{t \in [0,\tau_\infty(\omega))} \Vert u(t,\omega)\Vert_{H_A} =
\sup_{t \in [0,\tau_\infty(\omega))} \Vert \hat{u}(t,\omega)\Vert_{H_A}\leq \sup_{t \in [0,\hat{\tau}_n(\omega)]} \Vert \hat{u}(t,\omega)\Vert_{H_A}<1.\]
This contradicts assertion of the Corollary  \ref{cor-explosion} and therefore we infer that  $\mathbb{P}(\Omega_1 \cap \{\tau_\infty < \hat{\tau}\})=0$.

Now, we suppose that $\mathbb{P}(\Omega_2 \cap \{\tau_\infty < \hat{\tau}\})>0$.
Let $\bigl(\hat{\tau}_n\bigr)$ be the announcing sequence of $\hat{\tau}$. Since $\hat{\tau}_n \toup \hat{\tau}$ $\mathbb{P}$-a.s., we infer that
there exists  $n \in \mathbb{N}$ such that $\mathbb{P} \bigl( \Omega_2 \cap \{ \tau_\infty< \hat{\tau}_n< \hat{\tau}   \}\bigr)>0$.
Let us observe that if $\omega \in \Omega_2 \cap \{ \tau_\infty< \hat{\tau}_n< \hat{\tau}\}  $, then
\[
 \sup_{t \in [0,\tau_\infty(\omega))} \Vert u(\cdot,\omega)\Vert_{\Za_{t}} =
\sup_{t \in [0,\tau_\infty(\omega))} \Vert \hat{u}(\cdot,\omega)\Vert_{\Za_{t}}\leq \sup_{t \in [0,\hat{\tau}_n(\omega)]} \Vert \hat{u}(\cdot,\omega)\Vert_{\Za_{t}}< \infty.\]
 This is impossible in view of  the definition of the set $\Omega_2$ and therefore we infer that  $\mathbb{P}(\Omega_2 \cap \{\tau_\infty < \hat{\tau}\})=0$.

 Summing up, we proved that  $\mathbb{P}(\{\tau_\infty < \hat{\tau}\})=0$.
Thus, $(u,\tau_{\infty})$ is a  local maximal solution to problem \eqref{eqn-Stoch wave main}-\eqref{eqn-Stoch wave main}.

\end{proof}
Hence, the proof of Theorem \ref{thm-local existence} is complete.
\end{proof}

\appendix

\section{Stopped processes}\label{sec-stopped processes}
In this appendix we justify the choice of $I_{\tau}(G)$ process in the Definition \ref{del-local solution}. The proof of the next result to some extent is analogous to the proof of \cite[Lemma A.1]{Brz+M+S_2005}, where this is formulated in terms of a semigroup.
\begin{lemma} \label{lem-A.1}
	Assume that a process  $\xi$  belongs to $ \mathcal{M}^{p}(\mathcal{L}^{2}([0,T], \gamma(K,H)))$. Set \begin{equation*}
	I(t) := \int_{0}^{t} \frac{\sin((t-s) \sqrt{\nA})}{\sqrt{\nA}} \xi(s) \, dW(s),
	\end{equation*}
	and
	\begin{equation}\label{res-A.1-2}
	I_\tau(t) := \int_{0}^{t}  \frac{\sin((t-s) \sqrt{\nA})}{\sqrt{\nA}} \bigl( \mathds{1}_{[0,\tau)}(s)\xi(s) \bigr) \, dW(s).
	\end{equation}
	For any stopping time $\tau$ and for all $t \in [0,T]$, the following holds \begin{equation} \label{2DSWE-t1}
	I(t \wedge \tau) = I_\tau(t \wedge \tau), \qquad \tP\textrm{-a.s.}.
	\end{equation}
\end{lemma}
\begin{remark}
	 Let us note that, since $\tau$ is a stopping time, due to  \cite[Theorem 1.6 and Proposition 4.2]{Metivier_1982} the stochastic process $\mathds{1}_{[0,\tau)}$  is progressively measurable. In particular, the integrand  in \eqref{res-A.1-2} is progressively measurable.
\end{remark}
Note that, it follows from Lemma \ref{lem-A.1} that if $\xi = 0$ on $[0,\tau)$, then $I(t \wedge \tau)=0$ for all $t \in [0,T]$, $\tP$-a.s.. It is relevant to mention that the importance of such results goes back to \cite{Brz_1997}, \cite{Brz+G_1999}, and \cite{Carroll1999T}.
The next result is useful in the proof of Theorem \ref{thm-Uniqueness}. We ask the reader to see \cite[Corollary A.2]{Brz+H+R_2020-liquid crystals penalised} for the proof.
\begin{lemma}\label{lem-A.2}
	Let $\xi \in \mathcal{M}^{p}(\mathcal{L}^{2}([0,T], \gamma(K,H)))$ and $\tau,\sigma$ be two  stopping times such that $\sigma \leq \tau$, then
	\begin{equation}
	\label{eqn-A.6}  I_\tau(t\wedge \sigma) =I_\sigma(t\wedge \sigma) \mbox{ for all $t \in [0,T]$, $\tP$-a.s.}.
	\end{equation}
\end{lemma}

\section{On pointwise evaluation}\label{sec-pointwise evaluation}

Let us now formulate a special case of \cite[Proposition 1.2.25]{HVVW2016B}. This
result is a converse to \cite[Proposition 3.19]{DaPZ2014B} and is closely related to
   \cite[Proposition B.4]{Brz+Ondr_2013}.

\begin{proposition}\label{prop-pointwise evaluation}
Assume that $(\Omega,\mathscr{F},\mathbb{F},\tP)$, where $\mathbb{F}=(\mathscr{F}_t:t\geq 0)$ is a filtered probability space.
Assume also that $E$ is a separable Banach space, $T>0$ and $p\in [1,\infty)$. Assume finally that
\[
F: \Omega \to L^p(0,T;E) \mbox{ is Bochner integrable on $(\Omega,\mathscr{F}_T,\tP)$.}
\]
Then, there exists a $\mathscr{B}([0,T])\otimes \mathscr{F}_T$/$\mathscr{B}(E)$-measurable function
\[
f:[0,T] \times \Omega \to E,
\]
and there exists $\Omega^\prime \in \mathscr{F}_T$ such that $\tP(\Omega^\prime)=1$ and for every $\omega \in \Omega^\prime$, the following equality holds\footnote{
This equality is again imprecise. Rigorously, one should replace it by $[f(\cdot,\omega)]=F(\omega) \mbox{ in } L^p(0,T;E)$, where $[\cdot ]$ is the equivalence class w.r.t. the Lebesgue measure $\Leb$. But, see Remark \ref{remar-ito integral}, it is standard to use this imprecise formulation.}
\[
f(\cdot,\omega)=F(\omega) \mbox{ in } L^p(0,T;E).
\]
Moreover, if $g:[0,T] \times \Omega \to E$ is another $\mathscr{B}([0,T])\otimes \mathscr{F}_T$/$\mathscr{B}(E)$-measurable function such that
the above assertion  holds, then $f=g$, $\Leb \otimes \tP$-almost everywhere.

\end{proposition}

\section{About the definition of a solution}\label{sec-def sol}

Here we state an equivalence, without proof, between two natural definitions of a mild solution for SPDE \eqref{AbsWe}. We begin by recalling the framework from Section \ref{sec-localExistence}. In particular, we set
 \[H = L^2({\cO}); \;\; \rV = \Dom(\nA^{\frac12 });\;\;\; E = \Dom(\nA_q^{\frac{1-r}{2}})\, ,\]
where $(p,q,r)$ is any suitable triple which satisfy \eqref{InHomStrichartzIneqCondn}.
Let us recall that the linear (unbounded) operator $\mathcal{A}$ in the space $\mathfrak{H}:=  \rV \times H$ and the $C_0$-group on $\mathcal{H}$ generated by it
have been defined in formulae \eqref{eqn-group generator} and \eqref{eqn-group S(t)}, respectively.
Let us also recall that the space $\mathbb{M}^{p}(\Za_T)$ has been defined in \eqref{eqn-MpYT-Norm}.
We assume that the maps $F$ and $G$ satisfy \ref{Fassump} and \ref{Gassump}, respectively.

\begin{proposition}\label{prop-equivalence}
Suppose that $u_0\in  \rV$, $u_1\in H$, and $T>0$. \begin{itemize}
\item If
an $ \rV$-valued  process $u= \{ u(t): t \in [0,T]\}$ such that   $u\in \mathbb{M}^{p}(\Za_T)$, is a mild  to problem \eqref{eqn-Stoch wave main}-\eqref{eqn-Stoch wave main-IC}, i.e.
for all $t \in [0,T]$, $\tP$-a.s
 \begin{align}\label{eqn-equivalence-mildSoln}
		u(t ) & = \cos(t\sqrt{\nA})u_0 + \frac{\sin(t \sqrt{\nA})}{\sqrt{\nA}}u_1  \quad + \int_{0}^{t } \frac{\sin((t -s)  \sqrt{\nA})}{\sqrt{\nA}} F(u(s)) \, ds
	\nonumber\\
		&\quad + \int_{0}^{t}  \frac{\sin((t-s) \sqrt{\nA})}{\sqrt{\nA}} G(u(s)) \, dW(s).
		\end{align}
Then,  $\tP$-almost surely, the process $u$ is differentiable as $H$-valued process and  the $ \rV\times H$-valued $\mathbb{F}$-adapted process $\mathfrak{u}$ defined by  \[\mathfrak{u}(t)=(u(t),v(t)),\;\; t\in [0,T],\]
where $v(t)=u^\prime(t)$, $t \in [0,T]$, is  a continuous $H$-valued process and satisfies,  for all $t\in [0,T]$, $\tF$-almost surely, the following equation,
\begin{align}\label{prop-equivalence-firstOrderIntEqn1}
	\mathfrak{u}(t) = e^{\mathcal{A}t} \mathfrak{u}(0) + \int_{0}^{t} e^{\mathcal{A}(t-s)} \mathfrak{F}[\mathfrak{u}(s) ]\, ds + \int_{0}^{t} e^{\mathcal{A}(t-s)} \mathfrak{G}[\mathfrak{u}(s) ]\, dW(s),
\end{align}
where  $\mathfrak{u}(0)=(u_0,u_1) \in  \rV \times H$, and
\begin{align} \label{prop-equivalence-newNonlinearities}
\mathfrak{G}[\mathfrak{u}]=\bigl(
	0,
	G(u)
	\bigr),\;\; \mathfrak{F}[\mathfrak{u} ]=\bigl(
	0,
	F(u)
	\bigr), \;\;\;\mathfrak{u}=(u,v) \in  \rV \times H. \end{align}
\item Conversely, if an $\mathbb{F}$-adapted  continuous    $ \rV \times H$-valued process  \[\mathfrak{u}(t)=(u(t),v(t)),\;\; t\in [0,T],\]
such that $u\in \mathbb{M}^{p}(\Za_T)$,  is a solution to \eqref{prop-equivalence-firstOrderIntEqn1} with notation \eqref{prop-equivalence-newNonlinearities}, then the process $u$ is a solution to \eqref{eqn-equivalence-mildSoln}.
\end{itemize}
\end{proposition}

\begin{proposition}\label{prop-equivalence-local}
Suppose that $u_0\in  \rV$, $u_1\in H$, and $\tau$ is an accessible stopping time with an announcing sequence $\{\tau_k\}_{k\in \mathbb{N}}$.
\begin{itemize}
\item If an $ \rV$-valued  local process $u= \{ u(t): t \in [0,\tau)\}$  is a local mild  to problem \eqref{eqn-Stoch wave main}-\eqref{eqn-Stoch wave main-IC},
then,  $\tP$-almost surely, the trajectories of  $u$ are differentiable  and  the $ \rV \times H$-valued $\tF$-adapted local process $\mathfrak{u}$ defined by
\begin{equation}\label{eqn-u-mathfrak}
\mathfrak{u}(t)=(u(t),v(t)),\;\; t\in [0,\tau),
\end{equation}
where $v(t)=u^\prime(t)$, $t \in [0,\tau)$, is  continuous and,  for all $t\geq 0$,  satisfies  $\tP$-almost surely, the following equations, for every $k \in \mathbb{N}$,
\begin{align}\label{eqn-equivalence-firstOrderIntEqn} %{prop-equivalence-firstOrderIntEqn}
	\mathfrak{u}(t\wedge \tau_k) = e^{\mathcal{A}(t\wedge \tau_k)} \mathfrak{u}(0) + \int_{0}^{t \wedge \tau_k} e^{\mathcal{A}(t-s)} \mathfrak{F}[\mathfrak{u}(s) ]\, ds +
\mathfrak{I}(\tau_k,\mathfrak{G})(t \wedge \tau_k),
\end{align}
where  $\mathfrak{u}(0)=(u_0,u_1) \in  \rV \times H$ and
the process $\mathfrak{I}(\tau_k,\mathfrak{G})$ is defined by
\begin{equation}\label{eqn-vector I}
  \mathfrak{I}(\tau_k,\mathfrak{G})(t):=
\int_{0}^{t} e^{\mathcal{A}(t-s)} \mathds{1}_{[0,\tau_k)}(s) \mathfrak{G}[\mathfrak{u}(s) ]\, dW(s),\;\; t\geq 0,
\end{equation}
and $\mathfrak{F}$ and $\mathfrak{G}$ are defined in \eqref{prop-equivalence-newNonlinearities}.

\item Conversely, if an $\tF$-adapted  continuous    $ \rV \times H$-valued local process $\mathfrak{u}(t)=\bigl\{(u(t),v(t)),$ $ t\in [0,\tau)\bigr\}$, such that
	\begin{equation*}
		\mbox{
		$\mathfrak{u}$ belongs to $\mathbb{M}^{p}(\Za_{t \wedge \tau_k})$, for all $t\geq 0$ and every $k$,}		
		\end{equation*}
satisfies, for all $t\geq 0$,  $\tP$-almost surely,  for every $k \in \mathbb{N}$,
equation \eqref{eqn-equivalence-firstOrderIntEqn} with the notation \eqref{prop-equivalence-newNonlinearities} and \eqref{eqn-vector I}, then the process $u=\bigl\{u(t): t\in [0,\tau)\bigr\}$ is a solution to \eqref{eqn-equivalence-mildSoln}.
\end{itemize}
\end{proposition}
The next result corresponds to Remark \ref{rem-local solution-2}

\begin{remark}\label{rem-local solution-3} It can be shown, see Remark 2.22 in \cite{Brz+E_2000}, that if  an $\mathbb{F}$-adapted and $\Dom(\nA^{\frac12 })$-valued continuous process $u= \{  u(t): t \in [0,\tau) \}$ is a  local mild solution to problem \eqref{eqn-Stoch wave main}-\eqref{eqn-Stoch wave main-IC}, then equality \eqref{eqn-equivalence-firstOrderIntEqn} with \eqref{eqn-vector I}, hold for any stopping time of the
form $\tau_k\wedge \sigma$, where $\sigma$ is an accessible stopping time.
\end{remark}

\section{Stochastic Gronwall lemma}\label{sec-Gronwall}
The following result is a slight simplification of \cite[Lemma 5.3]{GHZ2009} which in turn is a generalization of \cite[Lemma 3.9]{Duan+Millet_2009}.

\begin{lemma} \label{lem-Gronwall}
Let us assume that $\tau$ is an accessible bounded stopping time.
Let $X$
and  $R$ be real valued non-negative local processes defined on $[0,\tau)\time \Omega$ such that for some  $\kappa>0$
\begin{align}
\label{ineq-E2}
& \int_0^\tau R(s)\, ds \leq \kappa, \qquad  \tP\mbox{-almost surely},
\end{align}
and  $\tE \int_0^\tau R(s)X(s) \, ds<\infty$.  Suppose also that there exists $C_0>0$ such that for  all pairs $(\tau_a,\tau_b)$ of accessible stopping times such that $0\leq \tau_a< \tau_b\leq \tau$,
\begin{align}\label{ineq-E3}
  \tE \left[\sup_{t \in [\tau_a,\tau_b]} X(t)
   \right]
  &\leq C_0\tE \left[  X(\tau_a)+\int_{\tau_a}^{\tau_b}  R(s) X(s) \, ds \right].
\end{align}
Then, there exists a constant $C=C(C_0,T,\kappa)$ such that
\begin{align}\label{ineq-E4}
  \tE \bigl[  X(\tau )  \bigr]
  &\leq C\tE\bigl[ X(0) \bigr].
\end{align}
\end{lemma}
\begin{proof}[\textbf{Proof of Lemma \ref{lem-Gronwall}}] The proof is almost identical to the proof of \cite[Lemma~5.3]{GHZ2009}. Since our assumptions are weaker, as we only consider pairs of accessible stopping times, we only need to observe that in view of the assumption \eqref{ineq-E2} there exist a natural number $N$ and a sequence $(\tau_i)_{i=0}^{N}$ of accessible stopping times defined as follows: $\tau_0=0$ and
\[
\tau_1=\inf \left\{ t \in [\tau_0,T]: \int_{\tau_0}^t R(s) \, ds \geq \frac{1}{2C_0}  \right\}.
\]
If $\tau_1=T$, then we put  $\tau_j=T$ for all $j\geq 2$. If $\tau_1<T$, then we proceed by induction.\\
Note that, it easily follows from the definition of $\tau_j$, $N$ indeed exists and satisfies $N \leq 2\kappa C_0$.
\end{proof}

\section{Extensions of functions}
\label{sec-brze-carr:app-A}

The following result is a slight generalization of a result that has  been proved in \cite{Brz+C_2003}.

\begin{lemma} \label{brze-carr:Lem-A1}
Let\/ $X$ and  $E$ be normed vector spaces with norms denoted respectively by\/ $\Vert \,.\,
\Vert_X$ and  $\Vert \,.\,\Vert_E$.
 For $R>0$  define\/ $\psi_R
:X \cap E  \to X \cap E$ through the formula
\begin{equation}
\label{brze-carr:Ae0}
    \psi_R(x) =
    \left\{
    \begin{array}{ll}
    x, &\quad \mbox{if\/ $\Vert_X \Vert_X \leq R$,}
    \\[1ex]
    %\frac{Rx}{\Vert x\Vert_X},
    Rx \Vert x\Vert_X^{-1},
    &\quad \mbox{if\/ $\Vert_X \Vert_X > R$.}
    \end{array}
    \right.
\end{equation}
Then for all\/ $x, x_1, x_2 \in X \cap E$
\begin{align} \label{ineq-Lipschitz-brz+rana}
 \Vert \psi_R (x_2) -
\psi_R (x_1) &\Vert_X \leq 3 \Vert x_2-x_1 \Vert_X,
\\
\Vert \psi_R(x) \Vert_X &\leq R \mbox{ for every } x \in X \cap E.
\\
\label{ineq-brz+ran-lin}
\Vert \psi_R(x) \Vert_E &\leq \Vert x \Vert_E \mbox{ for every } x \in X \cap E.
\end{align}
In particular, the range of $\psi_R$ is contained in the set $B_X(0,R) \cap E$.
\end{lemma}

The next result was also formulated and proved in  \cite{Brz+C_2003}.

\begin{corollary} \label{cor-Brze+Carr:Cor-A2} Let $X$ and $Y$
be normed vector spaces with norms denoted by\/ $\Vert \,.\,
\Vert$. Suppose that a map\/ $g:X \to Y$ is Lipschitz on the
closed ball  $\bar{B} (0, R)$, $R>0$, with Lipschitz
constant  $C$. Then, there exists a bounded map $\tilde{g}
:X \to Y$ such that $\tilde{g}=g$ on $\bar{B} (0, R)$ and
$\tilde{g}$ is Lipschitz on\/ $X$, with Lipschitz constant~$3C$.
\end{corollary}

We conclude this section with a new result which tells about the existence of a quasi-Lipschitz extension with the Lipschitz constant being $3$ times the Lipschitz constant of the original quasi-Lipschitz map. This  atypical extension type theorem is a generalization of \cite[Corollary 3]{Brz+C_2003}  and  to be best knowledge of the authours, this is a new result.

\begin{theorem}\label{thm-brz+rana_extension} Let $X$, $E$ and $Y$
be normed vector spaces with norms denoted respectively by  $\Vert \,.\, \Vert_X$,  $\Vert \,.\,\Vert_E$ and $\Vert \,.\, \Vert_Y$.
Suppose that a map
\[g:\bar{B}_X (0, R) \cap E  \to Y,\]
where $\bar{B}_X (0, R)$ is the closed ball in $X$ centered at the origin and of radius $R$,
satisfies the following ``Lipschitz" property. There exists $C>0$ and $p \geq 0$ such that
\begin{align}\label{eqn-Lipschitz-01}
   \Vert \, g(x_2)-g(x_1)\, \Vert_Y &\leq C (1+ \Vert \, x_2\, \Vert_E+ \Vert \,x_1\, \Vert_E )^p
   \Vert \, x_2-x_1\, \Vert_X, \;\; x_1,x_2 \in \bar{B} (0, R) \cap E.
\end{align}
 Then, there exists a  map
 \[\tilde{g}
:X \cap E \to Y\]  such that\/ $\tilde{g}=g$ on\/ $\bar{B}_X (0, R) \cap E$ and
\begin{align}\label{eqn-Lipschitz-02}
   \Vert \, \tilde{g}(x_2)-\tilde{g}(x_1)\, \Vert_Y &\leq 3C (1+ \Vert \, x_2\, \Vert_E+ \Vert \,x_1\, \Vert_E )^p
   \Vert \, x_2-x_1\, \Vert_X, \;\; x_1,x_2 \in X \cap E.
\end{align}

\end{theorem}

%\begin{theopargself}
\begin{proof}[Proof of Theorem \ref{thm-brz+rana_extension}]

Set
\begin{align*} \tilde{g} := g \circ \psi_R.\end{align*}

Since  the range of function $\psi_R$ is contained in the set $B_X(0,R) \cap E$,   the function
$\tilde{g} :X  \cap E\to Y$
is well defined and coincides with $g$ on the set
$\bar{B}_X (0, R) \cap E$.

Let us choose and fix   for $x_1, x_2 \in X \cap E$. Then   the generalized Lipschitz
property \eqref{eqn-Lipschitz-01} of $g$, the Lipschitz, w.r.t. the $X$-norm,  property  \eqref{ineq-Lipschitz-brz+rana} of the function $\psi_R$ and the linear growth, w.r.t. the $E$-norm,
property \eqref{ineq-brz+ran-lin} of the function $\psi_R$ imply that
\begin{align*}
    \Vert \tilde{g} (x_2) - \tilde{g} (x_1) \Vert_Y
    &= \Vert g  (\psi_R (x_2)) - g ( \psi_R (x_1)) \Vert_Y
    \\
    &\leq C (1+ \Vert \, \psi_R (x_2)\, \Vert_E+ \Vert \,\psi_R (x_1)\, \Vert_E )^p\Vert \psi_R (x_2) - \psi_R (x_1) \Vert_X
    \\
    &\leq  3C (1+ \Vert \, \psi_R (x_2)\, \Vert_E+ \Vert \,\psi_R (x_1)\, \Vert_E )^p  \Vert x_2-x_1 \Vert_X
    \\
    &\leq  3C (1+ \Vert \, x_2\, \Vert_E+ \Vert \,x_1\, \Vert_E )^p  \Vert x_2-x_1 \Vert_X.
\end{align*}
This implies the inequality \eqref{eqn-Lipschitz-02} and hence completes
the proof of Theorem $\ref{thm-brz+rana_extension}$.
\end{proof}

\textbf{Acknowledgements:} The authors would like to thank Jan van Neerven for discussion related to measurability part in the proof of auxiliary result Lemma \ref{lem-AuxRes1},
and to Mark Veraar for having useful conversations on the stochastic Fubini theorem  and to provide a suitable reference for the inequality \eqref{L2-LpConvergenceRel}.
They want to thank Martin Ondrej\'at for useful discussion, in particular  about stopping times and the proof of explosion part of Theorem \ref{thm-local existence}.
They are thankful to Gaurav Dhariwal and Anupam Gumber  for careful reading of manuscript.
The authors  are also grateful to Fabian Hornung and Tomasz Kosmala for providing comments which led to improvement of the presentation of this paper. Finally, they would like to thank an anonymous referee, whose critical comments  have contributed to better and clearer proofs.
The second author wishes to thank the York Graduate Research School, to award the Overseas scholarship (ORS), and the Department of Mathematics, University of York, to provide financial support and excellent research facilities during the period of this work. The second author also would like to acknowledge the German Science Foundation DFG to provide the financial support, through the Research Unit FOR 2402,  during the work on revision of this manuscript. Finally, the second named author would like to thanks the Isaac Newton Institue for Mathematical Sciences at Cambridge and the Simons Foundation for support during the  programme ``Fractional Differential Equations" when the final version of the manuscript was prepared.  \\

\end{document}